\ifx \labelsONmargin\undefined
  \ifx\RequirePackage\undefined
    \expandafter\ifx\csname amsart.sty\endcsname\relax
        \documentstyle[12pt,amssymb,amscd]{amsart}
    \fi

    \def\atSign{@@}

    \def\mathbb{\Bbb}
    \def\mathfrak{\frak}
    \def\mathbf{\bold}
    \ifx\boldsymbol\undefined	
      \def\boldsymbol#1{{\bold #1}}
    \fi

    \def\mathbit{\boldsymbol}

    \makeatletter
    \newenvironment{proof}{%
         \@ifnextchar[{%
                       \expandafter\let\expandafter\end@proof
                         \csname endpf*\endcsname
                         \my@proof
                      }{\let\end@proof\endpf\pf}%
        }{\end@proof}
    \def\my@proof[#1]{\@nameuse{pf*}{#1}}
    \makeatother

    \def\xrightarrow[#1]#2{@>{#2}>{#1}>}
    \def\xleftarrow[#1]#2{@<{#2}<{#1}<}

    \def\providecommand#1{\def#1}
    \def\emph#1{{\em #1}}
    \def\textbf#1{{\bf #1}}
  \else
      \ifx\usepackage\RequirePackage
      \else
        \documentclass[12pt]{amsart}
        \usepackage{amssymb,amscd}
      \fi

      \ifx\mathbit\undefined	
        \DeclareMathAlphabet{\mathbit}{OML}{cmm}{b}{it}
      \fi

      \def\atSign{@}

      \advance\oddsidemargin -0.1\textwidth
      \evensidemargin\oddsidemargin
      \def\Sb#1\endSb{_{\substack{#1}}}
      \def\Sp#1\endSp{^{\substack{#1}}}
  \fi
\makeatletter
\@ifundefined{DeclareFontShape} 
     {\@ifundefined{selectfont}
             {
                \message{We need AmS-LaTeX, this version of LaTeX does not 
                        support it}\@@end
             }{
                \def\mathcal{\cal}
                
                \def\pcyr{%
                        \def\default@family{UWCyr}%
                        \let\oldSl@\sl
                        \def\sl{\def\default@shape{it}\oldSl@}%
                        \cyracc
                        \language\Russian\family{UWCyr}\selectfont
                }
             }}{
                \DeclareFontEncoding{OT2}{}{\noaccents@}
                \DeclareFontSubstitution{OT2}{cmr}{m}{n}
                \DeclareFontFamily{OT2}{cmr}{\hyphenchar\font45 }
                \DeclareFontShape{OT2}{cmr}{m}{n}{%
                     <5><6><7><8><9><10>gen*wncyr %
                     <10.95><12><14.4><17.28><20.74><24.88> wncyr10 %
                }{}
                \DeclareFontShape{OT2}{cmr}{m}{it}{%
                     <5><6><7><8><9><10> gen * wncyi%
                     <10.95><12><14.4><17.28><20.74><24.88> wncyi10%
                }{}
                \DeclareFontShape{OT2}{cmr}{bx}{n}{%
                     <5><6><7><8><9><10> gen * wncyb%
                     <10.95><12><14.4><17.28><20.74><24.88> wncyb10%
                }{}
                \DeclareFontShape{OT2}{cmr}{m}{sl}{%
                     <-> ssub * cmr/m/it%
                }{}
                \DeclareFontShape{OT2}{cmr}{m}{sc}{%
                     <5><6><7><8><9><10>%
                     <10.95><12><14.4><17.28><20.74><24.88> wncysc10%
                }{}
                \DeclareFontFamily{OT2}{cmss}{\hyphenchar\font45 }
                \DeclareFontShape{OT2}{cmss}{m}{n}{%
                     <8><9><10> gen * wncyss%
                     <10.95><12><14.4><17.28><20.74><24.88> wncyss10%
                }{}
                
                \def\cyrencodingdefault{OT2}
                \def\pcyr{%
                        \cyracc
                        \let\encodingdefault\cyrencodingdefault
                        \language\Russian\fontencoding{OT2}\selectfont
                }
             }   

\@ifundefined{theorembodyfont} 
     {
        \def\theorembodyfont#1{\relax}
        \makeatletter
          \let\@@th@plain\th@plain
          \def\th@plain{ \@@th@plain \slshape }
        \makeatother
        \let\normalshape\relax
     }{}   

\ifx\cprime\undefined
     \def\cprime{$'$}
\fi
\makeatletter
  \def\@sect@my#1#2#3#4#5#6[#7]#8{%
\ifnum #2>\c@secnumdepth
   \let\@svsec\@empty
 \else
   \refstepcounter{#1}%
\edef\@svsec{\ifnum#2<\@m
             \@ifundefined{#1name}{}{\csname #1name\endcsname\ }\fi
\noexpand\rom{\csname the#1\endcsname.}\enspace}\fi
 \@tempskipa #5\relax
 \ifdim \@tempskipa>\z@ 
   \begingroup #6\relax
   \@hangfrom{\hskip #3\relax\@svsec}{\interlinepenalty\@M #8\par}%
   \endgroup
   \if@article\else\csname #1mark\endcsname{%
        \ifnum \c@secnumdepth >#2\relax\csname the#1\endcsname. \fi#7}\fi
\ifnum#2>\@m \else
       \let\@tempf\\ \def\\{\protect\\}\addcontentsline{toc}{#1}%
{\ifnum #2>\c@secnumdepth \else
             \protect\numberline{%
               \ifnum#2<\@m
               \@ifundefined{#1name}{}{\csname #1name\endcsname\ }\fi
               \csname the#1\endcsname.}\fi
           #8}\let\\\@tempf
     \fi
 \else
  \def\@svsechd{#6\hskip #3\@svsec
    \@ifnotempty{#8}{\ignorespaces#8\unskip
       \ifnum\spacefactor<1001.\fi}%
        \ifnum#2>\@m \else
          \let\@tempf\\ \def\\{\protect\\}\addcontentsline{toc}{#1}%
            {\ifnum #2>\c@secnumdepth \else
              \protect\numberline{%
                \ifnum#2<\@m
                \@ifundefined{#1name}{}{\csname #1name\endcsname\ }\fi
                \csname the#1\endcsname.}\fi
             #8}\let\\\@tempf\fi}%
 \fi
\@xsect{#5}}
  \ifx\RequirePackage\undefined
  \let\@sect\@sect@my             
  \fi
  \def\th@remark@my{\theorempreskipamount6\p@\@plus6\p@
    \theorempostskipamount\theorempreskipamount
    \def\theorem@headerfont{\it}\normalshape}
  \ifx\RequirePackage\undefined
  \let\th@remark\th@remark@my
  \fi
\let\myLabel\@gobble
\def\labelsONmargin{\@mparswitchfalse\def\myLabel##1{\@bsphack\marginpar
                                  {\normalshape\tiny\rm Label ##1}\@esphack}}
\ifx \url\undefined
  \def\url#1{{\tt #1}}%
\fi
\def\cyracc{\def\u##1{
                \if \i##1\char"1A%
                \else \if I##1\char"12%
                \else \accent"24 ##1\fi\fi }%
\def\"##1{\if e##1{\char"1B}%
                \else \if E##1{\char"13}%
                \else \accent"7F ##1\fi\fi }%
\def\9##1{\if##1z\char"19 
\else\if##1Z\char"11 
\else\if##1E\char"03 
\else\if##1e\char"0B 
\else\if##1u\char"18 
\else\if##1U\char"10 
\else\if##1A\char"17 
\else\if##1a\char"1F 
\else\if##1p\char"7E 
\else\if##1P\char"5E 
\else\if##1Q\char"5F 
\else\if##1q\char"7F 
\else\if##1i\char"1A 
\else\if##1I\char"12 
\else\if##1N\char"7D 
\fi
\fi
\fi
\fi
\fi
\fi
\fi
\fi
\fi
\fi
\fi
\fi
\fi
\fi
\fi
}%
\def\cydot{{\kern0pt}}}%
\def\cydot{$\cdot$}

\ifx\Russian\undefined
        \def\Russian{0\relax
    \message{Don't know the hyphenation rules for Russian^^J
                        Please do INITeX with `input  russhyph' in the 
                        command line}%
                \gdef\Russian{0\relax}%
        }
\fi
\def\@putname#1#2#3#4{\def\@@ref{#3}\let\old@bf\bf	
	\let\old@reset@font\reset@font			
        \def\bf##1{\old@bf\if?\noexpand##1?{#4}\else##1\fi}%
	\def\reset@font##1##2{\old@reset@font##1\if?\noexpand##2?{#4}\else##2\fi}#1{#2}%
        \let\bf\old@bf\let\reset@font\old@reset@font}
\let\my@ref=\ref
\def\ref#1{\@putname\my@ref{#1}{#1}{\tiny\rm\@@ref}}
\let\my@pageref=\pageref
\def\pageref#1{\@putname\my@pageref{#1}{#1}{\tiny\rm\@@ref}}
\let\my@cite=\cite
\def\cite#1{\@putname\my@cite{#1}{\@citeb}{\tiny\rm\@@ref}}
\makeatother
\fi
\relax 
\theoremstyle{plain} 
\theorembodyfont{\sl}

\numberwithin{equation}{section}
\theorembodyfont{\rm}
\theoremstyle{definition}

\newtheorem{definition}{Definition}[section]
\newtheorem{conjecture}[definition]{Conjecture}
\newtheorem{example}[definition]{Example}

\theoremstyle{remark}

\newtheorem{remark}[definition]{Remark} 

\theoremstyle{plain} 
\theorembodyfont{\sl}

\newtheorem{theorem}[definition]{Theorem}
\newtheorem{lemma}[definition]{Lemma}
\newtheorem{corollary}[definition]{Corollary}
\newtheorem{proposition}[definition]{Proposition}
\newtheorem{amplification}[definition]{Amplification}

\begin{document}
\bibliographystyle{amsplain}
\relax 

\title[Webs and bihamiltonian structures]{ Kronecker webs, bihamiltonian
structures, \\
and the method of argument translation }

\author{ Ilya Zakharevich }

\address{ Department of Mathematics, Ohio State University, 231 W.~18~Ave,
Columbus, OH, 43210 }

\email {ilya\atSign{}math.ohio-state.edu}

\date{ August 1999 (Revision III: March 2000)\quad Archived as
\url{math.SG/9908034}\quad Printed: \today }

\setcounter{section}{-1}

\maketitle
\begin{abstract}
We show that manifolds which parameterize values of first
integrals of integrable finite-dimensional bihamiltonian systems carry a
geometric structure which we call a {\em Kronecker web}. We describe two
opposite-direction functors between Kronecker webs and integrable
bihamiltonian structures, one is left inverse to the other.
Conjecturally, these two functors are mutually inverse (for ``small'' open
subsets of the manifolds in question).

The conjecture above is proven here when the bihamiltonian structure
allows an anti-involution of a particular form. This implies the
conjecture of \cite{GelZakh99Web} that on a dense open subset the
bihamiltonian structure on $ {\mathfrak g}^{*} $ is flat if $ {\mathfrak g} $ is semisimple, or if
$ {\mathfrak g}={\mathfrak G}\ltimes \operatorname{ad}_{{\mathfrak G}} $ and $ {\mathfrak G} $ is semisimple, and for some other Lie algebras of
mappings.

\end{abstract}
\tableofcontents

\section{Basic notions }\label{h002}\myLabel{h002}\relax 

We postpone the informal discussion of what is done in this paper
until Section~\ref{h0}, and start with introduction of notations and
conventions used throughout this text. People familiar with basic notions
and terminology of bihamiltonian geometry may jump directly to Section
~\ref{h0}, looking up Examples~\ref{ex00.30} and~\ref{ex002.45} on the ``when needed''
basis.

Many results of this paper may be stated in greater generality, but for
simplicity we assume that all the vector spaces we consider here are
finite-dimensional\footnote{With obvious exceptions of vector spaces of functions on manifolds.} vector spaces over a field $ {\mathbb K} $ which is either $ {\mathbb R} $ or $ {\mathbb C} $. A
{\em manifold\/} is a $ C^{\infty} $-manifold or a real-analytic manifold in the case $ {\mathbb K}={\mathbb R} $, and
an analytic manifold in the case $ {\mathbb K}={\mathbb C} $. We use the word {\em smooth\/} to mean
$ C^{\infty} $-smooth, real-analytic, or complex-analytic correspondingly.

For a vector space $ V $ over $ {\mathbb K} $ denote by $ V^{*} $ the space of $ {\mathbb K} $-linear
functionals on $ V $. Note that throughout this paper we do {\em not\/} consider
semilinear functionals or Hermitian forms on complex vector spaces.

We start by recalling some basic notions and notations of Poisson
geometry (see \cite{Arn89Math,GuiSte77Geo,ArnGiv85Sym}). In what follows
if $ f $ is a function or a tensor field on $ M $, $ f|_{m} $ denotes the value of $ f $ at
$ m\in M $.

\begin{remark} Throughout the paper we use standard idioms of differential
geometry. Say the phrase ``{\em at generic points\/}'' means ``at points of an
appropriate open dense subset''. Similarly, a ``{\em small open subset\/}'' is used
instead of ``an appropriate neighborhood of any given point''. A {\em local
isomorphism\/} between two geometric structures on $ M $ and $ M' $ is an
isomorphism of a neighborhood of a given point on $ M $ with a neighborhood
of a given point on $ M' $. Two geometric structures on $ M $ and on $ M' $ are
{\em locally isomorphic\/} if for any $ m\in M $ and $ m'\in M' $ there is a local isomorphism
which sends $ m $ to $ m' $. \end{remark}

\begin{definition} A {\em bracket\/} on a manifold $ M $ is a $ {\mathbb K} $-bilinear skew-symmetric
mapping which sends a pair of smooth functions\footnote{In the complex-analytic case one should consider functions on open
subsets $ U\subset M $ and require that the brackets on these subsets are compatible
on intersections.} $ f $, $ g $ on $ M $ to a smooth
function $ \left\{f,g\right\} $ on $ M $. This mapping is required to satisfy the Leibniz
identity $ \left\{f,gh\right\}=g\left\{f,h\right\}+h\left\{f,g\right\} $. A bracket is {\em Poisson\/} if it satisfies the
Jacobi identity $ \left\{f,\left\{g,h\right\}\right\}=\left\{\left\{f,g\right\},h\right\}+\left\{g,\left\{f,h\right\}\right\} $ (thus defines a structure
of a Lie algebra on functions on $ M $).

A {\em Poisson structure\/} is a manifold $ M $ equipped with a Poisson
bracket. \end{definition}

\begin{remark} \label{rem002.20}\myLabel{rem002.20}\relax  Leibniz identity implies $ \left\{f,g\right\}|_{m}=0 $ if $ f $ has a zero of
second order at $ m\in M $, or if $ f\equiv \operatorname{const} $. Thus a bracket is uniquely determined
by describing the functions $ \left\{f_{i},f_{j}\right\} $; here $ \left\{f_{i}\right\}_{i\in I} $ is an arbitrary
collection of functions on $ M $ which separates points of $ M $. Here we say
that a collection $ \left\{f_{i}\right\}_{i\in I} $ of smooth functions on $ M $ {\em separates points\/} of $ M $
if for any $ m\in M $ the collection $ \left\{df_{i}|_{m}\right\}_{i\in I} $ of vectors in $ {\mathcal T}_{m}^{*}M $ spans $ {\mathcal T}_{m}^{*}M $ as
a vector space. \end{remark}

\begin{definition} \label{def01.120}\myLabel{def01.120}\relax  Consider a bracket $ \left\{,\right\} $ on a manifold $ M $. The
{\em associated bivector\/}\footnote{A {\em bivector field\/} is a skew-symmetric contravariant tensor of valence 2.} {\em field\/} $ \eta $ is the section of $ \Lambda^{2}{\mathcal T}M $ given by $ \left\{f,g\right\}|_{m}=\left<
\eta|_{m},df\wedge dg|_{m} \right> $, $ m\in M $; here $ \left<, \right> $ denotes the canonical pairing between
$ \Lambda^{2}{\mathcal T}_{m}M $ and $ \Omega_{m}^{2}M $.

Given $ m_{0}\in M $, the {\em associated pairing\/} (,) in $ {\mathcal T}_{m_{0}}^{*}M $ is defined as
$ \left(\alpha,\beta\right)=\left\{f,g\right\}|_{m_{0}} $ if $ \alpha=df|_{m_{0}} $, $ \beta=dg|_{m_{0}} $. \end{definition}

Obviously, the associated bivector field uniquely determines the
bracket and visa versa. The associated pairing is a {\em skew-symmetric\/}
bilinear pairing.

\begin{definition} Given a skew-symmetric bilinear form (,) in a vector space
$ V $, put $ \operatorname{Ker}\left(,\right) \buildrel{\text{def}}\over{=} \left\{v\in V \mid \left(v,v'\right)=0\,\,\,\forall v'\in V\right\} $, and call $ \dim  \operatorname{Ker}\left(,\right) $ the
{\em corank\/} of (,). The {\em rank\/} of (,) is $ \dim  V-\dim  \operatorname{Ker}\left(,\right) $.

The {\em rank\/} of the bracket $ \left\{,\right\} $ at $ m\in M $ is $ r $ if the associated
skew-symmetric bilinear pairing on $ {\mathcal T}_{m}^{*}M $ has rank $ r $. In this case the
{\em corank\/} of the bracket is $ \dim  M-r $.

A bracket has a {\em constant (co)rank\/} if its rank does not depend on the
point $ m\in M $. A Poisson bracket is {\em symplectic\/} if the corank is constant and
equal to 0. \end{definition}

The associated tensor field $ \eta $ of a bracket on $ M $ can be considered as
a mapping $ {\text H}\colon {\mathcal T}^{*}M \to {\mathcal T}M $ (the {\em Hamiltonian mapping\/} of a bracket). If the
bracket is symplectic, this mapping is invertible, and the inverse
mapping $ {\text H}^{-1}\colon {\mathcal T}M \to {\mathcal T}^{*}M $ can be considered as a bilinear pairing on $ {\mathcal T}M $, or
as a tensor field. This tensor field is a section $ \omega $ of $ \Omega^{2}M $ of corank 0,
called the {\em symplectic\/} $ 2 $-{\em form\/} of the symplectic bracket. In local
coordinates the tensors $ \eta $ and $ \omega $ are given by mutually inverse matrices.

In the other direction, given a section $ \omega $ of $ \Omega^{2}M $ of corank 0,
setting $ \eta=\omega^{-1} $ gives a bracket on $ M $. It is easy to check that $ \eta $ is Poisson
iff\footnote{This condition is linear in $ \omega $, as opposed to the quadratic condition (of
Jacobi identity) on $ \eta $. This linearity makes it much easier to study symplectic
Poisson brackets.} $ d\omega=0 $.

\begin{example} \label{ex002.41}\myLabel{ex002.41}\relax  Given a manifold $ N $, put $ M={\mathcal T}^{*}N $, and let $ \pi\colon M \to N $ be the
natural projection. Given $ m\in M $, one can write $ m=\left(n,\nu\right) $, $ n\in N $, $ \nu\in{\mathcal T}_{n}^{*}N $. Apply
$ \pi^{*}\colon {\mathcal T}_{n}^{*}N \to {\mathcal T}_{m}^{*}M $ to $ \nu $, and note that $ \pi^{*}\nu $ is an element of $ {\mathcal T}_{m}^{*}M $ which
depends on $ m $ only. One can write $ \pi^{*}\nu=\alpha\left(m\right) $; here $ \alpha $ is a canonically
defined section of $ \Omega^{1}M $.

Local coordinates $ \left(n_{1},\dots ,n_{d}\right) $ on $ N $ define local coordinates
$ \left(n_{1},\dots ,n_{d},\nu_{1},\dots ,\nu_{d}\right) $ on $ {\mathcal T}^{*}N $. In these coordinates $ \alpha=\sum_{i}\nu_{i}dn_{i} $. Taking
$ \omega=d\alpha $, one obviously gets $ d\omega=0 $. In local coordinates $ \omega=\sum_{i}d\nu_{i}\wedge dn_{i} $, hence $ \omega $
is of corank 0, consequently defines a (symplectic) Poisson structure $ \eta $ on
$ M={\mathcal T}^{*}N $. \end{example}

Recall that any symplectic Poisson structure is locally isomorphic
to the structure of Example~\ref{ex002.41}.

\begin{definition} Call two Poisson brackets $ \left\{,\right\}_{1} $ and $ \left\{,\right\}_{2} $ on $ M $ {\em compatible\/} if
the bracket $ \lambda_{1}\left\{,\right\}_{1} +\lambda_{2}\left\{,\right\}_{2} $ is Poisson for any $ \lambda_{1} $, $ \lambda_{2} $.

A {\em bihamiltonian structure\/} is a manifold $ M $ with a pair of compatible
Poisson brackets. \end{definition}

Given that the Jacobi identity is ``quadratic'' in $ \left\{,\right\} $, one can show
that if {\em one\/} linear combination $ \lambda_{1}\left\{,\right\}_{1} +\lambda_{2}\left\{,\right\}_{2} $ of two Poisson brackets is
Poisson, and $ \lambda_{1}\not=0 $, $ \lambda_{2}\not=0 $, then {\em any\/} linear combination $ \lambda_{1}\left\{,\right\}_{1} +\lambda_{2}\left\{,\right\}_{2} $ is
Poisson. Even if $ M $ is a $ C^{\infty} $-manifold, the coefficients $ \lambda_{1} $, $ \lambda_{2} $ may be
taken to be complex numbers. Indeed, if $ M $ is a $ C^{\infty} $-manifold with a
bracket, one may consider the extension of the bracket to the $ {\mathbb C} $-vector
space of complex-valued functions on $ M $. In this case $ \lambda_{1}\left\{,\right\}_{1} +\lambda_{2}\left\{,\right\}_{2} $ is
well-defined even for complex values of $ \lambda_{1},\lambda_{2} $. By the remarks above,
complex linear combinations of brackets of a bihamiltonian structure are
also Poisson. In what follows we can always consider brackets as acting
on the spaces of complex-valued functions.

Given a pair of brackets $ \left\{,\right\}_{1} $ and $ \left\{,\right\}_{2} $, one obtains two bivector
fields $ \eta_{1} $, $ \eta_{2} $. Analogously, one obtains two skew-symmetric bilinear
pairings $ \left(,\right)_{1} $, $ \left(,\right)_{2} $ on $ {\mathcal T}_{m}^{*}M $, so that $ \left(\alpha,\beta\right)_{j}=\left\{f,g\right\}_{j}|_{m} $ if $ \alpha=df|_{m} $, $ \beta=dg|_{m} $,
$ j=1,2 $.

\begin{definition} \label{def00.26}\myLabel{def00.26}\relax  A Poisson structure $ \left\{,\right\} $ on a vector space $ V $ is
{\em translation-invariant\/} if for any parallel translation $ {\mathfrak T}\colon V \to V $ and any
two functions $ f $, $ g $ on $ V $ one has $ {\mathfrak T}^{*}\left\{f,g\right\}=\left\{{\mathfrak T}^{*}f,{\mathfrak T}^{*}g\right\} $. A bihamiltonian
structure on $ V $ is {\em translation-invariant\/} if both Poisson brackets are
translation-invariant.

A bihamiltonian structure on a manifold $ M $ is {\em flat\/} if it is {\em locally
isomorphic\/} to a translation-invariant bihamiltonian structure. In other
words, for any $ m\in M $ there is a neighborhood $ U\ni m $ such that the restriction
of the bihamiltonian structure to $ U $ is isomorphic to the restriction of an
appropriate translation-invariant bihamiltonian structure to an
appropriate open subset. \end{definition}

The tensor field $ \eta $ of a translation-invariant Poisson bracket on $ V $
has constant coefficients in any vector-space coordinate system on $ V $. For
us one particular example of a translation-invariant bihamiltonian
structure is of special interest.

\begin{example} \label{ex00.30}\myLabel{ex00.30}\relax  Consider a vector space $ V $ with coordinates $ x_{0},\dots ,x_{2k-2} $.
Define Poisson brackets $ \left\{,\right\}_{1} $ and$ \left\{,\right\}_{2} $ on the coordinate functions by
\begin{equation}
\left\{x_{2l},x_{2l+1}\right\}_{1}=1,\qquad \left\{x_{2l+1},x_{2l+2}\right\}_{2}=1,\qquad 0\leq l\leq k-2,
\label{equ45.20}\end{equation}\myLabel{equ45.20,}\relax 
and setting any other bracket of the coordinate functions $ x_{0},\dots ,x_{2k-2} $ to
be zero. This pair of brackets is in fact a translation-invariant
bihamiltonian structure. \end{example}

The following example is the simplest of the classical examples of
bihamiltonian structures arising in the theory of integrable systems.

\begin{example} \label{ex002.45}\myLabel{ex002.45}\relax  Given a Lie algebra $ {\mathfrak g} $ and an element $ c_{1}\in{\mathfrak g}^{*} $, define a
bihamiltonian structure on $ {\mathfrak g}^{*} $ as in \cite{MishFom78Eul,Bol91Com}. An
element $ X\in{\mathfrak g} $ defines a linear function $ f_{X} $ on $ {\mathfrak g}^{*} $. Due to Remark~\ref{rem002.20},
to define a bihamiltonian structure on $ {\mathfrak g}^{*} $ it is enough to describe
the brackets $ \left\{f_{X},f_{Y}\right\}_{j} $, $ j=1,2 $, $ X,Y\in{\mathfrak g} $.

Let $ \left\{f_{X},f_{Y}\right\}_{1} $ be the {\em constant\/} function on $ {\mathfrak g}^{*} $ and $ \left\{f_{X},f_{Y}\right\}_{2} $ be the
{\em linear\/} function on $ {\mathfrak g}^{*} $ given by the formulae
\begin{equation}
\left\{f_{X},f_{Y}\right\}_{1}=f_{\left[X,Y\right]}\left(c_{1}\right),\qquad \left\{f_{X},f_{Y}\right\}_{2}=f_{\left[X,Y\right]}.
\notag\end{equation}
The bracket $ \left\{,\right\}_{2} $ is the natural Lie--Kirillov--Kostant--Souriau Poisson
bracket on $ {\mathfrak g}^{*} $. The bracket $ \left\{,\right\}_{1} $ is translation-invariant. The bracket
$ \left\{,\right\}_{2} $ is translation-invariant only if $ {\mathfrak g} $ is abelian (and then $ \left\{,\right\}_{1,2} $
vanish).

In fact, instead of taking $ c_{1}\in{\mathfrak g}^{*} $ one can consider any $ 2 $-cocycle
$ c_{2}\in\Lambda^{2}{\mathfrak g}^{*} $, and define $ \left\{f_{X},f_{Y}\right\}_{1}=c_{2}\left(X,Y\right) $. The definition above is recovered
if one puts $ c_{2}=\partial c_{1} $. If $ H^{2}\left({\mathfrak g},{\mathbb C}\right)=0 $ then these two versions of the
definition are equivalent. \end{example}

\begin{definition} A smooth function $ F $ on a manifold $ M $ with a Poisson bracket
$ \left\{,\right\} $ is a {\em Casimir\/} function if $ \left\{F,f\right\}=0 $ for any smooth function $ f $ on $ M $. \end{definition}

Obviously, any function $ \varphi\left(F_{1},F_{2},\dots ,F_{k}\right) $ of several Casimir functions
is again Casimir.

\begin{definition} A collection of smooth functions $ F_{1},\dots ,F_{r} $ on $ M $ is {\em dependent\/}
if $ \varphi\left(F_{1},\dots F_{r}\right)\equiv 0 $ for an appropriate smooth function $ \varphi $ such that $ \varphi\not\equiv 0 $ on
any open subset of $ {\mathbb K}^{r} $. \end{definition}

\begin{remark} \label{rem00.51}\myLabel{rem00.51}\relax  Consider a manifold $ M $ with a Poisson bracket $ \left\{,\right\} $. The
local classification of Poisson structures of {\em constant rank\/} \cite{Kir76Loc,%
Wei83Loc} shows that for an arbitrary Poisson bracket there is an open
(and in interesting cases dense) subset $ U\subset M $ and $ k\in{\mathbb Z}_{\geq0} $ such that on $ U $
there are $ k $ independent Casimir functions $ F_{1},\dots ,F_{k} $, and any Casimir
function on $ U $ may be written as a function of $ F_{1},\dots ,F_{k} $ (we do not
exclude the case $ k=0 $). The common level sets $ F_{1}=C_{1},\dots ,F_{k}=C_{k} $ form an
invariantly defined foliation on $ U $, which is called the {\em symplectic
foliation}. Moreover, locally one can construct additional functions $ n_{i} $,
$ \nu_{i} $, $ i=1,\dots ,d $, on $ U $ such that $ \left\{n_{i},n_{j}\right\}=\left\{\nu_{i},\nu_{j}\right\}=0 $, $ \left\{n_{i},\nu_{j}\right\}=\delta_{ij} $, and the
functions $ F_{\bullet} $, $ n_{\bullet} $, $ \nu_{\bullet} $ form a coordinate system on $ U $. \end{remark}

This shows that any Poisson structure of constant rank is {\em flat},
i.e., locally isomorphic to a translation-invariant Poisson structure.
Consequently, any analytic Poisson structure is flat on an open dense
subset. Moreover, the leaves of the symplectic foliation of such Poisson
structures can be simultaneously equipped with coordinates as in Example
~\ref{ex002.41}.

\begin{definition} Consider a foliation $ {\mathcal F} $ on $ B $. Define the {\em tangent bundle\/} $ {\mathcal T}{\mathcal F} $ to
$ {\mathcal F} $ to be the vector subbundle $ E $ of $ {\mathcal T}B $ such that $ E_{b} $ coincides with $ {\mathcal T}_{b}L_{b} $
for any $ b\in B $; here $ L_{b} $ is the leaf of $ {\mathcal F} $ which passes through $ b $. Let the
{\em normal bundle\/} $ {\mathcal N}{\mathcal F}\subset{\mathcal T}^{*}B $ to $ {\mathcal F} $ be the orthogonal complement of $ {\mathcal T}{\mathcal F} $, and the
{\em cotangent bundle\/} $ {\mathcal T}^{*}{\mathcal F} $ to $ {\mathcal F} $ be the vector bundle dual to $ {\mathcal T}{\mathcal F} $. \end{definition}

The total space of the bundle $ {\mathcal T}^{*}{\mathcal F} $ is a union of the total spaces of
the cotangent bundles $ {\mathcal T}^{*}L $ of the leaves $ L $ of the foliation. Since each
$ {\mathcal T}^{*}L $ carries a natural symplectic Poisson structure (Example~\ref{ex002.41}),
one obtains

\begin{proposition} The total space of the cotangent bundle to a foliation
carries a natural Poisson structure of constant rank. \end{proposition}

Moreover, the foliation on $ {\mathcal T}^{*}{\mathcal F} $ specified above is the symplectic
foliation of this Poisson structure. Due to Remark~\ref{rem00.51}, any Poisson
structure of constant rank is locally isomorphic to the Poisson structure
on $ {\mathcal T}^{*}{\mathcal F} $ of a foliation $ {\mathcal F} $ with the appropriate dimension and rank.

\section{Introduction }\label{h0}\myLabel{h0}\relax 

Among other approaches to integrable systems\footnote{An {\em integrable system\/} is a notion which is used in different senses in
different situations. For the sake of this introduction, one can
understand it as a {\em Liouville integrable system}. In other words, it is a
system of ODE which can be solved explicitly, and such that the
solutions demonstrate no mixing property: trajectories stay in closed
submanifolds of ``small'' dimension.} the so-called
{\em bihamiltonian approach\/} is especially interesting from the geometric point
of view. In this approach all the properties of an integrable
bihamiltonian system are deduced from the bihamiltonian structure on the
phase manifold\footnote{There is a widespread belief that most (or all) integrable systems which
arise in problems of mathematical physics allow a natural bihamiltonian
structure.}. (The principal tool for this deduction is {\em Lenard scheme},
which is outside of the scope of this paper, see \cite{Lax76Alm,Mag78Sim,%
GelDor79Ham,FokFuch80Str,KosMag96Lax,GelZakh99Web}.) Since the
structure of a bihamiltonian system is nothing more than a pair of tensor
fields satisfying some invariantly-defined conditions, this approach puts
the integrable system into the standard framework of differential
geometry.

A powerful approach to study a particular object of differential
geometry is to classify {\em all possible\/} objects up to isomorphism, and to
describe the automorphisms of every object. After this what remains is to
describe the position of the given object inside the classification. For
example, in the case of symplectic structures or Poisson structures of
constant rank, the local classification is ``trivial'': locally there are
only some discrete parameters, thus any structure is locally isomorphic
to one from a finite list (for a given dimension). This reduces all
questions on geometry of symplectic manifolds to questions of global
nature. Similarly, all questions on geometry of Poisson manifolds are
reduced to questions of global nature, and to questions related to
subsets where the rank drops.

{\bf Classification in general position}. Starting from the end of the
80s, there were many results on local classifications of bihamiltonian
structures ``in general position''. The spirit of these classification
efforts follows the results on classification of the corresponding
linearized problem: classification of pairs of skew-symmetric bilinear
forms in a vector space. The pairs in general position behave differently
depending on the dimension: an even-dimensional pair in general position
is a direct sums of $ 2 $-dimensional blocks, while an odd-dimensional pair
in general position is an indecomposable block (a {\em Kronecker block\/}; for
details see Section~\ref{h37}).

For a long time the case of general position on even-dimensional
manifolds was considered to be the most relevant to the theory of
integrable systems. This classification was done in \cite{Tur89Cla,%
Mag88Geo,Mag95Geo,McKeanPC,GelZakh93} (under different
assumptions). The case of general position in the odd-dimensional case
was analyzed in \cite{GelZakhWeb,GelZakh93,Rig95Tis,Rig95Geo,%
Rig98Sys,Tur99Equi,Tur99MemA}. However, until recently these
results had little direct impact on the theory of integrable systems,
since it was not known which classical integrable systems are subject to
these conditions of general position. To answer this question, one would
need to classify a {\em linearized\/} bihamiltonian structure in a generic point
of a particular classical integrable system.

{\bf Applications to integrable systems}. It was a very old idea of
I.~M.~Gelfand that a key to the bihamiltonian geometry might be uncovered
by studying the linearized bihamiltonian structure of the periodic KdV
system. In the mid-80s (i.e., prior to all of the works mentioned above), the
author started a joint work with I.~M.~Gelfand on this problem. This
turned out to be a very hard problem of functional analysis, and the
first results \cite{GelZakhFAN} were achieved only after several years of
intensive work. However, the fruits of this work deserved the effort, as the
results of \cite{GelZakhFAN,GelZakh94Spe} were completely orthogonal to the
paradigms of the time.

It turned out that the periodic KdV system should be considered as
an infinite-dimensional analogue of an {\em odd-dimensional\/} bihamiltonian
structure in general position. The analogous finite-dimensional pairs of
skew-symmetric forms are exactly the Kronecker blocks mentioned above,
and the Lenard scheme mentioned above works for any such pair. It was
this analogy which fueled us to investigate the geometry of bihamiltonian
structures linearizations of which at any point are Kronecker blocks;
this we carried out in \cite{GelZakhWeb,GelZakh93}.

Already at that time we knew an example of a classical bihamiltonian
system which is subject to the classification above: the open Toda
lattice (for definition see, for example, \cite{FadTakh87Ham}). To some
extent, Toda lattices are finite-dimensional analogues of the KdV system,
so one could have expected that the periodic Toda lattice might have
properties similar to a periodic KdV system. Unfortunately, the periodic
Toda lattice is an even-dimensional bihamiltonian system, thus the direct
analogy could not work. However, an open Toda lattice {\em is\/} an
odd-dimensional bihamiltonian system, and it is relatively easy to check
that the linearization at a point in general position is a Kronecker
block.

{\bf Cases of non}-{\bf general position}. We had one example where Kronecker
blocks played a vital role, but this was not enough to change a paradigm.
What was needed was to have more examples, in particular, to investigate
what happens with the periodic Toda lattice. Again, it took a lot of time
until an understanding of the situation emerged in \cite{GelZakh99Web}. Near
points in general position the periodic Toda lattice turned out to be
locally isomorphic to a {\em product of two copies\/} of odd-dimensional
bihamiltonian manifolds, each of them subject to the classification
above.\footnote{In particular, a linearization at a point of the periodic Toda lattice
is {\em not\/} in general position.} Additionally, we could demonstrate similar results for
linearizations of other classical bihamiltonian systems: the
linearizations break down into a direct sum of Kronecker blocks, i.e.,
into a direct sum of odd-dimensional components.

{\bf Flatness}. These results show that Toda lattices are subject to the
local classification of odd-dimensional bihamiltonian systems of
\cite{GelZakhWeb,GelZakh93} (possibly after splitting into a direct
product). Let us stress that there are infinitely many non-isomorphic
bihamiltonian systems of these types (with parameters being several
functions of two variables), thus it is meaningful to ask {\em which\/} system of
the classification above is the open Toda lattice (or a factor of the
periodic Toda lattice).

Recall that the simplest possible odd-dimensional bihamiltonian
structure (in general position) is the structure of Example~\ref{ex00.30}. One
of the principal results of \cite{GelZakh99Web} is that the bihamiltonian
structure of the open Toda lattice {\em is locally isomorphic\/}\footnote{Contrast this with the property of having isomorphic linearizations
which we discussed before.} (on an open
dense subset) to the structure of Example~\ref{ex00.30} (here $ k $ is the number
of ``atoms'' in the Toda lattice). Moreover, the bihamiltonian structure of
the periodic Toda lattice is locally isomorphic (on an open dense subset)
to a direct product of two copies of the structure of Example~\ref{ex00.30}.
(If the periodic Toda lattice has $ n $ atoms, then one copy has $ k=n $, another
$ k=1 $.)

Extending this observation, one of the conjectures of \cite{GelZakh99Web}
says that other integrable systems of mathematical physics are also
locally isomorphic to a direct product of several copies of Example
~\ref{ex00.30}. In other words, these systems are flat on an open dense
subsets. By the reasons which will be explained shortly, here we call
this conjecture the {\em metaconjecture}.

There are similarities and differences of this metaconjecture with
the result on flatness on an open dense subset of analytic Poisson
manifolds. As in the Poisson case, this metaconjecture reduces questions
on geometry of these integrable systems to two questions: the description
of the behaviour at the points outside the dense open subsets mentioned
above (points where the bihamiltonian structure {\em degenerates\/}), and the
description of the gluing process. The question of local geometry in
generic points mostly disappear: the system is locally isomorphic to one
from a finite list\footnote{For example, due to the result above they are locally isomorphic to
a direct product of several open Toda lattices.} (for a fixed dimension).

The principal difference when comparing with the Poisson case is
that---as pointed out already---{\em in general\/} bihamiltonian geometry is not
locally flat. The metaconjecture above is a {\em selection principle\/}\footnote{For example, any proof of this metaconjecture (if possible) would need
to concentrate on the question {\em why\/} mathematical physicists study some
systems and do not study some other systems.}: out of
a huge variety of different integrable bihamiltonian systems the systems
studied in mathematical physics fall into a very thin subclass of {\em flat
bihamiltonian systems}.

The paper \cite{GelZakh99Web} also lists some specializations of the
metaconjecture above: it is conjectured that some particular
bihamiltonian systems of mathematical physics are flat on an open dense
subset. This list includes the full Toda lattice (for definition see
\cite{Kos79Sol}), the multi-dimensional Euler top (for definition see
\cite{MorPiz96Eul}), and the semisimple case of Example~\ref{ex002.45}, and some
other examples related to Toda lattices.

Here is the current status of these special-case conjectures: in
addition to the already proven case of various Toda lattices
(\cite{GelZakh99Web}), in this paper we establish the semisimple case of
Example~\ref{ex002.45} (see Corollary~\ref{cor11.80}). This provides several
established examples of ``classical'' integrable bihamiltonian systems
which locally look like a product of the structures of Example~\ref{ex00.30}.

{\bf Integrability}. The examples above show that Kronecker blocks play an
important role in geometry of classical integrable systems. Compare this
with another reference point: \cite{Bol91Com} formalized the notion of
integrability of a bihamiltonian system under the name of {\em completeness}.
This notion is very close to the property of being {\em micro-Kronecker\/}\footnote{Different faces of the relationship of Kronecker blocks with integrability
were independently discovered in \cite{Pan99Ver} and in \cite{GelZakh99Web}.} that
we introduce below, in Section~\ref{h4}: a bihamiltonian system is {\em complete\/} if
linearizations at points of an open dense subset break into a direct sum
of Kronecker blocks.

All of the examples above (as most of the other bihamiltonian systems
of mathematical physics) happen to be integrable in this strict sense. In
this informal introduction we always use the word {\em integrable\/} in the sense
of being micro-Kronecker. A system is complete if it integrable on an
open dense subset. This makes it very important to investigate the local
geometry of arbitrary integrable bihamiltonian systems. This
investigation is one of the principal targets of this paper.

{\bf Webs}. We mentioned two kinds of bihamiltonian structures with
relationship to Kronecker blocks: odd-dimensional bihamiltonian
structures in general position, and micro-Kronecker (or integrable)
bihamiltonian structures. The former ones are particular cases of the
latter ones, the condition being of having exactly one Kronecker block.
One of principal results of \cite{GelZakhWeb} is the introduction of a notion
of {\em Veronese web\/} (which is a family of foliations with appropriate
compatibility conditions). As \cite{GelZakhWeb} and \cite{Tur99Equi,Tur99MemA}
show, odd-dimensional bihamiltonian structures in general position are
locally classified by a Veronese web on a manifold of (approximately)
half the dimension of the initial manifold.

In \cite{Pan99Ver} (and independently---more generally, but in less
detail---in \cite{GelZakh99Web}) the construction of the Veronese web was
generalized to the case of structures which are integrable in the sense
of \cite{Bol91Com}. This paper starts with introduction of a geometric
structure of a {\em Kronecker web}, which is simultaneously a generalization of
the construction of a {\em Veronese web of higher codimension\/} of \cite{Pan99Ver}
and a more structured variant of the construction of a {\em web\/} of
\cite{GelZakh99Web}. Each Kronecker web has a {\em rank}, and the Veronese webs of
\cite{GelZakhWeb} coincide with Kronecker webs of rank 1.

Though the definition of a Kronecker web has no similarity with the
definition of Veronese webs, results of Sections~\ref{h33} imply that
Kronecker webs can also be described as families of foliations with
appropriate compatibility conditions.

{\bf Two functors and the principal conjecture}. Similarly to what was
done in \cite{GelZakhWeb} in the case of rank 1, we show how to associate to
any integrable bihamiltonian structure its Kronecker web, and show how to
construct an integrable bihamiltonian structure from an arbitrary
Kronecker web. Conjecturally (Conjecture~\ref{con8.50}), these functors are
mutually inverse (as in the case of rank 1, see \cite{GelZakhWeb,%
Tur99Equi,Tur99MemA}), but here we prove only that one is left
inverse to another.

If these functors were mutually inverse, the question of local
classification of integrable bihamiltonian structures would be reduced to
the question of local classification of Kronecker webs. Analogously to
Definition~\ref{def00.26}, one can define a {\em translation-invariant\/} Kronecker
web, and a {\em flat\/} Kronecker web. In particular, to show that a given
bihamiltonian structure is flat, it would be enough to show that its
Kronecker web is flat\footnote{As \cite{GelZakhWeb,GelZakh93,GelZakh99Web} show, there are plenty of
examples of non-flat Kronecker webs, even in the case of rank 1.}.

For example, in \cite{Pan99Ver} it is shown that in the semisimple case
of Example~\ref{ex002.45} the corresponding Kronecker web is flat on an open
dense subset. Together with our construction of two functors the
conjecture above would immediately show that the bihamiltonian structure
of Example~\ref{ex00.30} is flat on an open dense subset. Another example, the
(known) case of rank 1 of this conjecture is used in \cite{GelZakh99Web} to
demonstrate flatness of Toda lattices.

Moreover, note that proofs of flatness of webs of particular
bihamiltonian systems are very simple (see \cite{Pan99Ver,GelZakh99Web},
and Theorem~\ref{th11.70}). Thus the theory of Kronecker webs allows one to
condense all the problems of geometric classification of bihamiltonian
structures into a proof of Conjecture~\ref{con8.50} (or an appropriate
particular case of this conjecture).

{\bf Anti-involutions and Lie algebras}. In this paper we prove only a
very special case of Conjecture~\ref{con8.50}. Section~\ref{h91} introduces a
special subclass of bihamiltonian structures, structures which allow an
anti-involution of a special form. In Section~\ref{h10} we show that in the
case of such structures Conjecture~\ref{con8.50} holds (on a large open
subset).

Section~\ref{h11} uses this approach to show that the semisimple case of
Example~\ref{ex002.45} is flat on an open dense subset. As already explained,
the flatness of the Kronecker web is already proved in \cite{Pan99Ver}. We use
standard tools of the theory of semisimple Lie algebras to construct
appropriate anti-involutions. This provides an $ \operatorname{ad} $ hoc way to prove this
particular case of Conjecture~\ref{con8.50}.

{\bf Conclusions}. This paper presents several important steps on the road
to understanding geometry of bihamiltonian systems of mathematical
physics. We introduce a notion of a {\em Kronecker web\/} (Section~\ref{h34}), and of
a {\em micro-Kronecker\/} bihamiltonian structure (Section~\ref{h4}). Given a
Kronecker web $ B $, there is a canonically defined vector bundle $ \Phi $ over $ B $
such that the total space $ \Phi\left(B\right) $ of this bundle carries a canonically
defined micro-Kronecker bihamiltonian structure (Corollary~\ref{cor33.99}).
Given a micro-Kronecker bihamiltonian structure on a manifold $ M $, there is
a canonically defined foliation $ {\mathcal F} $ on $ M $ (the {\em action foliation\/}) such that
the local base $ {\mathcal B}_{M} $ of this foliation has a canonically defined structure
of a Kronecker web (Theorem~\ref{th4.55}). The leaves of the action foliation
on $ \Phi\left(B\right) $ are fibers of the projection $ \Phi\left(B\right) \to B $; thus $ {\mathcal B}_{\Phi\left(B\right)}\simeq B $. Moreover,
this isomorphism is compatible with the Kronecker web structures on $ {\mathcal B}_{\Phi\left(B\right)} $
and on $ B $ (Proposition~\ref{prop8.27}).

We conjecture that the ``inverse'' is also true, at least in the
following sense. Given a point $ m $ on a manifold $ M $ with a micro-Kronecker
bihamiltonian structure, there is a Kronecker structure on $ {\mathcal B}_{M} $, thus a
micro-Kronecker bihamiltonian structure on $ \Phi\left({\mathcal B}_{M}\right) $. We conjecture that
there is a local isomorphism between bihamiltonian structures on a
neighborhood of $ m $ and a neighborhood of $ \left(b,0\right) $; here $ b $ is the image of $ m $
in $ {\mathcal B}_{M} $, and $ \left(b,0\right) $ is the point on the $ 0 $-section of $ \Phi\left({\mathcal B}_{M}\right) $ over $ b $. This
constitutes Conjecture~\ref{con8.50}. (As the results of \cite{GelZakhWeb} show,
this local isomorphism cannot be canonical in non-trivial cases.) This
conjecture would reduce local classification of bihamiltonian structures
up to isomorphism to a local classification of Kronecker webs up to
isomorphism.

We prove a particular case of Conjecture~\ref{con8.50}, the case when the
bihamiltonian structure on $ M $ allows an anti-involution of a special form
(Section~\ref{h10}). We apply this knowledge to several important particular
cases of Example~\ref{ex002.45}: the cases of semisimple $ {\mathfrak g} $, and cases when $ {\mathfrak g} $
is a Lie algebra of mappings from a self-dual Artin scheme (i.e., a
Gorenstein scheme) to a semisimple Lie algebra. We show that in these
cases the bihamiltonian structure of Example~\ref{ex002.45} is locally
isomorphic to a direct product of several copies of Example~\ref{ex00.30}
(when restricted on a dense open subset of $ {\mathfrak g}^{*} $). This establishes one of
conjectures of \cite{GelZakh99Web}.

We believe that the same anti-involution trick will also work with
complete Toda lattice and the multi-dimensional Euler top. However, we do
not know whether this approach would work with the open and periodic Toda
lattices. (Recall that the proof of \cite{GelZakh99Web} used the established
case of rank 1 of Conjecture~\ref{con8.50}.)

Section~\ref{h92} shows that the complex variant of Conjecture~\ref{con8.50}
implies the real-analytic case as well. The appendix (Section~\ref{h99})
describes which geometric structures one can invariantly associate to a
Kronecker web.

With pleasure we thank V.~Serganova for fruitful discussions and for
important contributions: the proof of Proposition~\ref{prop11.75}, and a
significant simplification of the proof of Amplification~\ref{amp33.20}. This
paper could not have appeared without the long-term joint research with
I.~M.~Gelfand, which molded the mindset of this paper. Special thanks go
to A.~Panasyuk for providing preprints of \cite{Pan98Sym} and \cite{Pan99Ver}, and
to B.~Okun for patient discussions of results of \cite{GelZakh99Web}, which
lead to a cristallization of crucial ideas of the current paper, and to
Y.~Flicker for tireless suggestions for improvement.

{\bf Revision history. }The revisions of this paper are archived as
\url{math.SG/9908034} on the math preprint server archive
\url{arXiv.org/abs/math}. In addition to cosmetic changes, revision II of
this paper contains a major simplification of arguments in Section~\ref{h33},
adds Remarks~\ref{rem1.50} and~\ref{rem1.60}. Section~\ref{h11} is expanded beyond
Corollary~\ref{cor11.80} by adding the discussion of the case of groups of
mappings. Revision III reworked the introduction. The numeration of
statements did not change.

\section{Linear relations and pencils }\label{h1}\myLabel{h1}\relax 

\begin{definition} A {\em linear relation\/} $ W $ between vector spaces $ V_{1} $ and $ V_{2} $ is a
vector subspace $ W\subset V_{1}\oplus V_{2} $. Call $ W $ {\em bisurjective\/} if both projections of $ W $
to $ V_{1} $ and $ V_{2} $ are surjective. The {\em left kernel\/} $ \operatorname{Ker}_{L}W $ of $ W $ is the
kernel of the projection $ W \to V_{2} $ considered as a vector subspace of
$ V_{1} $. The {\em right kernel\/} $ \operatorname{Ker}_{R}W\subset V_{2} $ is defined similarly.

A {\em linear relation in\/} $ V $ is a linear relation between $ V $ and $ V $. \end{definition}

Obviously, if $ W $ is a linear relation, then $ W\supset\operatorname{Ker}_{L}W\oplus\operatorname{Ker}_{R}W $, thus $ W $
induces a vector subspace $ \widetilde{W} $ of $ \left(V_{1}/\operatorname{Ker}_{L}W\right)\oplus\left(V_{2}/\operatorname{Ker}_{R}W\right) $. If $ W $ is bisurjective,
then $ \widetilde{W} $ is a graph of a bijective linear mapping $ V_{1}/\operatorname{Ker}_{L}W \overset{\sim}\to V_{2}/\operatorname{Ker}_{R}W $. In
particular, $ \dim  V_{1}-\dim  V_{2}=\dim  \operatorname{Ker}_{L}W-\dim  \operatorname{Ker}_{R}W $.

\begin{definition} \label{def1.20}\myLabel{def1.20}\relax  Fix once and for all a two-dimensional vector space $ {\mathcal S} $
with a basis $ {\mathbit s}_{1} $, $ {\mathbit s}_{2} $. A {\em pencil\/} of linear operators between vector spaces $ V $
and $ W $ is a linear mapping $ V\otimes{\mathcal S} \xrightarrow[]{{\mathcal P}} W $. This induces two linear mappings
$ {\mathcal P}_{1,2}\colon V \to W $ defined as $ {\mathcal P}_{i}\left(v\right)\buildrel{\text{def}}\over{=}{\mathcal P}\left(v\otimes{\mathbit s}_{i}\right) $. \end{definition}

Given a pencil $ {\mathcal P} $, one obtains a linear $ 2 $-parametric family of linear
mappings $ \lambda_{1}{\mathcal P}_{1}+\lambda_{2}{\mathcal P}_{2} $ between $ V $ and $ W $. Inversely, any two linear mappings
$ {\mathcal P}_{1,2}\colon V \to W $ correspond to a pencil and to a linear $ 2 $-parametric family
of linear mappings.

Analyze possible connections between pencils and linear relations in
vector spaces. Given a linear relation $ W $ in $ V $, two projections $ \pi_{1,2} $ of
$ W\subset V\oplus V $ to $ V $ define a pencil $ \pi $ of operators $ W \to V $. In the other direction,
given any pencil $ \Pi\colon W \to V $, one can construct a linear relation $ \widetilde{W}\subset V\oplus V $ as
the image of $ \Pi_{1}\oplus\Pi_{2}\colon W \to V\oplus V $. The pencils $ \Pi $ which may obtained from
linear relations are those for which $ \operatorname{Ker}\Pi_{1}\cap\operatorname{Ker}\Pi_{2}=0 $. For such pencils
these two constructions are mutually inverse.

Introduce {\em another\/}\footnote{In fact, we will not use the previous connection between linear
relations and pencils.} connection between bisurjective linear relations
in $ V $ and pencils of linear operators $ V \to V' $ (here $ V' $ is an arbitrary
vector space). Given a bisurjective linear relation
$ W $ in $ V $, one obtains an identification $ \alpha\colon V/\operatorname{Ker}_{R}W \to V/\operatorname{Ker}_{L}W $. Let
$ V' $ be either $ V/\operatorname{Ker}_{R}W $ or $ V/\operatorname{Ker}_{L}W $ (or any other vector space identified
with both $ V/\operatorname{Ker}_{R}W $ and $ V/\operatorname{Ker}_{L}W $ in such a way that the identifications
commute with $ \alpha $). To simplify notations, assume $ V'=V/\operatorname{Ker}_{L}W $.

Denote by $ \pi_{1,2} $ the natural projections of $ V $ to $ V/\operatorname{Ker}_{L}W $ and to
$ V/\operatorname{Ker}_{R}W $. Then $ {\mathcal P}_{1}\buildrel{\text{def}}\over{=}\pi_{1} $, $ {\mathcal P}_{2}=\alpha\circ\pi_{2} $ give a pair of operators $ V \to V'=V/\operatorname{Ker}_{L}W $.
In the other direction, given a pencil $ {\mathcal P} $ of linear operators $ V \to V' $, one
can consider a linear relation $ W={\mathcal P}_{2}^{-1}{\mathcal P}_{1} $ (in other words, $ \left(v_{1},v_{2}\right)\in W $ iff
$ {\mathcal P}_{1}v_{1}={\mathcal P}_{2}v_{2} $). Obviously, this gives two mutually inverse operations
between bisurjective relations to pencils.

Obviously, a pencil $ {\mathcal P} $ of operators $ V \to V' $ can be obtained from a
bisurjective linear relation in $ V $ iff both linear operators $ {\mathcal P}_{1} $, $ {\mathcal P}_{2} $ of the
pencil are surjective. Call such pencils {\em bisurjective}. In fact, one can
obtain a much stronger result:

\begin{definition} Given a linear relation $ W_{1} $ in $ V_{1} $ and a linear relation $ W_{2} $ in
$ V_{2} $, call a linear mapping $ \varphi\colon V_{1} \to V_{2} $ an {\em isomorphism\/}\footnote{One could consider mappings which are not isomorphisms, but we will not
consider such mappings as morphisms.} between $ W_{1} $ and $ W_{2} $
if $ \varphi $ is an isomorphism, and $ \left(\varphi\oplus\varphi\right)\left(W_{1}\right)=W_{2} $.

Given a pencil $ {\mathcal P} $ of operators $ V \to W $ and a pencil $ {\mathcal P}' $ of operators $ V'
\to W' $ say that linear mappings $ \varphi\colon V \to V' $ and $ \psi\colon W \to W' $ form a {\em morphism\/}
between $ {\mathcal P} $ and $ {\mathcal P}' $ if $ {\mathcal P}'_{1}\varphi=\psi{\mathcal P}_{1} $, $ {\mathcal P}'_{2}\varphi=\psi{\mathcal P}_{2} $. Call $ \left(\varphi,\psi\right) $ an {\em isomorphism\/} if $ \varphi $
and $ \psi $ form a morphism and $ \varphi $ and $ \psi $ are isomorphisms of vector spaces. \end{definition}

\begin{proposition} \label{prop1.50}\myLabel{prop1.50}\relax  Consider a category $ bs{\mathfrak R} $ of bisurjective linear
relations in vector spaces with isomorphisms of relations as $ \operatorname{Mor} bs{\mathfrak R} $, and
the category $ bs{\mathfrak P} $ of bisurjective pencils of linear mappings with
isomorphisms of pencils as $ \operatorname{Mor} bs{\mathfrak P} $. The mappings $ bs{\mathfrak R} \to bs{\mathfrak P} $ and $ bs{\mathfrak P} \to
bs{\mathfrak R} $ defined above give an equivalence of these categories. Moreover, let
$ \operatorname{Vect} $ be the category of vector spaces, consider the functor $ bs{\mathfrak R} \to \operatorname{Vect} $
which sends a relation $ W $ in $ V $ to the vector space $ V $, and the functor $ bs{\mathfrak P}
\to \operatorname{Vect} $ which sends a pencil $ {\mathcal P} $ of operators $ V \to W $ to $ V $. The equivalence
of categories defined above commutes with the mappings to $ \operatorname{Vect} $. \end{proposition}

In plain words, it is ``the same'' to consider bisurjective linear
relations in $ V $ and bisurjective pencils of operators $ V \to V' $ up to
isomorphisms of $ V' $.

\begin{remark} \label{rem1.50}\myLabel{rem1.50}\relax  We will use the equivalence of categories in the following
way: given a bundle $ {\mathcal V} $ over a manifold $ M $ with the fibers being vector
spaces with bisurjective relations, one can construct a vector bundle $ {\mathcal V}' $
the whose fibers are vector spaces $ V' $ of the previous construction. The
vector bundle $ {\mathcal V}' $ is canonically defined. \end{remark}

\begin{remark} \label{rem1.60}\myLabel{rem1.60}\relax  Consider a linear relation $ W $ in a vector space $ V $, i.e.,
$ W\subset V\oplus V $. Note that $ V\oplus V=V\otimes{\mathcal S} $; here $ {\mathcal S} $ is the ``coordinate'' two-dimensional
vector space. The group $ \operatorname{GL}\left({\mathcal S}\right) $ acts on $ V\otimes{\mathcal S} $, thus acts in the set of linear
relations in $ V $. This action can be reduced to the action of $ \operatorname{PGL}\left({\mathcal S}\right) $.

On the other hand, $ \operatorname{GL}\left({\mathcal S}\right) $ also acts in the set of pencils $ {\mathcal P} $ of
operators $ V \to V' $, $ \left({\mathcal P}_{1},{\mathcal P}_{2}\right) \to \left(a{\mathcal P}_{1}+b{\mathcal P}_{2},c{\mathcal P}_{1}+d{\mathcal P}_{2}\right) $. Obviously, these two
actions are compatible with the mappings given above. \end{remark}

\begin{remark} One can make the construction of $ V' $ a little bit more symmetric
by using the definition $ V'=\left(V\oplus V\right)/W $. \end{remark}

\section{Kronecker relations }

Recall that to any bisurjective linear relation $ W $ in $ V $ we associated
a pencil $ {\mathcal P} $ of operators $ V \to V' $ (for an appropriate vector space $ V' $).
Thus one can consider the corresponding linear family $ \lambda_{1}{\mathcal P}_{1}+\lambda_{2}{\mathcal P}_{2} $ of
mappings $ V \to V' $. In particular, for any $ \lambda_{1} $ and $ \lambda_{2} $ one obtains a vector
subspace $ \operatorname{Ker}\left(\lambda_{1}{\mathcal P}_{1}+\lambda_{2}{\mathcal P}_{2}\right)\subset V $. Denote this subspace by $ \operatorname{Ker}_{\lambda_{1}:\lambda_{2}}W $. Obviously,
for $ \left(\lambda_{1},\lambda_{2}\right)\not=\left(0,0\right) $ this subspace depends only on the ratio $ \left(\lambda_{1}:\lambda_{2}\right)\in{\mathbb P}^{1} $, and
$ \operatorname{Ker}_{L}=\operatorname{Ker}{\mathcal P}_{1}=\operatorname{Ker}_{1:0} $, $ \operatorname{Ker}_{R}=\operatorname{Ker}_{0:1} $. Moreover, if the vector spaces we
consider are defined over a field $ {\mathbb K} $, in fact one can consider $ \lambda_{1} $ and $ \lambda_{2} $
to be in any extension $ {\mathbb E} $ of $ {\mathbb K} $, then $ \operatorname{Ker}_{\lambda_{1}:\lambda_{2}}W\subset V\otimes_{{\mathbb K}}{\mathbb E} $.

In particular, if $ \bar{{\mathbb K}} $ is the algebraic closure of $ {\mathbb K} $, then for any
$ \lambda=\left(\lambda_{1}:\lambda_{2}\right)\in\bar{{\mathbb K}}{\mathbb P}^{1} $ one can consider a correctly defined number $ \dim  \operatorname{Ker}_{\lambda}W $.

\begin{definition} A bisurjective linear relation $ W $ in $ V $ is {\em Kronecker\/}
if $ \dim  \operatorname{Ker}_{\lambda}W $ does not depend on $ \lambda\in\bar{{\mathbb K}}{\mathbb P}^{1} $. Call this common dimension the
{\em rank\/} of the relation. \end{definition}

\begin{example} Assume that $ W $ is a graph of a linear mapping $ V \to V $. Then $ W $
is Kronecker iff $ \dim  V=0 $. \end{example}

\begin{example} \label{ex2.30}\myLabel{ex2.30}\relax  Let $ V={\mathbb K}^{n} $. Define $ W\subset V\oplus V $ by $ \left(v,v'\right)\in W $ iff $ v_{k}=v_{k+1}' $,
$ k=1,\dots ,n-1 $. This is a Kronecker linear relation of rank 1. \end{example}

\begin{definition} A {\em Kronecker block\/} is a linear relation isomorphic to
the relation of Example~\ref{ex2.30}. \end{definition}

\begin{definition} Given a linear relation $ W $ in $ V $ and a linear relation $ W' $ in
$ V' $, $ W\oplus W' $ can be considered as a linear relation in $ V\oplus V' $. Call this linear
relation a {\em direct sum\/} of $ W $ and $ W' $. \end{definition}

\begin{definition} For $ \lambda\in{\mathbb K} $ say that a relation $ W $ in $ V $ is a {\em Jordan block\/} with
eigenvalue $ \lambda $ if $ W $ is a graph of a mapping $ V \to V $ which is a Jordan block
with eigenvalue $ \lambda $. Say that $ W $ is a {\em Jordan block\/} with eigenvalue $ \infty $ if $ W^{-1} $
is a Jordan block with eigenvalue 0; here $ W^{-1} $ is the image of $ W $ under
transposition of summands in $ V\oplus V $. \end{definition}

\begin{theorem} \label{th2.50}\myLabel{th2.50}\relax  Suppose that $ {\mathbb K} $ is algebraically closed and $ \dim  V<\infty $.
\begin{enumerate}
\item
any Kronecker linear relation in $ V $ of rank 1 is a Kronecker block.
\item
Any Kronecker linear relation in $ V $ is isomorphic to a direct sum of
Kronecker blocks.
\item
Any bisurjective linear relation is isomorphic to a direct sum of
Kronecker and Jordan blocks.
\end{enumerate}
The collection of dimensions (and---for Jordan blocks---eigenvalues) of
these blocks is uniquely determined by $ V $. \end{theorem}

\begin{proof} See classification of pencils of finite-dimensional linear
operators, say in \cite{Gan59The} or \cite{TurAith61Int}. \end{proof}

Given a Kronecker relation $ W $ in $ V $ of rank $ r $, one obtains a natural
mapping $ \operatorname{Ker}_{\bullet}\colon {\mathbb P}^{1} \to \operatorname{Gr}_{r}V $. Call this parameterized curve in $ \operatorname{Gr}_{r}V $ the
{\em spectral curve\/} of $ W $.

\begin{remark} \label{rem2.70}\myLabel{rem2.70}\relax  Given a bisurjective relation $ W $ in a finite-dimensional
vector space $ V $, $ \dim  \operatorname{Ker}_{\lambda}W $ is constant on a Zariski open subset of $ \bar{{\mathbb K}}{\mathbb P}^{1} $,
i.e., outside of a finite subset $ \Lambda\subset\bar{{\mathbb K}}{\mathbb P}^{1} $. In particular, there is a natural
mapping $ \operatorname{Ker}_{\bullet}\colon {\mathbb P}^{1}\smallsetminus \Lambda \to \operatorname{Gr}_{r}V $. Since $ \operatorname{Gr}_{r}V $ is complete, one can extend this
mapping to a mapping $ {\mathbb P}^{1} \to \operatorname{Gr}_{r}V $, which one can call a {\em spectral curve\/} as
well. \end{remark}

\begin{proposition} \label{prop2.90}\myLabel{prop2.90}\relax  Consider two matrices $ A_{1},A_{2}\in\operatorname{Mat}\left(m,n\right) $. They define a
pencil of mappings $ \bar{{\mathbb K}}^{m} \to \bar{{\mathbb K}}^{n} $. Suppose that $ m\geq n $ and $ \operatorname{rk}\left(\lambda_{1}A_{1}+\lambda_{2}A_{2}\right)=n $ for any
$ \left(\lambda_{1},\lambda_{2}\right)\in\bar{{\mathbb K}}^{2}\smallsetminus\left(0,0\right) $. Then the same is true for any pair $ \left(A_{1}',A_{2}'\right) $ which is
close to $ \left(A_{1},A_{2}\right) $. \end{proposition}

\begin{proof} Consider the projectivization $ {\mathbb P}\operatorname{Mat}\left(m,n\right) $ of the vector space of
matrices. Let $ Z\subset{\mathbb P}\operatorname{Mat}\left(m,n\right) $ be the projection of the set of matrices of
rank $ <n $ to $ {\mathbb P}\operatorname{Mat}\left(m,n\right) $, denote projections of $ A_{1} $, $ A_{2} $ on $ {\mathbb P}\operatorname{Mat}\left(m,n\right) $ by $ \alpha_{1} $,
$ \alpha_{2} $. Then $ Z $ is a closed subset, and the line through $ \alpha_{1} $, $ \alpha_{2} $ does not
intersect $ Z $. Due to compactness of $ {\mathbb P}\operatorname{Mat}\left(m,n\right) $, nearby lines do not
intersect $ Z $ as well. \end{proof}

\begin{remark} \label{rem2.90}\myLabel{rem2.90}\relax  Due to the correspondence of Proposition~\ref{prop1.50}, one
can restate Proposition~\ref{prop2.90} as the fact that Kronecker relations in
a vector space $ V $ form an open subset of all relations in $ V $. Here we
identify the set of relations with the union of Grassmannians of
subspaces of different dimensions in $ V\oplus V $. \end{remark}

\section{Kronecker webs }\label{h34}\myLabel{h34}\relax 

\begin{definition} A {\em preweb\/} on a manifold $ B $ is a bisurjective linear relation
in $ {\mathcal T}^{*}B $, in other words, it is a vector bundle $ {\mathcal W} $ on $ B $ with an inclusion $ {\mathcal W}
\hookrightarrow {\mathcal T}^{*}B\oplus{\mathcal T}^{*}B $ which makes a bisurjective linear relation in each fiber of
$ {\mathcal T}^{*}B $. \end{definition}

Given a preweb $ {\mathcal W} $ and $ \lambda\in{\mathbb K}{\mathbb P}^{1} $, one can consider $ \operatorname{Ker}_{\lambda}{\mathcal W} $, which is a
collection of vector subspaces $ \operatorname{Ker}_{\lambda}{\mathcal W}_{b}\subset{\mathcal T}_{b}^{*}B $, $ b\in B $. Recall that $ {\mathcal N}{\mathcal F} $, $ {\mathcal T}{\mathcal F} $, and
$ {\mathcal T}^{*}{\mathcal F} $ for a foliation $ {\mathcal F} $ were defined in Section~\ref{h002}.

\begin{definition} Consider a vector subbundle $ E $ of $ {\mathcal T}^{*}B $. Call $ E $ {\em integrable\/} if
there is a foliation $ {\mathcal F} $ on $ B $ such that $ E={\mathcal N}{\mathcal F} $. \end{definition}

\begin{definition} \label{def3.09}\myLabel{def3.09}\relax  A preweb $ {\mathcal W} $ is a {\em web\/} if for any given $ \lambda\in{\mathbb K}{\mathbb P}^{1} $
the number $ \dim  \operatorname{Ker}_{\lambda}{\mathcal W}_{b} $ does not depend on $ b\in B $, and this collection of
subspaces is integrable. In other words, there is a foliation $ {\mathcal F}_{\lambda} $ on $ B $
such that $ \operatorname{Ker}_{\lambda}{\mathcal W}_{b} $ coincides with the normal spaces to $ {\mathcal F}_{\lambda} $ at $ b\in B $. Call this
foliation the $ \lambda $-{\em integrating foliation\/} (or just the {\em integrating foliation}, if
$ \lambda $ is clear from the context) of $ {\mathcal W} $.

A preweb $ {\mathcal W} $ is {\em Kronecker\/} of rank $ r $ if for any $ b\in B $ the linear relation
$ {\mathcal W}_{b} $ in $ {\mathcal T}_{b}^{*}B $ is Kronecker of rank\footnote{Since rank of a Kronecker relation $ W $ in $ V $ equals $ \dim  W- \dim  V $, rank of
$ {\mathcal W}_{m} $ does not depend on $ m $.} $ r $. \end{definition}

\begin{remark} \label{rem33.92}\myLabel{rem33.92}\relax  Note that Theorem~\ref{th33.10} implies that the notion of
Kronecker webs is a generalization of the notion of {\em Veronese webs of higher
codimension\/} introduced in \cite{Pan99Ver}. Since the definition of \cite{Pan99Ver}
cannot not stated in a coordinate-independent form, Definition~\ref{def3.09}
is much more convenient to work with. \end{remark}

Due to Proposition~\ref{prop1.50}, to describe a preweb $ {\mathcal W} $ on $ B $ is ``the
same'' as to define a vector bundle $ \Phi\left({\mathcal W}\right) $ on $ B $ and a pencil $ {\mathcal P} $ of
bisurjective mappings of vector bundles $ {\mathcal T}^{*}B \to \Phi\left({\mathcal W}\right) $. If the preweb $ {\mathcal W} $ is
clear from the context, we may denote $ \Phi\left({\mathcal W}\right) $ as $ \Phi $ as well.\footnote{The choice of notation is related to the fact that in applications the
coordinates on fibers of $ \Phi $ are angle-coordinates of an integrable system.
(The {\em action\/} variables are coordinates on the global space of $ \Phi $ coming from
coordinates on $ M $.)} Recall that
given $ {\mathcal W} $, $ \Phi\left({\mathcal W}\right)={\mathcal T}^{*}B/\operatorname{Ker}_{\lambda}{\mathcal W} $ (here $ \lambda $ is any number fixed in advance), given $ \Phi $
and $ {\mathcal P} $, the vector bundle $ {\mathcal W}\subset{\mathcal T}^{*}B\oplus{\mathcal T}^{*}B $ can be described by the condition
$ {\mathcal P}_{1,b}v_{1}={\mathcal P}_{2,b}v_{2} $, $ \left(v_{1},v_{2}\right)\in{\mathcal T}_{b}^{*}B $.

By definition, this preweb is Kronecker if $ \operatorname{Ker}\left(\lambda_{1}{\mathcal P}_{1,b}+\lambda_{2}{\mathcal P}_{2,b}\right) $, $ b\in B $,
is an family of subspaces of $ {\mathcal T}^{*}B $ (or $ {\mathcal T}^{*}B\otimes{\mathbb C} $) of the same dimension for any
$ b\in B $ and $ \left(\lambda_{1}:\lambda_{2}\right)\in{\mathbb C}{\mathbb P}^{1} $. If $ {\mathcal W} $ is a Kronecker preweb, $ {\mathcal W} $ is a Kronecker web if
$ \operatorname{Ker}\left(\lambda_{1}{\mathcal P}_{1,b}+\lambda_{2}{\mathcal P}_{2,b}\right) $ is an integrable family of subspaces of $ {\mathcal T}^{*}B $ of the
same dimension for any $ b\in B $ and $ \left(\lambda_{1}:\lambda_{2}\right)\in{\mathbb K}{\mathbb P}^{1} $. Due to Remark~\ref{rem2.90}, the
condition on a (pre)web to be Kronecker is a condition of being in
general position.

For a given $ b\in B $ one can define $ \Lambda_{b}\subset{\mathbb C}{\mathbb P}^{1} $ as in Remark~\ref{rem2.70}.
Obviously, if $ \lambda_{0}\notin\Lambda_{b_{0}} $, then there are neighborhoods $ U\ni b_{0} $ and $ {\mathbit U}\ni\lambda_{0} $ such
that the vector subspaces $ \operatorname{Ker}_{\lambda}{\mathcal W}_{b} $ depend smoothly on $ b\in U $ and $ \lambda\in{\mathbit U} $. For any
given $ \left(\lambda_{1}:\lambda_{2}\right)\in{\mathbit U}\subset{\mathbb C}{\mathbb P}^{1} $ the bundle $ {\mathcal T}^{*}B/\operatorname{Ker}_{\left(\lambda_{1}:\lambda_{2}\right)}{\mathcal W}|_{U} $ is canonically
isomorphic to $ \Phi|_{U} $. On the other hand, if $ \lambda\in{\mathbb K}{\mathbb P}^{1} $, and $ {\mathcal W} $
is a web, then $ {\mathcal T}^{*}B/\operatorname{Ker}_{\left(\lambda_{1}:\lambda_{2}\right)}{\mathcal W}={\mathcal T}^{*}{\mathcal F}_{\left(\lambda_{1}:\lambda_{2}\right)} $.

We see that given $ \left(\lambda_{1},\lambda_{2}\right)\in{\mathbb K}^{2}\smallsetminus\left(0,0\right) $ such that $ \lambda=\left(\lambda_{1}:\lambda_{2}\right)\in{\mathbit U} $ the vector
bundle $ \Phi|_{U} $ is canonically identified with $ {\mathcal T}^{*}{\mathcal F}_{\lambda}|_{U} $. Since $ {\mathcal T}^{*}{\mathcal F} $ has a natural
Poisson structure (see Section~\ref{h002}), we see that

\begin{proposition} The total space $ \Phi|_{U} $ carries a natural Poisson structure
for any $ \left(\lambda_{1},\lambda_{2}\right)\in{\mathbb K}^{2}\smallsetminus\left(0,0\right) $ such that $ \lambda=\left(\lambda_{1}:\lambda_{2}\right)\in{\mathbit U} $. This Poisson structure
depends smoothly on $ \left(\lambda_{1},\lambda_{2}\right) $. \end{proposition}

Call this Poisson structure $ \eta_{\lambda_{1},\lambda_{2}} $. Let $ \widetilde{U}=\left\{\left(\lambda_{1},\lambda_{2}\right) \mid \lambda_{1}:\lambda_{2}\in U\right\} $.

\begin{proposition} The Poisson structure $ \eta_{\lambda_{1},\lambda_{2}} $ is homogeneous in
$ \left(\lambda_{1},\lambda_{2}\right)\in\widetilde{U} $ of degree 1. \end{proposition}

\begin{proof} Indeed, multiplication of $ \left(\lambda_{1},\lambda_{2}\right) $ by a constant changes the
identification of $ \Phi|_{U} $ with $ {\mathcal T}^{*}{\mathcal F}_{\lambda}|_{U} $ by the same constant. For a
diffeomorphism $ \alpha\colon M_{1} \to M_{2} $ and a Poisson structure $ \eta $ on $ M_{1} $ denote by $ \alpha_{*}\eta $
the Poisson structure $ \eta $ transferred to $ M_{2} $ via $ \alpha $,
$ \alpha_{*}\eta\left(f,g\right)\buildrel{\text{def}}\over{=}\eta\left(f\circ\alpha,g\circ\alpha\right)\circ\alpha^{-1} $. Apply this
definition in the case when $ \alpha=\mu_{c} $ is multiplication by $ c $ in $ {\mathcal T}^{*}N $.

Now the only thing to prove is that $ \left(\mu_{c}\right)_{*}\eta=c\eta $ if $ \eta $ is the canonical
Poisson structure on $ {\mathcal T}^{*}N $, and $ N $ is an arbitrary manifold. One can check
it momentarily in local coordinates on $ N $ (as in Example~\ref{ex002.41}). \end{proof}

\begin{corollary} \label{cor3.50}\myLabel{cor3.50}\relax  If $ {\mathbb K}={\mathbb C} $, and the web $ {\mathcal W} $ is Kronecker, then the Poisson
structure $ \eta_{\lambda_{1},\lambda_{2}} $ on the total space of $ \Phi $ depends linearly on $ \lambda_{1},\lambda_{2} $. \end{corollary}

\begin{proof} Since for a Kronecker web $ \Lambda_{b}=\varnothing $, one can take $ {\mathbit U}={\mathbb C}{\mathbb P}^{1} $, and $ U=B $,
thus the Poisson structure $ \eta_{\lambda_{1},\lambda_{2}} $ is defined on the whole total space of
$ \Phi $ for any $ \left(\lambda_{1},\lambda_{2}\right)\not=\left(0,0\right) $. Fix $ x\in\Phi $. Since a Poisson structure is a bivector
field, one can associate to $ \left(\lambda_{1},\lambda_{2}\right)\not=\left(0,0\right) $ an element $ \eta_{\lambda_{1},\lambda_{2}}|_{x} $ of $ \Lambda^{2}{\mathcal T}_{x}\Phi $.
We know that this element depends analytically on $ \left(\lambda_{1},\lambda_{2}\right)\not=\left(0,0\right) $ and is
homogeneous of degree 1 in $ \left(\lambda_{1},\lambda_{2}\right) $. However, any mapping of $ {\mathbb C}^{2}\smallsetminus\left(0,0\right) $ to a
vector space which is of homogeneity degree 1 is linear, which finishes
the proof. \end{proof}

\begin{theorem} \label{th3.60}\myLabel{th3.60}\relax  If a web $ {\mathcal W} $ is Kronecker, then the Poisson structure
$ \eta_{\lambda_{1},\lambda_{2}} $ on the total space of $ \Phi $ depends linearly on $ \lambda_{1},\lambda_{2} $. \end{theorem}

\begin{proof} The only case which remains to be proven is $ {\mathbb K}={\mathbb R} $. If a Kronecker
web $ {\mathcal W} $ is real-analytic, then one can consider the complex-analytic
continuation, and one returns to the case $ {\mathbb K}={\mathbb C} $. Thus in real-analytic case
Corollary~\ref{cor3.50} implies the theorem.

Consider now $ C^{\infty} $-case. First, note that the theorem follows from
some ``abstract nonsense'' remarks.
Recall that one of the contributions of algebraic geometry to
differential geometry is the understanding of the importance of
considering ``formal'' objects as tools for investigation of ``geometric''
objects. In particular, given a $ C^{\infty} $-manifold, one can consider ``an $ \infty $-jet
of the complex-analytic neighborhood'' (or a ``formal neighborhood'') of
this manifold. This formal neighborhood is canonically defined, and
carries many properties of complex-analytic manifolds.

Readers familiar with the language we used above may immediately
recognize that all the objects needed for the proof of Corollary~\ref{cor3.50}
(foliations, tensors, cotangent bundles, Poisson structures) make sense
in settings of ``formal geometry'', thus one can finish the proof in
$ C^{\infty} $-case in the same way we did it in real-analytic case (when we had a
no-nonsense complex neighborhood instead of a formal one). For other
readers we provide a stripped-down version of the proof below. This
version will not use ``hard'' notions of formal geometry (such as jets of
manifolds etc.), but will use only notions of (finite-order) jets of
sections of ``real'' bundles on ``real'' manifolds.\footnote{As opposed to ``formal'' manifolds.} This simplified version
of the proof goes until the end of this section.

Consider a vector bundle $ E $ over $ B $, fix a point $ b\in B $. To describe a
vector subbundle $ F $ of $ E $ of rank $ d $ is the same as to describe a section of
the bundle $ \operatorname{Gr}_{d}\left(E\right) $ (the fibers of this bundle are Grassmannians of fibers
of $ F $). Recall that a $ k $-jet near $ b $ of a section $ f $ of a bundle is the
collection of Taylor coefficients of $ f $ of the order $ k $ or less.\footnote{In a coordinate-less form, $ k $-jets are equivalence classes of sections;
here the equivalence is having the same Taylor coefficients of order $ k $ or
less in any coordinate system.} If a
geometric object can be described locally by a section of some bundle,
one can define a $ k $-jet of this object near $ b $ as being a $ k $-jet of a
section of this bundle\footnote{However, it may happens that there are different descriptions of the
same object which lead to mutually shifted enumerations of jets. For
example, locally one can describe a closed $ 1 $-form $ \alpha $ by a function $ f $ such
that $ df=\alpha $, or by components of $ \alpha $ considered as a tensor field. The $ k $-jet
in the first description is a $ k-1 $-jet in the second description.}. Thus one can define a $ k $-jet of a vector
subbundle $ F $ of $ E $ near $ b $.

\begin{remark} In what follows we will also need to consider $ E/F $; here $ F $ is a
($ k $-jet of) a subbundle of $ E $. To avoid introducing some new ``abstract
nonsense'', equip $ E $ with a norm, and identify $ E/F $ with the orthogonal
complement to $ F $. Obviously, given a $ k $-jet of $ F $, the $ k $-jet of $ E/F $ is
well-defined as a $ k $-jet of a vector (sub)bundle. \end{remark}

In particular, one can define a $ k $-jet $ {\mathcal W} $ of a preweb. Consider two
prewebs $ {\mathcal W} $ and $ {\mathcal W}' $ which have the same $ k $-jet near $ b $. Then the subbundles
$ \operatorname{Ker}_{\lambda}{\mathcal W}\subset{\mathcal T}^{*}B $ and $ \operatorname{Ker}_{\lambda}{\mathcal W}'\subset{\mathcal T}^{*}B $ have the same $ k $-jet near $ b $. Thus there is a
canonically defined $ k $-jet of an isomorphism $ \iota $ between vector bundles
$ {\mathcal T}^{*}B/\operatorname{Ker}_{\lambda}{\mathcal W} $ and $ {\mathcal T}^{*}B/\operatorname{Ker}_{\lambda}{\mathcal W}' $. This leads to a $ k $-jet\footnote{It is possible to define what is a $ k $-jet of something {\em near a
submanifold}, and the $ k $-jet in question is defined near two fibers over $ b $ of
these two vector bundles.} of an isomorphism
between the total spaces of these bundles. On the other hand, given a
$ k $-jet of an isomorphism $ \iota $ between two manifolds $ F $ and $ F' $, and a tensor
field $ \tau $ on $ F $, one can define a $ k-1 $-jet of the tensor field $ \iota_{*}{\mathcal T} $ on $ F' $.

Suppose that both prewebs $ {\mathcal W} $ and $ {\mathcal W}' $ are in fact webs. As we know,
total spaces of both bundles $ {\mathcal T}^{*}B/\operatorname{Ker}_{\lambda}{\mathcal W} $ and $ {\mathcal T}^{*}B/\operatorname{Ker}_{\lambda}{\mathcal W}' $ carry natural
Poisson structures, and one can consider Poisson structures as bivector
(thus tensor) fields $ \eta $ and $ \eta' $. Thus given a $ k $-jet $ \iota $ of an isomorphism
between these manifolds, one can consider $ \iota_{*}\eta $, which is a $ k-1 $-jet of a
Poisson structure on $ {\mathcal T}^{*}B/\operatorname{Ker}_{\lambda}{\mathcal W}' $. Obviously, $ \iota_{*}\eta $ coincides with the
$ k-1 $-jet of $ \eta' $. This immediately implies

\begin{lemma} \label{lm3.100}\myLabel{lm3.100}\relax  Consider two Kronecker webs $ {\mathcal W} $ and $ {\mathcal W}' $ on $ B $ which have the
same $ k $-jet near $ b\in B $. Then there is a canonically defined $ k $-jet of an
isomorphism $ \iota $ between $ \Phi\left({\mathcal W}\right) $ and $ \Phi\left({\mathcal W}'\right) $, and this isomorphism identifies
$ k-1 $-jets of Poisson structures $ \eta_{\lambda_{1},\lambda_{2}} $ (on $ \Phi $) and $ \eta'_{\lambda_{1},\lambda_{2}} $ (on $ \Phi' $) near
fibers of these bundles over $ b $. \end{lemma}

Thus

\begin{lemma} In the conditions of Lemma~\ref{lm3.100} suppose that $ k=1 $, and the
web $ {\mathcal W}' $ is in fact real-analytic. Then the tensor field $ \eta_{\lambda_{1},\lambda_{2}} $ on the
fiber of $ \Phi $ over $ b $ depends linearly on $ \left(\lambda_{1},\lambda_{2}\right) $. \end{lemma}

In particular, the theorem follows from the case $ k=1 $ of

\begin{conjecture} \label{con3.140}\myLabel{con3.140}\relax  Fix $ k\geq0 $. Given a Kronecker web $ {\mathcal W} $ on an open subset
$ U\subset{\mathbb R}^{n} $, $ 0\in U $, one can find a real-analytic Kronecker web $ {\mathcal W}' $ on $ U'\subset U $, $ 0\in U' $,
such that $ k $-jets of $ {\mathcal W} $ and $ {\mathcal W}' $ at 0 coincide. \end{conjecture}

In fact, in the case $ k=1 $ we propose the following amplification:

\begin{conjecture} Given a Kronecker web $ {\mathcal W} $ on an open subset $ U\subset{\mathbb R}^{n} $, $ 0\in U $, one
can find $ U'\subset U $, $ 0\in U' $, and a diffeomorphism $ U' \xrightarrow[]{f} V\subset{\mathbb R}^{n} $, $ f\left(0\right)=0 $, such that
$ 1 $-jet of $ f_{*}{\mathcal W} $ (the transfer of $ {\mathcal W} $ via $ f $) coincides with $ 1 $-jet of an
appropriate translation-invariant\footnote{A Kronecker web $ {\mathcal W} $ on a vector space $ V $ is {\em translation-invariant\/} if the
vector subspace $ {\mathcal W}_{v}\subset{\mathcal T}_{v}^{*}V\oplus{\mathcal T}_{v}^{*}V=V^{*}\oplus V^{*} $ does not depend on $ v\in V $.} Kronecker web $ {\mathcal W}' $. \end{conjecture}

Not only this conjecture implies the case $ k=1 $ of Conjecture
~\ref{con3.140}, but it would also directly imply the theorem we are proving.
However, since these conjectures are not settled yet\footnote{As Proposition~\ref{prop3.71} implies, the condition on a preweb of being a web
can be written as a system of differential equations in partial
derivatives. These systems are non-linear, and in all but a handful of
cases they are enormously overdetermined.
\endgraf
This makes the conjectures above so hard to tackle.}, one needs an
alternative way to prove the theorem.

We are going to do this without the extension-properties for webs
which are mentioned above by introducing integrability conditions on
$ k $-{\em jets of prewebs}. Since jets are just collections of numbers, thus are
not sensible to a change of the base field, we will be able to consider a
$ k $-jet of a $ C^{\infty} $-web as a $ k $-jet of a complex-analytic web, thus we will be
able to extend $ \eta_{\lambda_{1},\lambda_{2}} $ to $ \left(\lambda_{1},\lambda_{2}\right)\in{\mathbb C}^{2}\smallsetminus\left(0,0\right) $, then apply the arguments of
Corollary~\ref{cor3.50}.

\begin{definition} A $ k $-jet $ K $ of a vector subbundle of $ {\mathcal T}^{*}B $ is $ k $-{\em integrable\/} near
$ b_{0}\in B $ if there is a foliation $ {\mathcal F} $ such that the $ k $-jet of the normal bundle
to $ {\mathcal F} $ coincides with $ K $. \end{definition}

\begin{definition} A $ k $-jet $ {\mathcal W} $ of a preweb is a $ k $-{\em jet-web\/} if for any $ \lambda\in{\mathbb K}{\mathbb P}^{1} $ the
$ k $-jet $ \operatorname{Ker}_{\lambda}{\mathcal W} $ of a vector subbundle of $ {\mathcal T}^{*}B $ is $ k $-integrable. A $ k $-jet $ {\mathcal W} $ of a
web near $ b $ is {\em Kronecker\/} if the relation $ {\mathcal W}_{b} $ in $ {\mathcal T}_{b}^{*}B $ is Kronecker. \end{definition}

Due to Remark~\ref{rem2.90}, if a $ k $-jet near $ b $ of preweb $ {\mathcal W} $ on $ B $ is
Kronecker, then the preweb $ {\mathcal W} $ is Kronecker in a neighborhood of $ b $. Thus
the definition of a Kronecker $ k $-jet-web above is compatible with taking a
$ k $-jet of a Kronecker web.

Consider a $ k $-jet-web $ {\mathcal W} $. It is a $ k $-jet of a preweb, denote this
preweb by $ {\mathcal W}' $. Given $ \left(\lambda_{1},\lambda_{2}\right)\in{\mathbb K}^{2}\smallsetminus\left(0,0\right) $, the $ k $-jet of $ {\mathcal T}^{*}B/\operatorname{Ker}_{\lambda}{\mathcal W}' $, $ \lambda=\lambda_{1}:\lambda_{2} $,
is naturally identified with the $ k $-jet of the cotangent bundle of the
corresponding foliation $ {\mathcal F}_{\lambda} $, in other words, one can construct a $ k $-jet of
an identification $ \alpha $ of $ {\mathcal T}^{*}B/\operatorname{Ker}_{\lambda}{\mathcal W}' $ with $ {\mathcal T}^{*}{\mathcal F}_{\lambda} $. Recall that a $ k $-jet of a
diffeomorphism acts on $ k-1 $-jets of tensor fields. $ {\mathcal T}^{*}{\mathcal F}_{\lambda} $ carries a natural
Poisson structure, thus it carries the tensor field which describes the
bracket of the Poisson structure. Using the identification $ \alpha $ mentioned
above, one obtains a $ k-1 $-jet of a tensor field $ \eta_{\lambda_{1},\lambda_{2}} $ on $ {\mathcal T}^{*}B/\operatorname{Ker}_{\lambda}{\mathcal W}' $.
Obviously, one can consider this $ k-1 $-jet as living on $ {\mathcal T}^{*}B/\operatorname{Ker}_{\lambda}{\mathcal W} $. It does
not depend on the choice of $ {\mathcal W}' $.

If $ {\mathbb K}={\mathbb C} $, and $ k=1 $, then one obtains a $ 0 $-jet\footnote{Which is a tensor field on the total space of $ \Phi $ defined on the fiber
over $ b $ only.} of a tensor field $ \eta_{\lambda_{1},\lambda_{2}} $
for any $ \left(\lambda_{1},\lambda_{2}\right)\in{\mathbb C}^{2}\smallsetminus\left(0,0\right) $, and the same arguments as in the proof of
Corollary~\ref{cor3.50} show that

\begin{proposition} Let $ {\mathbb K}={\mathbb C} $, consider a Kronecker $ 1 $-jet-web $ {\mathcal W} $ near $ b_{0}\in B $. The
family of tensors fields $ \eta_{\lambda_{1},\lambda_{2}} $ on the total space of $ \Phi $ is well defined
on the fiber over $ b_{0} $ of the projection $ \Phi \to B $. This family depends
linearly on $ \lambda_{1} $, $ \lambda_{2} $. \end{proposition}

\begin{proof} The only thing to prove is that $ \eta_{\lambda_{1},\lambda_{2}} $ depends smoothly on
$ \lambda_{1},\lambda_{2} $. The definition of $ \eta_{\lambda_{1},\lambda_{2}} $ depended on the foliation $ {\mathcal F}_{\lambda} $, and a
priori we have no conditions on how $ {\mathcal F}_{\lambda} $ depends on $ \lambda $. It is enough to show
that the $ k $-jet of $ {\mathcal F}_{\lambda} $ depends smoothly on $ \lambda $. Any foliation $ {\mathcal F} $ on a manifold
$ M $ can be locally described by equations $ Y=F\left(X,T\right) $; here
$ \left(x_{1},\dots ,x_{k},y_{1},\dots ,y_{r}\right) $ is a local coordinate system, $ X=\left(x_{1},\dots ,x_{k}\right) $,
$ Y=\left(y_{1},\dots ,y_{r}\right) $, $ T=\left(t_{1},\dots ,t_{r}\right) $, and $ F\left(0,T\right)\equiv T $. If we know the subbundle
$ {\mathcal N}{\mathcal F}\subset{\mathcal T}^{*}M $, then we know that
\begin{equation}
\frac{\partial F}{\partial X}=A\left(x,F\right),
\label{equ3.150}\end{equation}\myLabel{equ3.150,}\relax 
here $ A $ is a matrix function of $ k+r $ variables which can be deduced from
equations of $ {\mathcal N}{\mathcal F} $ in $ {\mathcal T}^{*}M $. Given $ A $, one can find the unique solution of a
{\em part\/} of~\eqref{equ3.150} with initial conditions $ F\left(0,T\right)\equiv T $: Require that $ \partial F/\partial x_{l} $
is given by~\eqref{equ3.150} if $ x_{l+1}=\dots =x_{k}=0 $. Moreover, given a $ k $-jet of $ A $,
this uniquely determines $ k+1 $-jet of $ F $.

This implies that if there is a family of foliations $ {\mathcal F}^{\left(\mu\right)} $ and the
subbundle $ {\mathcal N}{\mathcal F}^{\left(\mu\right)}\subset{\mathcal T}^{*}M $ depends smoothly on $ \mu $, then $ {\mathcal F}^{\left(\mu\right)} $ depends smoothly on
$ \mu $. In turn, this implies that $ \eta_{\lambda_{1},\lambda_{2}} $ depends smoothly on $ \lambda_{1},\lambda_{2} $. \end{proof}

To finish the proof of the theorem, it is enough to show that given
a Kronecker $ 1 $-jet-web over $ {\mathbb R} $, one can consider it as a Kronecker
$ 1 $-jet-web over $ {\mathbb C} $. Recall that a jet of a preweb is just a collection of
Taylor coefficients of a section of a bundle, thus a collection of
numbers. Any collection of real numbers can be considered as a collection
of complex numbers, thus a $ 1 $-jet of a preweb over $ {\mathbb R} $ can be considered as
$ 1 $-jet of a preweb over $ {\mathbb C} $. Since the condition of being Kronecker does not
change when we change field of scalars, the only thing we need to prove
is the integrability condition for complex $ \lambda_{1}:\lambda_{2} $.

Recall Frobenius integrability condition:

\begin{proposition} \label{prop3.71}\myLabel{prop3.71}\relax  Consider a vector subbundle $ E $ of $ {\mathcal T}^{*}B $. Suppose that
for any small open subset $ U\subset B $ and any section $ \alpha $ of $ E $ over $ U $ its de Rham
differential $ d\alpha $ (which is a $ 2 $-form on $ B $) can be written as $ d\alpha=\sum\alpha_{i}\wedge\beta_{i} $
with $ \alpha_{i} $ being sections of $ E $ and $ \beta_{i} $ being arbitrary differential forms on
$ B $. Then $ E $ is integrable. \end{proposition}

By constructing jet-solutions of ordinary differential equations,
one can easily prove the following jet-analogue of Proposition~\ref{prop3.71}:

\begin{proposition} \label{prop3.80}\myLabel{prop3.80}\relax  Consider a $ k $-jet of a vector subbundle $ E $ of $ {\mathcal T}^{*}B $ near
$ b\in B $. Suppose that for any open subset $ U\subset B $ and for any $ k $-jet $ \alpha $ of a
section of $ E $ near $ b\in B $ its de Rham differential $ d\alpha $ (which is a $ k-1 $-jet of
a $ 2 $-form on $ B $) can be written as $ \sum\alpha_{i}\wedge\beta_{i} $ with $ \alpha_{i} $ being $ k-1 $-jets of
sections of $ E $ and $ \beta_{i} $ being arbitrary differential forms on $ B $. Then $ E $ is
$ k $-integrable. \end{proposition}

In fact one can do more. Given a section $ \alpha $ of $ E $, consider the image
of $ d\alpha $ under projection $ \Omega^{2}B=\Lambda^{2}{\mathcal T}^{*}B \to \Lambda^{2}\left({\mathcal T}^{*}B/E\right) $. This defines a mapping $ \delta $
from sections of $ E $ to sections of $ \Lambda^{2}\left({\mathcal T}^{*}B/E\right) $. A priori it is a
differential operator of order 1, but in fact it has order 0, thus is a
linear mapping between bundles. Indeed, if $ \alpha|_{b}=0 $, then one can easily
check that $ \delta\alpha|_{b}=0 $.

Thus $ \delta $ defines a section of the vector bundle $ E^{*}\otimes\Lambda^{2}\left({\mathcal T}^{*}B/E\right) $, call
this section the {\em torsion\/} of $ E $. Given a $ k $-jet of a vector subbundle $ E\subset{\mathcal T}^{*}B $, $ \delta $
is defined as $ k-1 $-jet of a section of $ E^{*}\otimes\Lambda^{2}\left({\mathcal T}^{*}B/E\right) $. If a vector subbundle
$ E_{t} $ depends smoothly on a parameter $ t\in T $, then $ \delta_{t} $ depends smoothly on $ t $
(i.e., it is a smooth section of the vector bundle $ E_{t}^{*}\otimes\Lambda^{2}\left({\mathcal T}^{*}B/E_{t}\right) $ over $ B\times T $).

Suppose $ k\geq1 $. Apply the construction of $ \delta $ above to the vector
subbundle $ \operatorname{Ker}_{\lambda}{\mathcal W}\subset{\mathcal T}^{*}B $ of a Kronecker preweb $ {\mathcal W} $ in the case $ {\mathbb K}={\mathbb C} $. Here $ \lambda\in{\mathbb C}{\mathbb P}^{1} $,
thus $ \delta $ is a section of $ \left(\operatorname{Ker}_{\lambda}{\mathcal W}\right)^{*}\otimes\Lambda^{2}\left({\mathcal T}^{*}B/\operatorname{Ker}_{\lambda}{\mathcal W}\right) $ over $ B\times{\mathbb C}{\mathbb P}^{1} $. Restrict this
section to $ \left\{b\right\}\times{\mathbb C}{\mathbb P}^{1} $, $ b\in B $. We obtain a vector bundle over $ {\mathbb C}{\mathbb P}^{1} $, and a
regular section $ \delta^{\left(b\right)} $ of this vector bundle. There may be only two
different cases: either $ \delta^{\left(b\right)}\equiv 0 $, or $ \delta^{\left(b\right)} $ vanishes at a finite number of
points of $ {\mathbb C}{\mathbb P}^{1} $.

\begin{proposition} Consider a Kronecker $ 1 $-jet-web $ {\mathcal W} $ near $ 0\in{\mathbb R}^{n} $ (and $ {\mathbb K}={\mathbb R} $). The
complexification of this $ 1 $-jet-web is a $ 1 $-jet $ {\mathcal W}_{{\mathbb C}} $ of a Kronecker preweb
with $ {\mathbb K}={\mathbb C} $. Then $ {\mathcal W}_{{\mathbb C}} $ is a Kronecker $ 1 $-jet-web. \end{proposition}

\begin{proof} Arbitrarily extend $ {\mathcal W} $ to a preweb $ {\mathcal W}' $ in neighborhood of 0 in $ {\mathbb R}^{n} $,
and arbitrarily extend $ {\mathcal W}_{{\mathbb C}} $ to a preweb $ {\mathcal W}_{{\mathbb C}}' $ in a neighborhood of 0 in $ {\mathbb C}^{n} $.
Consider the torsion $ \delta_{{\mathbb C}} $ of the subbundles $ \operatorname{Ker}_{\lambda}{\mathcal W}_{{\mathbb C}}' $ restricted to $ \left\{0\right\}\times{\mathbb C}{\mathbb P}^{1} $,
and the torsion $ \delta $ of the subbundles $ \operatorname{Ker}_{\lambda}{\mathcal W}' $ restricted to $ \left\{0\right\}\in{\mathbb R}{\mathbb P}^{1} $.
Obviously, $ \delta $ is a restriction of $ \delta_{{\mathbb C}} $ from $ {\mathbb C}{\mathbb P}^{1} $ to $ {\mathbb R}{\mathbb P}^{1} $. On the other hand, $ \delta $
vanishes, since $ 1 $-jet $ {\mathcal W} $ of $ {\mathcal W}' $ is a $ 1 $-jet-web. Thus $ \delta_{{\mathbb C}} $ vanishes on $ {\mathbb R}{\mathbb P}^{1} $,
thus at infinitely many points of $ {\mathbb C}{\mathbb P}^{1} $, thus $ \delta_{{\mathbb C}}=0 $. This implies that the
conditions of Proposition~\ref{prop3.80} are satisfied, which finishes the
proof. \end{proof}

This finishes the proof of Theorem~\ref{th3.60}. \end{proof}

In the proof above we ignored the question of possibility of
extension of jets of webs completely. Let us state

\begin{conjecture} Given a Kronecker $ k $-jet-web $ {\mathcal W} $ near $ 0\in{\mathbb C}^{n} $, one can find
$ U\subset{\mathbb C}^{n} $, $ 0\in U $, and a Kronecker web $ {\mathcal W}' $ on $ U $ such that $ k $-jet of $ {\mathcal W}' $ at 0
coincides with $ {\mathcal W} $. \end{conjecture}

Anyway, we proved the following

\begin{corollary} \label{cor33.99}\myLabel{cor33.99}\relax  Given a Kronecker web $ {\mathcal W} $, the total space of the
vector bundle $ \Phi\left({\mathcal W}\right) $ is equipped with a natural bihamiltonian structure. \end{corollary}

\section{Pairs of skew-symmetric forms }\label{h37}\myLabel{h37}\relax 

Recall the classification of pairs of skew-symmetric bilinear
pairings from \cite{GelZakhFAN} (see also \cite{GelZakhWeb,GelZakh93}). For $ k\in{\mathbb N} $
consider the identity $ k\times k $ matrix $ I_{k} $. For $ \mu\in{\mathbb C} $ consider the Jordan block
$ J_{k,\mu} $ of size $ k $ and eigenvalue $ \mu $. The pair of matrices
\begin{equation}
{\text H}_{1}^{\left(\mu\right)}= \left( 
\begin{matrix}
0 & J_{k,\mu}
\\
-J_{k,\mu}^{t} & 0
\end{matrix}
\right),\qquad {\text H}_{2}^{\left(\mu\right)}=\left( 
\begin{matrix}
0 & I_{k}
\\
-I_{k} & 0
\end{matrix}
\right)
\notag\end{equation}
defines a pair of skew-symmetric bilinear pairings on vector space $ {\mathbb C}^{2k} $. The
limit case of $ \mu \to \infty $ may be deformed to
\begin{equation}
{\text H}_{1}^{\left(\infty\right)}= \left( 
\begin{matrix}
0 & I_{k}
\\
-I_{k} & 0
\end{matrix}
\right),\qquad {\text H}_{2}^{\left(\infty\right)}= \left( 
\begin{matrix}
0 & J_{k,0}
\\
-J_{k,0}^{t} & 0
\end{matrix}
\right).
\notag\end{equation}
Denote the pair $ \left({\text H}_{1}^{\left(\mu\right)},{\text H}_{2}^{\left(\mu\right)}\right) $ of skew-symmetric bilinear pairings by
$ {\mathcal J}_{2k,\mu} $, $ k\in{\mathbb N} $, $ \mu\in{\mathbb C}{\mathbb P}^{1} $.

Add to this list the so-called Kroneker pair $ {\mathcal K}_{2k-1} $. This is a pair
in a vector space $ {\mathbb C}^{2k-1} $ with a basis $ \left({\mathbit w}_{0},{\mathbit w}_{1},\dots ,{\mathbit w}_{2k-2}\right) $. The only non-zero
pairings are
\begin{equation}
\left({\mathbit w}_{2l},{\mathbit w}_{2l+1}\right)_{1}=1,\qquad \left({\mathbit w}_{2l+1},{\mathbit w}_{2l+2}\right)_{2}=1,
\label{equ2.10}\end{equation}\myLabel{equ2.10,}\relax 
for $ 0\leq l\leq k-2 $. Obviously, different pairs from this list are not isomorphic.

\begin{theorem} \label{th4.35}\myLabel{th4.35}\relax  (\cite{GelZakhFAN,Thom91Pen}) Any pair of skew-symmetric
bilinear pairings on a finite-dimensional complex vector space can be
decomposed into a direct sum of pairs of the pairings isomorphic to
$ {\mathcal J}_{2k,\mu} $, $ k\in{\mathbb N} $, $ \mu\in{\mathbb P}^{1} $, and $ {\mathcal K}_{2k-1} $, $ k\in{\mathbb N} $. The types of the components of this
decomposition are uniquely determined. \end{theorem}

Though this simple statement was known for a long time (for example,
the preprint of \cite{Thom91Pen} existed in 1973), we do not know whether it
was published before it was used in \cite{GelZakhFAN}. The discussions in
\cite{Gan59The} and \cite{TurAith61Int} come very close, but do not state this
result.

\begin{definition} A pair of bilinear skew-symmetric forms $ \left(,\right)_{1} $ and
$ \left(,\right)_{2} $ in a vector space $ V $ over $ {\mathbb C} $ is {\em Kronecker\/} if it has no Jordan blocks.
The {\em rank\/} of a Kronecker pair is the number of Kronecker blocks in the
decomposition of the pair into indecomposable components.\footnote{Note a conflict in these notations: for a pair $ \left(,\right)_{1,2} $ of {\em rank\/} $ r $, both
forms have {\em corank\/} $ k $.}

Given a Kronecker pair of bilinear skew-symmetric forms in $ V $ let the
{\em action subspace\/}\footnote{The reason for the name is that in applications this subspace is spanned
by differentials of {\em action function\/} of a system of action-angle
variables.} of $ V $ be spanned by vector subspaces $ \operatorname{Ker}\left(,\right)_{\lambda_{1},\lambda_{2}} $,
$ \left(\lambda_{1},\lambda_{2}\right)\in{\mathbb C}^{2}\smallsetminus\left(0,0\right) $. Here $ \left(v,v'\right)_{\lambda_{1},\lambda_{2}}\buildrel{\text{def}}\over{=}\lambda_{1}\left(v,v'\right)_{1}+\lambda_{2}\left(v,v'\right)_{2} $. \end{definition}

\begin{proposition} \label{prop4.40}\myLabel{prop4.40}\relax  The action subspace of a Kronecker pair of bilinear
forms in $ V $ of rank $ r $ is isotropic with respect to $ \left(,\right)_{1} $ and $ \left(,\right)_{2} $, and has
dimension $ \frac{\dim  V+r}{2} $. It is a maximal isotropic subspace for any form
$ \left(,\right)_{\lambda_{1},\lambda_{2}} $, $ \left(\lambda_{1},\lambda_{2}\right)\in{\mathbb C}^{2}\smallsetminus\left(0,0\right) $. \end{proposition}

\begin{proof} This follows immediately from the explicit form of a Kronecker
block. \end{proof}

\begin{proposition} \label{prop4.50}\myLabel{prop4.50}\relax  The action subspace $ {\mathcal A} $ of a Kronecker pair of
bilinear forms in $ V $ of rank $ r $ has a natural Kronecker linear relation of
rank $ r $. \end{proposition}

\begin{proof} Since $ {\mathcal A} $ is isotropic with respect to $ \left(,\right)_{1} $, $ \left(,\right)_{1} $ induces a
natural pairing of $ {\mathcal A} $ with $ V/{\mathcal A} $, or a mapping $ \alpha_{1}\colon {\mathcal A} \to \left(V/{\mathcal A}\right)^{*} $. Similarly,
$ \left(,\right)_{2} $ induces a mapping $ \alpha_{2}\colon {\mathcal A} \to \left(V/{\mathcal A}\right)^{*} $. Consider the relation $ \alpha_{2}^{-1}\alpha_{1} $ in
$ {\mathcal A} $. Looking on the explicit form of a Kronecker block of a pair of
skew-symmetric bilinear forms, one can easily recognize in the relation
$ \alpha_{2}^{-1}\alpha_{1} $ a direct sum of Kronecker blocks. \end{proof}

\begin{proposition} \label{prop4.60}\myLabel{prop4.60}\relax  Consider two families $ \left(,\right)_{1,t} $, $ \left(,\right)_{2,t} $, $ t\in T $, of
skew-symmetric bilinear forms in a vector space $ V $, parameterized by a
manifold $ T $. Suppose that there is $ r\in{\mathbb N} $ such that for any $ t\in T $ the pair
$ \left(,\right)_{1,t} $, $ \left(,\right)_{2,t} $ is Kronecker of rank $ r $. Let $ {\mathcal A}_{t} $ be the action subspace of
the pair $ \left(,\right)_{1,t} $, $ \left(,\right)_{2,t} $. Then $ {\mathcal A}_{t} $ depends smoothly on $ t $. \end{proposition}

\begin{proof} This follows immediately from the following

\begin{lemma} Consider a vector space $ V $ and numbers $ n_{1},\dots ,n_{k},N\in{\mathbb N} $. Consider
the product $ X $ of Grassmannians $ \prod_{i=1}^{k}\operatorname{Gr}_{n_{i}}\left(V\right) $ and the subset $ Z\subset X $
consisting of $ k $-tuples of subspaces such that the linear span
of such a $ k $-tuple has dimension $ N $. Consider the mapping
\begin{equation}
\iota\colon Z \to \operatorname{Gr}_{N}\left(V\right)\colon \left(V_{1},\dots ,V_{k}\right) \mapsto V_{1}+\dots +V_{k},
\notag\end{equation}
here $ V_{i} $ is a vector subspace of $ V $ of dimension $ n_{i} $. Then the mapping $ \iota $ is
smooth. \end{lemma}

\begin{proof} Let $ x_{0}\in Z $, $ U_{0} $ be a small neighborhood of $ x_{0} $ in $ X $. For $ x\in X $ denote
by $ V_{i}\left(x\right) $ the $ i $-th component of $ x $, $ V_{i}\left(x\right)\in\operatorname{Gr}_{n_{i}}\left(V\right) $. Choose a basis $ v_{il}\left(x\right) $ in
$ V_{i}\left(x\right) $, $ i=1,\dots ,k $, $ l=1,\dots ,n_{i} $, which depends regularly on $ x\in U $. Pick up a
basis $ \left\{v_{\alpha}\left(x_{0}\right)\right\}_{\alpha\in A} $ out of vectors $ v_{il}\left(x_{0}\right) $, $ i=1,\dots ,k $, $ l=1,\dots ,n_{i} $ (each $ \alpha $
has a form $ il $). Then in an open subset $ U_{1} $ of $ U_{0} $ the vectors $ v_{\alpha}\left(x\right) $ remain
linearly independent, thus on $ U_{1}\cap Z $ they span $ V_{1}+\dots +V_{k} $. This implies\footnote{A similar argument can show that $ Z $ is a submanifold of $ X $.}
that $ V_{1}\left(z\right)+\dots +V_{k}\left(z\right) $ depends regularly on $ z\in Z $. \end{proof}

This finishes the proof of Proposition~\ref{prop4.60}. \end{proof}

\section{Micro-Kronecker bihamiltonian structures }\label{h4}\myLabel{h4}\relax 

Recall that a bihamiltonian structure on $ M $ induces a pair $ \left(,\right)_{1} $, $ \left(,\right)_{2} $
of bilinear skew-symmetric forms in $ {\mathcal T}^{*}M $. As a corollary, it also induces
a pair of bilinear skew-symmetric forms in $ {\mathcal T}^{*}M\otimes{\mathbb C} $, which we will denote by
the same symbols.

\begin{definition} A bihamiltonian structure on $ M $ is {\em micro-Kronecker\/} at $ m\in M $ if
the pair of bilinear skew-symmetric forms $ \left(,\right)_{1} $ and $ \left(,\right)_{2} $ in $ {\mathcal T}_{m}^{*}M $ has no
Jordan blocks. The {\em rank\/} of the bihamiltonian structure at $ m\in M $ is the
number of Kronecker blocks. A bihamiltonian structure on $ M $ is
{\em micro-Kronecker\/} of {\em rank\/} $ r $ if it is micro-Kronecker at all the points $ m\in M $
of the same rank $ r $. \end{definition}

The definition above is almost identical to the definition of
{\em completeness\/} in \cite{Bol91Com} (see also \cite{Pan99Ver}), a similar but more
restrictive definition of {\em homogeneity\/} was introduced in \cite{GelZakh99Web}.

\begin{proposition} \label{prop4.90}\myLabel{prop4.90}\relax  Consider a micro-Kronecker bihamiltonian structure
of rank $ r $ on a manifold $ M $. There is a foliation $ {\mathcal F} $ on $ M $ such that for any
$ m\in M $ the action subspace of $ {\mathcal T}_{m}^{*}M $ is the normal space to the leaf of $ {\mathcal F} $
through $ m $. \end{proposition}

\begin{proof} By Proposition~\ref{prop4.60} action subspaces from a subbundle of
$ {\mathcal T}^{*}M $. Consider $ \left(\lambda_{1},\lambda_{2}\right)\in{\mathbb K}^{2}\smallsetminus\left(0,0\right) $ and the Poisson structure $ \lambda_{1}\left\{,\right\}_{1}+\lambda_{2}\left\{,\right\}_{2} $ on
$ M $. This Poisson structure has corank $ r $, consider its symplectic foliation
$ \widetilde{{\mathcal F}}_{\lambda_{1},\lambda_{2}} $. Then the kernel of the bilinear form $ \left(,\right)_{\lambda_{1},\lambda_{2}} $ in $ {\mathcal T}_{m}^{*}M $ is the
normal space to the leaf of $ \widetilde{{\mathcal F}}_{\lambda_{1},\lambda_{2}} $ through $ m $. Now the statement follows
from the following lemmas:

\begin{lemma} \label{lm35.30}\myLabel{lm35.30}\relax  Consider a pair of skew-symmetric bilinear forms in $ V $
which has no Jordan blocks. Suppose that dimensions of all the
Kronecker components of $ V $ are $ \leq k $. If $ \left\{\lambda_{i}\right\} $, $ i\in K $, is a finite subset of $ {\mathbb K}{\mathbb P}^{1} $
with $ k $ elements, then $ V=\sum_{i\in K}\operatorname{Ker}_{\lambda_{i}}W $. \end{lemma}

\begin{proof} Direct corollary of the explicit form of a Kronecker block. \end{proof}

\begin{lemma} Consider vector subbundles $ E,E_{1},\dots ,E_{k} $ of $ {\mathcal T}^{*}M $, such that $ E=\sum E_{i} $.
Suppose that $ E_{i} $ coincides with the normal bundle to a foliation $ {\mathcal F}_{i} $ on
$ M $, $ i=1,\dots ,k $. Then there is a foliation $ {\mathcal F} $ on $ M $ such that $ E $ coincides with
the normal bundle to $ {\mathcal F} $. Leaves of $ {\mathcal F} $ can be described as intersections of
leaves of foliations $ {\mathcal F}_{i} $. \end{lemma}

\begin{proof} This follows immediately for the Frobenius integrability
criterion (Proposition~\ref{prop3.71}). \end{proof}

This finishes the proof of Proposition~\ref{prop4.90}. \end{proof}

Obviously, the foliation $ {\mathcal F} $ is canonically defined by the
bihamiltonian structure. A (most important) particular case of Proposition
~\ref{prop4.90} was announced in \cite{Pan99Ver}.

\begin{definition} \label{def4.51}\myLabel{def4.51}\relax  Given a micro-Kronecker bihamiltonian structure on a
manifold $ M $, call the foliation $ {\mathcal F} $ of Proposition~\ref{prop4.90} the
{\em action foliation\/}\footnote{In applications given a system of action-angle variables $ \left(H_{i},\varphi_{i}\right) $, the
foliation can be written as $ H_{i}\equiv \operatorname{const}_{i} $.} of the bihamiltonian structure. \end{definition}

\begin{definition} Given a micro-Kronecker bihamiltonian structure on $ M $ such that
the action foliation $ {\mathcal F} $ is a fibration, let $ {\mathcal B}_{M} $ be the base of this
fibration. If $ M $ is clear from the context, denote the base by $ {\mathcal B} $. \end{definition}

From now on suppose that the foliation $ {\mathcal F} $ is in fact a fibration of $ M $
over the base $ {\mathcal B} $. One can always achieve this by decreasing $ M $.

\begin{theorem} \label{th4.55}\myLabel{th4.55}\relax  The base of the action foliation of a micro-Kronecker
bihamiltonian structure has a canonically defined structure of a
Kronecker web. \end{theorem}

\begin{proof} Indeed, consider $ m\in M $ and the projection $ b $ of $ m $ to $ {\mathcal B} $. The
(co)differential $ \left(db|_{m}\right)^{*} $ of the mapping of projection identifies the
vector space $ {\mathcal T}_{b}^{*}{\mathcal B} $ with the action subspace at $ m $. However, by Proposition
~\ref{prop4.50} the action subspace at $ m $ is equipped with a Kroneker linear
relation, thus $ {\mathcal T}_{b}^{*}{\mathcal B} $ is equipped with such a relation as well. To show
that $ {\mathcal B} $ has a structure of a Kronecker preweb, it is enough to show that
this relation in $ {\mathcal T}_{b}^{*}{\mathcal B} $ does not depend on the choice of the point $ m $ over
$ b $.

Let $ m $, $ m' $ be two point of $ M $ over $ b\in{\mathcal B} $. Let $ {\mathcal A}_{m} $, $ {\mathcal A}_{m'} $ be the action
subspaces in $ {\mathcal T}_{m}^{*}M $ and $ {\mathcal T}_{m'}^{*}M $. Both $ {\mathcal A}_{m} $ and $ {\mathcal A}_{m'} $ are identified with $ {\mathcal T}_{b}^{*}{\mathcal B} $,
thus one with the other.

\begin{lemma} For any $ \lambda\in{\mathbb K}{\mathbb P}^{1} $ the identification between $ {\mathcal A}_{m} $ and $ {\mathcal A}_{m'} $ sends
$ \operatorname{Ker}_{\lambda}W_{m}\subset{\mathcal A}_{m} $ to $ \operatorname{Ker}_{\lambda}W_{m'}\subset{\mathcal A}_{m'} $. \end{lemma}

\begin{proof} Let $ \lambda=\lambda_{1}:\lambda_{2} $. Consider the Poisson structure $ \lambda_{1}\left\{,\right\}_{1}+\lambda_{2}\left\{,\right\}_{2} $ on $ M $.
Then $ \operatorname{Ker}_{\lambda}W $ is the normal bundle to the symplectic foliation $ \widetilde{{\mathcal F}}_{\lambda} $ for this
Poisson structure. Consider the leaf $ L $ of this foliation which passes
through $ m $. This leaf contains the leaf of the action foliation $ {\mathcal F} $ through
$ m $, thus $ L $ is a preimage of a submanifold $ \widetilde{L}\subset{\mathcal B} $. Thus $ \operatorname{Ker}_{\lambda}W_{m}\subset{\mathcal A}_{m} $ is the image
of the normal space to $ \widetilde{L} $ under the (co)differential of the projection
mapping.

By the same reason $ m' $ is in $ L $, and $ \operatorname{Ker}_{\lambda}W_{m'}\subset{\mathcal A}_{m'} $ is also the image the
same normal space. Thus the identification of $ {\mathcal A}_{m} $ and $ {\mathcal A}_{m'} $ via projection
to $ {\mathcal B} $ indeed sends $ \operatorname{Ker}_{\lambda}W_{m} $ to $ \operatorname{Ker}_{\lambda}W_{m'} $. \end{proof}

\begin{lemma} Consider vector spaces $ V $ and $ V' $ with Kronecker linear relations
$ W $ and $ W' $ in $ V $ and $ V' $ correspondingly. Consider an isomorphism $ \alpha\colon V \to V' $.
If $ \alpha $ sends $ \operatorname{Ker}_{\lambda}W $ to $ \operatorname{Ker}_{\lambda}W' $ for any $ \lambda\in{\mathbb K}{\mathbb P}^{1} $, then $ \alpha\oplus\alpha $ sends $ W\subset V\oplus V $ to
$ W'\subset V'\oplus V' $. \end{lemma}

This lemma is equivalent to Theorem~\ref{th33.10} proven in Section~\ref{h33}.

We had shown that the base $ {\mathcal B} $ of the action foliation has a
canonically defined structure of a Kronecker preweb. Denote by $ \widetilde{{\mathcal W}}_{b} $ the
relation in $ {\mathcal T}_{b}^{*}{\mathcal B} $, $ b\in{\mathcal B} $.

To show that this preweb is in fact a web, it is enough to describe
the $ \lambda $-integrating foliation $ {\mathcal F}_{\lambda} $ of $ \widetilde{{\mathcal W}} $. Fix $ \lambda_{1},\lambda_{2}\in{\mathbb K} $, consider the Poisson
structure $ \lambda_{1}\left\{,\right\}_{1}+\lambda_{2}\left\{,\right\}_{2} $ on $ M $. Since this is a Poisson structure of
constant rank, symplectic leaves form a foliation on $ M $. This foliation
depends only on $ \lambda_{1}:\lambda_{2}\in{\mathbb K}{\mathbb P}^{1} $, denote this foliation $ \widetilde{{\mathcal F}}_{\lambda_{1}:\lambda_{2}} $.

The normal space to this foliation at $ m\in M $ is the kernel of the
bilinear form $ \lambda_{1}\left(,\right)_{1}+\lambda_{2}\left(,\right)_{2} $ in $ {\mathcal T}_{m}^{*}M $. In particular, the normal bundle to
$ \widetilde{{\mathcal F}}_{\lambda_{1}:\lambda_{2}} $ is contained in the action bundle of the bihamiltonian structure.
Thus $ {\mathcal F} $ is a subfoliation of $ \widetilde{{\mathcal F}}_{\lambda_{1}:\lambda_{2}} $. In particular, $ \widetilde{{\mathcal F}}_{\lambda_{1}:\lambda_{2}} $ induces a
``quotient'' foliation $ {\mathcal F}_{\lambda_{1}:\lambda_{2}} $ on the base $ {\mathcal B} $ of the foliation $ {\mathcal F} $. Now one can
immediately see that the normal space to the foliation $ {\mathcal F}_{\lambda_{1}:\lambda_{2}} $ at $ b\in{\mathcal B} $
coincides with $ \operatorname{Ker}_{\lambda}\widetilde{{\mathcal W}}_{b} $, thus $ {\mathcal F}_{\lambda_{1}:\lambda_{2}} $ is the $ \lambda_{1}:\lambda_{2} $-integrating foliation of
the preweb on $ {\mathcal B} $. Thus the preweb structure on $ {\mathcal B} $ is indeed a Kronecker
web. \end{proof}

A statement which is parallel to a particular case of Theorem
~\ref{th4.55} was announced in \cite{Pan99Ver}.

\section{Lattice of kernels }\label{h33}\myLabel{h33}\relax 

Given a bisurjective relations $ W $ in $ V $, one obtains a collection of
vector subspaces $ \operatorname{Ker}_{\lambda}W $, $ \lambda\in{\mathbb K}{\mathbb P}^{1} $, of the vector space $ V $. The goal of this
section is to prove

\begin{theorem} \label{th33.10}\myLabel{th33.10}\relax  Let $ GR\left(V\right) $ denote the set of all vector subspaces
of a vector space $ V $. Associate to a bisurjective linear relation $ W $ in $ V $ a
mapping of sets $ {\mathcal K}_{W}\colon {\mathbb K}{\mathbb P}^{1} \to GR\left(V\right) $, $ {\mathcal K}_{W}\left(\lambda\right)=\operatorname{Ker}_{\lambda}W $.

If $ W $ and $ W' $ are two bisurjective linear relations in $ V $ such that $ {\mathcal K}_{W}={\mathcal K}_{W'} $
and $ W $ is Kronecker, then $ W=W' $. \end{theorem}

This theorem is a corollary of the following

\begin{amplification} \label{amp33.20}\myLabel{amp33.20}\relax  In the conditions of Theorem~\ref{th33.10} suppose that
all the Kronecker blocks of $ W $ have dimensions $ k $ or less. Let $ \Lambda\subset{\mathbb K}{\mathbb P}^{1} $ be a
collection of $ k+1 $ points. Then if $ {\mathcal K}_{W}|_{\Lambda}={\mathcal K}_{W'}|_{\Lambda} $, then $ W=W' $. \end{amplification}

\begin{proof}[Proof (V.~Serganova) ] Given a bisurjective relation $ {\mathbit W} $ in a vector
space $ {\mathbit V} $, suppose that $ {\mathbit K}=\operatorname{Ker}_{\lambda_{1}:\lambda_{2}}{\mathbit W} $. Define the {\em lifting\/} $ {\mathbit K}^{\left(2\right)} $ of $ {\mathbit K} $ into $ {\mathbit V}\times{\mathbit V} $
as the collection of vectors $ \left\{\left(\lambda_{1}k,\lambda_{2}k\right)\subset{\mathbit V}\times{\mathbit V} \mid k\in{\mathbit K}\right\} $. One can easily check that
$ {\mathbit K}^{\left(2\right)}\subset{\mathbit W} $.

Now suppose that $ {\mathbit W} $ is Kroneker. Let $ \Lambda=\left\{\lambda_{1},\dots ,\lambda_{k+1}\right\}\subset{\mathbb K}{\mathbb P}^{1} $, $ {\mathbit K}_{i}=\operatorname{Ker}_{\lambda_{i}}{\mathbit W} $,
and consider $ {\mathbit K}_{i}^{\left(2\right)}\subset{\mathbit W} $. It is sufficient to show that these vector
subspaces span $ {\mathbit W} $. Indeed, assume that this is true. Then the statement of
the amplification is obvious if $ W' $ is Kronecker. In general, in the
conditions of the amplification let $ V=V_{\text{Kron}}\oplus V_{\text{Jord}} $; here $ V_{\text{Kron}} $ is a sum of
Kroneker blocks of $ W' $, $ V_{\text{Jord}} $ is a sum of Jordan blocks of $ W $. Let
$ W_{\text{Kron}}'=W'\cap\left(V_{\text{Kron}}\oplus V_{\text{Kron}}\right) $ be restriction of $ W' $ to $ V_{\text{Kron}} $. We know that
vector subspaces $ {\mathbit K}_{i}^{\left(2\right)} $ span $ W $, that $ {\mathbit K}_{i}^{\left(2\right)} $ are subspaces of $ W' $, and that
$ W_{\text{Kron}}' $ is spanned by some subspaces of $ {\mathbit K}_{i}^{\left(2\right)} $, thus $ W_{\text{Kron}}'\subset W\subset W' $. If $ n $, $ n' $
are the numbers of Kronecker blocks in $ W $ and $ W'_{\text{Kron}} $, then $ \dim  W = \dim  V+n $
and $ \dim  W'=\dim  V+n' $. Hence $ n'\geq n $. On the other hand, $ n=\dim  \operatorname{Ker}_{\lambda}W $, $ n'=\dim 
\operatorname{Ker}_{\lambda}W'_{\text{Kron}}\leq\dim  \operatorname{Ker}_{\lambda}W'=\dim  \operatorname{Ker}_{\lambda}W=n $. Thus $ n'\leq n $. Hence $ \dim  W=\dim  W' $, thus $ W=W' $.

What remains to prove is that that vector subspaces $ {\mathbit K}_{i}^{\left(2\right)}\subset{\mathbit W} $ span $ {\mathbit W} $
if $ {\mathbit W} $ is Kronecker. Decomposing into a direct sum, we can restrict our
attention to the case when $ {\mathbit W} $ is one Kroneker block of dimension $ k $ or
less. By decreasing $ \Lambda $ we may assume that $ \dim  {\mathbit V}=k $. Then $ \dim  {\mathbit W}=k+1 $, $ \dim 
{\mathbit K}_{i}=\dim  {\mathbit K}_{i}^{\left(2\right)}=1 $ for any $ i $. Thus it is enough to show that $ {\mathbit K}_{i}^{\left(2\right)} $,
$ i=1,\dots ,k+1 $, are linearly independent.

Since all the Kroneker blocks of dimension $ k $ are isomorphic, it is
enough to do this for one particular Kroneker block of dimension $ k $.
Consider a $ 2 $-dimensional vector space $ {\mathcal S} $ with basis $ {\mathbit s}_{1} $, $ {\mathbit s}_{2} $ (as in
Definition~\ref{def1.20}). Let $ {\mathbit V} $ be the symmetric power $ \operatorname{Sym}^{k-1}{\mathcal S} $. Consider $ {\mathbit V} $ as
a subspace of the vector space of polynomials on $ {\mathcal S}^{*} $, then two partial
derivatives $ \partial_{1} $, $ \partial_{2} $ define two mappings $ V=\operatorname{Sym}^{k-1}{\mathcal S} \to \operatorname{Sym}^{k-2}{\mathcal S} $. Let $ {\mathbit W}\subset{\mathbit V}\times{\mathbit V} $
consists of pairs $ \left(v_{1},v_{2}\right) $ such that $ \partial_{1}v_{1}=\partial_{2}v_{2} $. One can momentarily check
that $ {\mathbit W} $ is a Kronecker block in $ {\mathbit V} $.

Obviously, $ \operatorname{Ker}_{\lambda_{1}:\lambda_{2}}{\mathbit W} $ is spanned by $ \left(\lambda_{1}{\mathbit s}_{2}-\lambda_{2}{\mathbit s}_{1}\right)^{k-1}\subset{\mathbit V} $. Let $ {\mathbit V}'=\operatorname{Sym}^{k}{\mathcal S} $,
introduce a mapping $ {\mathbit V}\times{\mathbit V} \xrightarrow[]{\alpha} {\mathbit V}' $ given by $ \left(p,p'\right) \mapsto {\mathbit s}_{2}p-{\mathbit s}_{1}p' $. One can
easily check that $ \alpha\left(\left(\operatorname{Ker}_{\lambda_{1}:\lambda_{2}}W\right)^{\left(2\right)}\right) $ is spanned by the vector
$ \left(\lambda_{1}{\mathbit s}_{2}-\lambda_{2}{\mathbit s}_{1}\right)^{k}\subset{\mathbit V}' $. Thus to show that $ {\mathbit K}_{i}^{\left(2\right)} $ are linearly independent, it is
enough to show that for $ k+1 $ non-proportional linear functions $ l_{1},\dots ,l_{k+1} $
on $ {\mathcal S}^{*} $ the polynomials $ l_{1}^{k},\dots ,v_{k+1}^{k} $ are linearly independent. In turn,
this is an obvious corollary of non-vanishing of Vandermont determinant.
\end{proof}

Later we will need the following lemma, which in the Kronecker case
is a corollary of the proof, and in Jordan case follows immediately from
the explicit description of Jordan blocks:

\begin{lemma} \label{lm35.40}\myLabel{lm35.40}\relax  Consider a relation $ W\subset V\oplus V $ which is a Kronecker block with
$ \dim  V\leq k $, or is a Jordan block. Let $ \left\{\lambda_{i}\right\} $, $ 1\leq i\leq k $, be a subset of $ {\mathbb K}{\mathbb P}^{1} $. Then
$ \operatorname{Ker}_{\lambda_{k}}W\cap\sum_{i=1}^{k-1}\operatorname{Ker}_{\lambda_{i}}W=0 $. \end{lemma}

\section{Lagrangian foliations }\label{h75}\myLabel{h75}\relax 

Here we recall more or less standard results of symplectic geometry
which will be useful in Section~\ref{h8}. See \cite{ArnGiv85Sym} for details.

\begin{definition} Given a bracket $ \left\{,\right\} $ on $ M $, a submanifold $ L\subset M $ is {\em involutive\/} if
$ \left\{f,g\right\}|_{L}=0 $ for any functions $ f $ and $ g $ on $ M $ such that\footnote{In complex-analytic situation one needs to consider functions on open
subsets of $ M $.} both $ f|_{L} $ and $ g|_{L} $ are
constant. Call $ L $ {\em Lagrangian\/} if it is involutive, and any submanifold $ L_{1}\subset L $
of codimension 1 or more is not involutive. A foliation $ {\mathcal F} $ on $ M $ is
{\em Lagrangian\/} if each leaf of $ {\mathcal F} $ is Lagrangian. \end{definition}

In the context of Example~\ref{ex002.41} the foliation on fibers of the
projection $ \pi $ is Lagrangian. The $ 0 $-section of $ \pi $ is a Lagrangian
submanifold. In fact, on a symplectic manifold $ M $ any Lagrangian
submanifold $ L $ has $ \dim  L=\frac{\dim  M}{2} $, an involutive submanifold is
Lagrangian if it has this dimension, and locally any Lagrangian foliation
can be reduced to the foliation on fibers of the projection $ \pi $ of Example
~\ref{ex002.41}:

\begin{definition} A {\em locally-affine\/} structure on a manifold $ L $ is a connection
in $ {\mathcal T}L $ with vanishing curvature and torsion. \end{definition}

\begin{proposition} \label{prop75.20}\myLabel{prop75.20}\relax  Consider a symplectic Poisson structure on $ M $, and a
Lagrangian foliation $ {\mathcal F} $. Suppose that leaves of $ {\mathcal F} $ are fibers\footnote{Locally any foliation can be represented in such a form.} of a
projection $ \pi\colon M \to N $. Let $ L $ be a leaf of the foliation $ {\mathcal F} $, $ \pi\left(L\right)=\left\{n\right\}\subset N $.
Then there is a canonical identification of $ {\mathcal T}_{l}L $ with with $ {\mathcal T}_{n}^{*}N $ for any
$ l\in L $, thus all the tangent spaces to $ L $ can be identified with each other.
This identification provides $ L $ with a locally-affine structure.

Consider a section $ s\colon N \to M $ of the projection $ \pi $. If the image $ \operatorname{Im}
s\subset M $ is a
Lagrangian submanifold, there is exactly one identification of a
neighborhood of $ \operatorname{Im} s $ with an open subset of $ {\mathcal T}^{*}N $ which is compatible with
the projections on $ N $, with Poisson structures on $ M $ and $ {\mathcal T}^{*}N $, and which
sends $ \operatorname{Im} s $ to the $ 0 $-section of $ {\mathcal T}^{*}N $.

The identification above sends the locally-affine structure on
leaves of $ {\mathcal F} $ to the tautological locally-affine structure on vector spaces
$ {\mathcal T}_{n}N $, $ n\in N $. \end{proposition}

\begin{remark} Obviously, a locally-affine structure on a
contractible set is isomorphic to the tautological connection in $ {\mathcal T}U $;
here $ U $ is
an appropriate open subset of a vector space. Thus another way to
define a locally affine structure is to introduce identifications of open
subsets $ U_{i} $ which cover $ L $ with open subsets of vector spaces, and require
that the transition functions on $ U_{i}\cap U_{j} $ correspond to affine mappings of
these vector spaces. \end{remark}

\begin{definition} \label{def8.35}\myLabel{def8.35}\relax  Consider a small open connected subset $ U $ of a manifold
$ L $ with a locally-affine structure. Then tangent spaces at different
points $ l\in U $ are identified with each other (such an identification may
depend on a choice of homotopy type of a curve which connects these
points). Call an open subset $ U\subset L $ {\em simple\/} if the identifications above do
not depend on the choice of homotopy class of connecting curves. If $ U $ is
simple, the corresponding set of equivalence classes of tangent vectors
has a vector space structure. Call this structure {\em the vector space\/}
associated to a simple locally-affine structure on $ U $. \end{definition}

\begin{corollary} Consider a Lagrangian foliation $ {\mathcal F} $ on a symplectic manifold.
Let $ B $ be a local base of $ {\mathcal F} $. By Proposition~\ref{prop75.20} leaves of $ {\mathcal F} $ have a
canonical locally-affine structure. The vector space associated to the
leaf over $ b\in B $ is canonically identified with $ {\mathcal T}_{b}^{*}B $. \end{corollary}

We will need an explicit construction of the identification of $ {\mathcal T}_{m}L $
with $ {\mathcal T}_{n}^{*}N $, $ n=\pi\left(m\right) $, from Proposition~\ref{prop75.20}. Given a Lagrangian
submanifold $ L $ in a symplectic manifold $ M $, the Hamiltonian mapping $ {\text H}\colon {\mathcal T}_{m}^{*}M
\to {\mathcal T}_{m}M $ can be restricted to a mapping $ {\mathcal N}_{m}^{*}L \to {\mathcal T}_{m}L $, $ m\in L $, which is a
bijection. On the other hand, for any leaf $ L $ of any foliation there is a
canonical flat connection on the normal bundle $ {\mathcal N}L $ to $ L $, thus on the dual
bundle $ {\mathcal N}^{*}L $. Identifying $ {\mathcal N}^{*}L $ with $ {\mathcal T}L $, one obtains a flat connection on $ {\mathcal T}L $.

Moreover, $ {\mathcal N}_{m}^{*}L $ is identified with $ {\mathcal T}_{\pi\left(m\right)}^{*}N $. Basing on these data, it is
easy to construct the identification of Proposition~\ref{prop75.20}.

Now extend this discussion to the case of Poisson structures of
constant rank.

\begin{proposition} \label{prop75.47}\myLabel{prop75.47}\relax  Any Lagrangian submanifold of a Poisson manifold $ M $ is
contained in a symplectic leaf\footnote{It is possible to define what is a symplectic leaf of an
arbitrary Poisson structure. However, for the purpose of our discussion,
it is enough to restrict attention to Poisson structures of constant
rank, as in Remark~\ref{rem00.51}.} of $ M $. If the Poisson structure on $ M $ is of
constant rank, Lagrangian foliations on $ M $ coincide with smooth families
of Lagrangian foliations, one per each symplectic leaf of $ M $. \end{proposition}

\begin{remark} \label{rem75.50}\myLabel{rem75.50}\relax  If the Lagrangian foliation is a fibration $ M \xrightarrow[]{\pi} B $, and
Poisson structure on $ M $ is of constant rank, each symplectic leaf $ S\subset M $
coincides with $ \pi^{-1}\pi S $, thus the symplectic foliation $ \widetilde{{\mathcal F}} $ on $ M $ is a preimage
of a foliation $ {\mathcal F}_{0} $ on $ B $. Due to Proposition~\ref{prop75.20}, given a section $ s $
of $ \pi $ such that intersection of $ \operatorname{Im} s $ with each symplectic leaf $ S $ is
Lagrangian in $ S $, one can identify a neighborhood of $ \operatorname{Im} s $ with a
neighborhood of the $ 0 $-section of $ {\mathcal T}^{*}{\mathcal F}_{0} $. Moreover, $ {\mathcal T}^{*}{\mathcal F}_{0} $ carries a natural
Poisson structure, and the identification above is compatible with
Poisson structures on $ {\mathcal T}^{*}{\mathcal F}_{0} $ and $ M $. \end{remark}

\section{Two functors }\label{h8}\myLabel{h8}\relax 

Presently we have two constructions: by Corollary~\ref{cor33.99}, given a
Kronecker web $ {\mathcal W} $ on $ B $ of rank $ r $, one can construct a bihamiltonian
structure on the total space of $ \Phi\left({\mathcal W}\right) $ (which has dimension $ 2 \dim  B-r $). By
Theorem~\ref{th4.55}, given a micro-Kronecker bihamiltonian structure $ M $ of
dimension $ d $ and rank $ r $, one can associate to small open subsets $ U\subset M $ a
Kronecker web structure on manifold $ {\mathcal B}_{U} $ (of dimension $ \frac{d+r}{2} $).
Investigate the relation of these two constructs.

\begin{proposition} \label{prop8.27}\myLabel{prop8.27}\relax  Consider a Kronecker web $ {\mathcal W} $ of rank $ r $ on a manifold $ B $.
\begin{enumerate}
\item
If $ \left(\lambda_{1},\lambda_{2}\right)\in{\mathbb K}^{2}\smallsetminus\left(0,0\right) $, then the symplectic foliation of the Poisson
bracket $ \lambda_{1}\left\{,\right\}_{1}+\lambda_{2}\left\{,\right\}_{2} $ on $ \Phi\left({\mathcal W}\right) $ is the preimage of the integrating
foliation $ {\mathcal F}_{\lambda_{1}:\lambda_{2}} $ on $ B $.
\item
The bihamiltonian structure on $ \Phi\left({\mathcal W}\right) $ is micro-Kronecker of rank $ r $.
\item
The action foliation on $ \Phi\left({\mathcal W}\right) $ coincides with foliation on fibers of
projection $ \Phi\left({\mathcal W}\right) \to B $.
\item
The structure of Kronecker web on $ B $ induced by the bihamiltonian
structure on $ \Phi\left({\mathcal W}\right) $ coincides with the initial Kronecker web structure $ {\mathcal W} $ on
$ B $.
\end{enumerate}
\end{proposition}

\begin{proof} The first statement follows from the definition of the Poisson
bracket $ \left\{,\right\}_{\lambda_{1}:\lambda_{2}} $ on $ \Phi\left({\mathcal W}\right) $. If $ {\mathbb K}={\mathbb C} $, the second statement is a direct
corollary of the first one. If $ {\mathbb K}={\mathbb R} $, the second statement follows from the
following

\begin{lemma} Consider a pair of skew-symmetric bilinear forms $ \left(,\right)_{1} $ and $ \left(,\right)_{2} $
in a vector space $ V $. Then $ \dim  \operatorname{Ker}\left(\lambda_{1}\left(,\right)_{1}+\lambda_{2}\left(,\right)_{2}\right) $ is constant for $ \left(\lambda_{1},\lambda_{2}\right) $
inside an open subset of $ {\mathbb K}^{2} $, call this common value $ r $. Suppose that for
any finite subset $ \Lambda_{0}\subset{\mathbb K}{\mathbb P}^{1} $ there is $ \Lambda\subset{\mathbb K}{\mathbb P}^{1}\smallsetminus\Lambda_{0} $ such the vector subspaces
$ \operatorname{Ker}\left(\lambda_{1}\left(,\right)_{1}+\lambda_{2}\left(,\right)_{2}\right) $, $ \lambda_{1}:\lambda_{2}\in\Lambda $, span a vector subspace of $ V $ of dimension
$ \frac{\dim  V + r}{2} $. Then the pair $ \left(,\right)_{1} $, $ \left(,\right)_{2} $ is Kronecker. \end{lemma}

\begin{proof} Lemma follows immediately from the classification of Theorem
~\ref{th4.35}. \end{proof}

The remaining statements of the proposition follow immediately from
the first two statements.\end{proof}

\begin{proposition} \label{prop8.40}\myLabel{prop8.40}\relax  Consider a micro-Kronecker bihamiltonian structure
on a manifold $ M $. Then
\begin{enumerate}
\item
the fibers of action foliations have a natural locally-affine
structure;
\item
suppose that the action foliation is a fibration with a base $ {\mathcal B} $, and
the locally-affine structures on fibers are simple. Denote by $ E $ the vector
bundle over $ {\mathcal B} $ which is associated\footnote{As in Definition~\ref{def8.35}.} to the bundle of locally-affine
structures mentioned above. This vector bundle is canonically isomorphic to
$ \Phi\left({\mathcal B}\right) $.
\end{enumerate}
\end{proposition}

\begin{proof} Consider Hamiltonian mappings $ {\text H}_{1,2}\colon {\mathcal T}_{m}^{*}M \to {\mathcal T}_{m}M $ of Poisson
structures on $ M $. Restrict these mappings to the action subspace $ {\mathcal A}_{m}\subset{\mathcal T}_{m}^{*}M $.
By Proposition~\ref{prop4.40} the action subspace is isotropic with respect to
both pairings in $ {\mathcal T}_{m}^{*}M $, thus the elements $ {\text H}_{1,2}a $, $ a\in{\mathcal A}_{m} $, are orthogonal\footnote{W.r.t. the standard duality between $ {\mathcal T}_{m}^{*}M $ and $ {\mathcal T}_{m}M $.}
to $ {\mathcal A}_{m} $. By definition, $ {\mathcal A}_{m}^{\perp} $ coincides with the tangent space to the fiber
$ L_{m} $ of the action foliation which passes through $ m $. In particular, there
are two mappings $ \widetilde{{\text H}}_{1,2}\colon {\mathcal A}_{m} \to {\mathcal T}_{m}L_{m} $. Note that these mappings form a pencil
which corresponds to the linear relation in $ {\mathcal A}_{m} $ (as in Proposition
~\ref{prop4.50} and in the proof of Theorem~\ref{th4.55}).

Since the statement of the proposition is local on $ M $, one can
suppose that the action foliation is a fibration. Let $ \pi\colon M \to {\mathcal B} $ be the
projection to the base of foliation, and $ b=\pi\left(m\right) $. Then $ \pi^{*} $ induces a
canonical isomorphism $ {\mathcal T}_{b}^{*}{\mathcal B}\simeq{\mathcal A}_{m} $. Now $ {\mathcal B} $ is a Kronecker web, thus $ {\mathcal T}_{b}^{*}{\mathcal B} $
carries a Kronecker relation, and this relation is compatible with the
relation on $ {\mathcal A}_{m} $ w.r.t.~the isomorphism above. On the other hand, let $ \Phi_{b} $ be
the fiber of $ \Phi\left({\mathcal B}\right) $ over $ b $. Then the linear relation in $ {\mathcal T}_{b}^{*}{\mathcal B} $ can be
described both by a pencil $ \widetilde{{\mathcal P}}_{1,2}=\widetilde{{\text H}}_{1,2}\circ\pi^{*}\colon {\mathcal T}_{b}^{*}{\mathcal B} \to {\mathcal T}_{m}L_{m} $ and by a pencil $ {\mathcal P}_{1,2}:
{\mathcal T}_{b}^{*}{\mathcal B} \to \Phi_{b} $. Proposition~\ref{prop1.50} gives a canonical isomorphism between
$ {\mathcal T}_{m}L_{m} $ and $ \Phi_{b} $, $ b=\pi\left(m\right) $. As a corollary there is a canonical identification
of tangent spaces to $ L_{m} $ at different points (since they all project to
the same point of $ {\mathcal B} $). This induces a flat connection $ \nabla $ on the tangent
bundle to $ L_{m} $.

What remains is to show that this connection on $ L_{m} $ is a
locally-affine structure. Pick up any $ \lambda\in{\mathbb P}^{1} $, for example, $ \lambda=1:0 $. Consider
the corresponding Poisson structure $ 1\left\{,\right\}_{1}+0\left\{,\right\}_{2} $ of the pencil.

By Proposition~\ref{prop4.40}, $ \left\{\pi^{*}\varphi_{1},\pi^{*}\varphi_{2}\right\}_{1}=0 $ for any functions $ \varphi_{1} $, $ \varphi_{2} $ on
$ {\mathcal B} $. This implies that fibers of $ {\mathcal F} $ are involutive submanifolds of $ M $ w.r.t.~
$ \left\{,\right\}_{1} $. On the other hand, if the rank of bihamiltonian structure on $ M $ is
$ r $, then symplectic leaves of $ M $ have dimension $ \dim  M -r $, and $ L_{m} $ has
dimension $ \frac{\dim  M-r}{2} $. Since $ L_{m} $ has half the dimension of the symplectic
leaf, and is involutive, it is Lagrangian. Thus $ {\mathcal F} $ is a Lagrangian
foliation w.r.t.~$ \left\{,\right\}_{1} $. Now Propositions~\ref{prop75.20} and~\ref{prop75.20} imply
that the leaves of $ {\mathcal F} $ are equipped with canonical locally-affine structures.

Thus it is enough to prove that the connection $ \nabla $ on $ L_{m} $ constructed
above coincides with the connection of this locally-affine structure.
Recall that the connection $ \nabla $ on $ L_{m} $ is constructed basing on the mappings
$ \widetilde{{\mathcal P}}_{1}\colon {\mathcal T}_{b}^{*}{\mathcal B} \to {\mathcal T}_{m}L_{m} $, in other words, on the isomorphisms $ {\mathcal T}_{b}^{*}{\mathcal B}/\operatorname{Ker}\left(\widetilde{{\mathcal P}}_{1}\right)\simeq{\mathcal T}_{m}L_{m} $
for different $ m $ with $ \pi\left(m\right)=b $. However, $ {\mathcal T}_{b}^{*}{\mathcal B}/\operatorname{Ker}\left(\widetilde{{\mathcal P}}_{1}\right) $ coincides with $ {\mathcal T}_{b}^{*}F_{b} $;
here $ F_{b} $ is the fiber of integrating foliation $ {\mathcal F}_{1:0} $ which passes through
$ b\in{\mathcal B} $. Let $ \widetilde{F}=\pi^{-1}F_{b} $, recall that $ \widetilde{F} $ is a symplectic leaf of the Poisson
structure $ \left\{,\right\}_{1} $ on $ M $. The foliation $ {\mathcal F} $ allows a restriction $ {\mathcal F}|_{\widetilde{F}} $ to $ \widetilde{F} $, this
restriction is a Lagrangian foliation on a symplectic manifold. Both the
identification $ {\mathcal T}_{b}^{*}F_{b}={\mathcal T}_{b}^{*}{\mathcal B}/\operatorname{Ker}\left(\widetilde{{\mathcal P}}_{1}\right)\simeq{\mathcal T}_{m}L_{m} $ and the identification of $ {\mathcal T}^{*}N $ with
$ {\mathcal T}_{m}L $ in Section~\ref{h75} are constructed as restrictions of Hamiltonian
mappings. Setting $ N=F_{b} $, $ L=L_{m} $, and replacing $ {\mathcal F} $ with $ {\mathcal F}|_{\widetilde{F}} $ shows that the
connection $ \nabla|_{\widetilde{F}} $ coincides with the connection of the locally-affine
structure of Section~\ref{h75}. This finishes the proof. \end{proof}

\begin{remark} \label{rem8.45}\myLabel{rem8.45}\relax  In the proof above we worked with two locally-affine
structures on $ L_{m} $: one constructed basing on the bihamiltonian structure,
another based on the Poisson structure $ \left\{,\right\}_{1} $. However, one could also
consider another Poisson structure $ \lambda_{1}\left\{,\right\}_{1}+\lambda_{2}\left\{,\right\}_{2} $. The fact that two
locally-affine structures of the proof coincide shows that the
locally-affine structures on $ L_{m} $ which correspond to the Poisson structure
of the pencil $ \lambda_{1}\left\{,\right\}_{1}+\lambda_{2}\left\{,\right\}_{2} $ do not depend on $ \left(\lambda_{1},\lambda_{2}\right)\not=\left(0,0\right) $. \end{remark}

\section{Conjecture on classification }\label{h9}\myLabel{h9}\relax 

\begin{conjecture} \label{con8.50}\myLabel{con8.50}\relax  Consider a micro-Kronecker bihamiltonian structure on
a manifold $ M $. Suppose that the restriction of the action foliation to
$ U_{0}\subset M $ is a fibration with a base $ {\mathcal B}_{U_{0}} $, $ m\in U_{0} $, and $ b $ is the projection of $ m $
to $ {\mathcal B} $. Let $ {\mathcal W} $ be the Kronecker web on $ {\mathcal B}_{U_{0}} $ induced by the bihamiltonian
structure on $ M $. Let $ m'\in\Phi\left({\mathcal W}\right) $ be the point on the $ 0 $-section of $ \Phi\left({\mathcal W}\right) $ over $ b $.
Then there is a local isomorphism of bihamiltonian structures on $ M $ near $ m $
and on $ \Phi\left({\mathcal W}\right) $ near $ m' $. \end{conjecture}

This conjecture is a generalization of a conjecture from \cite{GelZakh99Web}.
The case of analytic bihamiltonian structure on $ M $ of rank 1 is claimed in
\cite{GelZakhWeb}, the $ C^{\infty} $-case of rank 1 is claimed in \cite{Tur99Equi} and
proven in \cite{Tur99MemA}.

It looks like this conjecture would easily imply all the conjectures
of \cite{GelZakh99Web}. (Compare with the way we prove Theorem~\ref{th11.70}.)

\begin{remark} The conjecture implies that to check existence of a local
isomorphism between two micro-Kronecker bihamiltonian structures on $ M $ and
$ M' $, it is enough to check existence of a local isomorphism between
structures of Kronecker webs on local bases of action foliations. By
Amplification~\ref{amp33.20} to check an isomorphism of Kronecker webs on $ B $
and $ B' $ it is enough to restrict attention to a finite number of
foliations on $ B $ and $ B' $. It so happens that for webs which appear in
problems of mathematical physics the latter problem is very easy (much
easier than for generic Kronecker webs).

In what follows we establish a particular case of the conjecture,
and in Section~\ref{h11} we use this particular case to show flatness of some
bihamiltonian systems. After showing that the conjecture is applicable to
these systems, the only thing which remains to be proven is that the
structure of the Kronecker web on the base of the action foliation is
{\em flat}, i.e., locally isomorphic to a translation-invariant Kronecker web.
\end{remark}

\begin{remark} Yet another reformulation of the conjecture is that to check
that there is a local isomorphism between two micro-Kronecker
bihamiltonian structures on $ M $ and $ M' $ one should concentrate attention on
{\em action variables\/} (or first integrals) of these structures, and one can
completely ignore the question of {\em angle variables}. Note that this
enormously simplifies the question: for example, usually first integrals
may be explicitly written as polynomials in appropriate coordinates,
while all one can say about phase variables is that they satisfy some
partial differential equations, and that they may be expressed in terms
of $ \vartheta $-functions. \end{remark}

In this remark we used the following relationship between {\em action
variables\/} (i.e., functions which are constant on fibers of action
foliation) and {\em first integrals\/} (i.e., particular coordinate functions on
a bihamiltonian manifold which arise when one {\em integrates\/} the system): for
a micro-Kronecker bihamiltonian system first integrals obtained by Lenard
scheme can be taken as action variables. For details see \cite{GelZakh99Web}.

We prefer an abstract consideration of the (invariantly-defined)
action foliation to the language of first integrals, which (at least a
priori) may depend on the choice of the particular scheme of integration.

\section{Anti-involutions }\label{h91}\myLabel{h91}\relax 

\begin{definition} An {\em involution\/} is a mapping $ f\colon X \to X $ such that $ f\circ f=\operatorname{id} $.
An involution $ \alpha\colon M \to M $ is an {\em antiinvolution\/} of a Poisson structure
$ \left\{,\right\} $ on a manifold $ M $ if $ \left\{\alpha^{*}f,\alpha^{*}g\right\}=-\alpha^{*}\left\{f,g\right\} $ for any two functions $ f,g $ on $ M $.
An {\em antiinvolution\/} of a bihamiltonian structure on a manifold $ M $ is a
mapping $ M \to M $ which is an anti-involution of both Poisson brackets $ \left\{,\right\}_{1} $,
$ \left\{,\right\}_{2} $. \end{definition}

\begin{proposition} \label{prop9.20}\myLabel{prop9.20}\relax  Let $ {\mathcal W} $ be a Kronecker web on a manifold $ B $. The
bihamiltonian structure on $ \Phi\left({\mathcal W}\right) $ allows an anti-involution. \end{proposition}

\begin{proof} $ \Phi\left({\mathcal W}\right) $ is a vector bundle over $ B $. Let $ \alpha $ be the multiplication by
$ -1 $ in this bundle. To show that $ \alpha $ is an anti-involution, note that $ \left\{,\right\}_{1} $ is
isomorphic to the Poisson structure on $ {\mathcal T}^{*}{\mathcal F}_{1:0} $ (here $ {\mathcal F}_{1:0} $ is the
integrating foliation on $ B $), and this isomorphism is a mapping of vector
bundles. Thus it is enough to show that multiplication by $ -1 $ is an
anti-involution on the cotangent bundle to a foliation. In turn, it is
enough to prove this for the cotangent bundle to a manifold, which can be
checked immediately in local coordinates. \end{proof}

Consider an anti-involution $ \iota $ of a micro-Kronecker bihamiltonian
structure on a manifold $ M $. Let $ m\in M $ be a fixed point of $ \iota $. There is an
$ \iota $-invariant open subset $ U\ni m $ such that the restriction of the action
foliation $ {\mathcal F} $ on $ U $ is a fibration, denote the base of this fibration by $ {\mathcal B} $.
Since the action foliation is determined by the {\em set\/} of linear
combinations $ \lambda_{1}\left\{,\right\}_{1}+\lambda_{2}\left\{,\right\}_{2} $ of Poisson structures, and not by the
parameterization $ \left(\lambda_{1},\lambda_{2}\right) $ of this set, $ \iota $ sends leaves of $ {\mathcal F} $ to leaves of $ {\mathcal F} $.
Thus $ \iota $ induces an involution $ \iota_{{\mathcal B}} $ of $ {\mathcal B} $. Consider the submanifold\footnote{Recall that the space of fixed point of an involution $ i\colon M \to M $ on a
manifold $ M $ is a submanifold $ \operatorname{Fix}\left(i\right)\subset M $, and for $ f\in\operatorname{Fix}\left(i\right) $ the tangent space
to $ \operatorname{Fix}\left(i\right) $ at $ f $ coincides with the invariant subspace of the
differential $ i_{*}|_{f} $ of $ i $. Indeed,
consider an $ i $-invariant Riemannian metric on $ M $, then the exponential
mapping of this mapping commutes with $ i $.} $ \operatorname{Fix}\left(\iota_{{\mathcal B}}\right)\subset{\mathcal B} $
of fixed points of $ \iota_{{\mathcal B}} $.

\begin{definition} The {\em defect\/} of anti-involution $ \iota $ at the fixed point $ m\in M $ is
the codimension of $ \operatorname{Fix}\left(\iota_{{\mathcal B}}\right)\subset{\mathcal B} $ at the image $ b\in{\mathcal B} $ of $ m $. \end{definition}

Note that the defect is locally constant on $ \operatorname{Fix}\left(\iota\right) $.

Obviously, the anti-involution of Proposition~\ref{prop9.20} has defect
0. Now Conjecture~\ref{con8.50} immediately implies

\begin{conjecture} \label{con9.45}\myLabel{con9.45}\relax  Consider a micro-Kronecker bihamiltonian structure on
a manifold $ M $, $ m\in M $. There is a neighborhood $ U\ni m $ and an anti-involution $ \alpha $
of the bihamiltonian structure on $ U $ such that $ \alpha\left(m\right)=m $ and the defect of $ \alpha $
at $ m $ is 0. \end{conjecture}

Anti-involutions with defect 0 are important for us because of the
following

\begin{theorem} \label{th9.50}\myLabel{th9.50}\relax  Suppose that a bihamiltonian structure on a manifold $ M $
allows an anti-involution $ \iota $ such that $ \iota\left(m\right)=m $, $ m\in M $, and the defect of $ \iota $ at
$ m $ is 0. Then Conjecture~\ref{con8.50} holds for $ m $ and $ M $. \end{theorem}

The theorem is an immediate corollary of Propositions~\ref{prop9.51},
~\ref{prop9.52} below.

\begin{definition} Consider a micro-Kronecker bihamiltonian structure on $ M $ with
the action foliation $ {\mathcal F} $. A submanifold $ S\subset M $ is a {\em cross-section\/} if
$ S $ is a transversal section of foliation $ {\mathcal F} $, and for any $ \left(\lambda_{1},\lambda_{2}\right)\in{\mathbb K}^{2}\smallsetminus\left(0,0\right) $
and any symplectic leaf $ L $ of $ \lambda_{1}\left\{,\right\}_{1}+\lambda_{2}\left\{,\right\}_{2} $ the intersection $ L\cap S $ is
Lagrangian in $ L $. \end{definition}

\begin{proposition} \label{prop9.51}\myLabel{prop9.51}\relax  Given a micro-Kronecker bihamiltonian structure on $ M $
with a base $ {\mathcal B} $ of the action foliation, and a cross-section $ S $, there is a
local isomorphism of bihamiltonian structures $ M $ and $ \Phi\left({\mathcal B}\right) $ which sends $ S $ to
$ 0 $-section of $ \Phi\left({\mathcal B}\right) $. \end{proposition}

\begin{proof} Given a locally-affine structure $ L $ and a point $ l\in L $, the
exponential mapping of the connection on $ {\mathcal T}L $ gives a canonical
identification of a neighborhood of $ l $ in $ L $ and a neighborhood of 0 in
$ {\mathcal T}_{l}L $. Similarly, given a bundle $ {\mathcal L} \to B $ with fibers having a locally-affine
structure, and a section $ s $ of this bundle, one obtains a canonical
identification of a neighborhood of $ \operatorname{Im} s\subset{\mathcal L} $ with a neighborhood of
$ 0 $-section in an appropriate vector bundle $ {\mathcal E} \to B $. Obviously, fibers of $ {\mathcal E} $
are vertical tangent space of $ {\mathcal L} $ at points of $ \operatorname{Im} s $. In conditions of the
proposition $ {\mathcal E}=\Phi\left(B\right) $.

In particular, any submanifold $ S\subset M $ which is transversal to leaves of
$ {\mathcal F} $, and such that $ \dim  S=\operatorname{codim} {\mathcal F} $, gives a natural identification $ \gamma $ of a
neighborhood of $ S $ with a neighborhood of $ 0 $-section of $ \Phi\left({\mathcal B}\right) $ (one can
identify $ S $ with the local base $ {\mathcal B} $). We need to show that if $ S $ is a
cross-section, then this identification is compatible with bihamiltonian
structures on $ M $ and on $ \Phi\left({\mathcal B}_{U}\right) $.

It is enough to show that $ \gamma $ is compatible with $ \left\{,\right\}_{1} $ (the same
argument will be applicable to $ \left\{,\right\}_{2} $). Consider a symplectic leaf $ \widetilde{L}\subset M $ for
$ \left\{,\right\}_{1} $. It is a preimage of a submanifold $ L\subset{\mathcal B}_{M} $. The restriction foliation
$ {\mathcal F}|_{\widetilde{L}} $ is well-defined and is a Lagrangian foliation. Moreover, $ L\cap S $ is
Lagrangian. Consider the identification of a neighborhood of $ \widetilde{L}\cap S $ in $ \widetilde{L} $
with a neighborhood of $ 0 $-section of $ {\mathcal T}^{*}L $ (see Proposition~\ref{prop75.20}).
Since $ \widetilde{L}\cap S $ is Lagrangian w.r.t.~$ \left\{,\right\}_{1} $, this identification is compatible
with the Poisson structure $ \left\{,\right\}_{1} $. Composing this identification with $ \gamma $,
one obtains an identification $ \widetilde{\gamma} $ of $ {\mathcal T}^{*}L $ with a submanifold of $ \Phi\left({\mathcal B}_{M}\right) $.

Obviously, $ \widetilde{\gamma} $ sends a fiber of $ {\mathcal T}^{*}L $ to a fiber of $ \Phi\left({\mathcal B}_{M}\right) $, and $ 0 $-section
of $ {\mathcal T}^{*}L $ into a $ 0 $-section of $ \Phi\left({\mathcal B}_{M}\right) $. Due to Remark~\ref{rem8.45}, $ \widetilde{\gamma} $ is compatible
with affine structures on fibers of $ {\mathcal T}^{*}L $ and of $ \Phi\left({\mathcal B}_{M}\right) $, and since it sends
0 to 0, it is a linear mapping of vector bundles $ {\mathcal T}^{*}L \xrightarrow[]{\widetilde{\gamma}} \Phi|_{L} $. However,
given $ \left(\lambda_{1},\lambda_{2}\right)\in{\mathbb K}^{2}\smallsetminus\left(0,0\right) $, $ \Phi|_{L} $ is canonically isomorphic to $ {\mathcal T}^{*}{\mathcal F}_{\lambda_{1}:\lambda_{2}} $ (as in
Section~\ref{h34}), and by Proposition~\ref{prop8.27} $ L $ is a leaf of the
foliation $ {\mathcal F}_{1:0} $. Thus there is a canonical isomorphism $ i\colon {\mathcal T}^{*}L \to \Phi|_{L} $ which
is compatible with the Poisson structure $ \left\{,\right\}_{1} $ on $ \Phi $.

It is enough to show that $ \widetilde{\gamma}=i $. But this is a direct corollary of the
proof of Proposition~\ref{prop8.40}. \end{proof}

\begin{proposition} \label{prop9.52}\myLabel{prop9.52}\relax  Consider a micro-Kronecker bihamiltonian structure
on $ M $. If $ \iota $ is an anti-involution of $ M $ of defect 0, then any connected
component of $ \operatorname{Fix}\left(\iota\right) $ with index 0 is a cross-section. \end{proposition}

\begin{proof} This proposition follows immediately from the following

\begin{lemma} \label{lm9.90}\myLabel{lm9.90}\relax  Consider a symplectic manifold $ M $ with a Lagrangian
foliation $ {\mathcal F} $ and an anti-involution $ \iota\colon M \to M $ which sends each leaf of $ {\mathcal F} $
into itself. Then fixed points of $ \iota $ form a Lagrangian submanifold of $ M $
which is transversal to leaves of $ {\mathcal F} $. \end{lemma}

\begin{proof} It is enough to restrict our attention to a neighborhood of
$ \operatorname{Fix}\left(\iota\right) $, thus one may assume that $ M $ is an open subset of $ {\mathcal T}^{*}B $, and $ {\mathcal F} $ is the
foliation on fibers of projection $ M \to B $. Since the locally-affine
structure on fibers of a Lagrangian foliation does not change if one
multiplies the Poisson bracket by a number, the restriction of $ \iota $ on any
leaf of $ {\mathcal F} $ induces an affine transformation of this leaf (which is an open
subset of $ {\mathcal T}_{b}B $ for an appropriate $ b\in B $).

Recall that the set of fixed points of an involution $ \iota $ of $ M $ forms a
manifold $ F\subset M $, and if $ f\in F $, then $ {\mathcal T}_{f}F\subset{\mathcal T}_{f}M $ coincides with the invariant
subspace of $ \left(d\iota\right)|_{{\mathcal T}_{f}M} $.

\begin{lemma} Consider a symplectic vector space $ V $, and a Lagrangian subspace
$ l\subset V $. Consider a linear involution $ i\colon V \to V $ such that $ \left[i v_{1},i v_{2}\right]=-\left[v_{1},v_{2}\right] $
for any $ v_{1},v_{2}\in V $. Suppose that $ il=l $, and that $ i $ induces an identity
mapping of $ V/l $ into itself. Then $ i|_{l} $ is multiplication by $ -1 $, and $ l $ is
a complement to the invariant subspace of $ i $ in $ V $, which is a Lagrangian
subspace of $ V $. \end{lemma}

\begin{proof} Let $ V_{1} $ is the vector subspace of fixed points of $ i $, $ V_{-1} $ be the
eigenspace of $ i $ with eigenvalue $ -1 $. Obviously $ V=V_{1}\oplus V_{-1} $, and due to the
conditions on $ i $ both $ V_{1} $ and $ V_{-1} $ are isotropic. Thus both $ V_{1} $ and $ V_{-1} $ are
Lagrangian. Similarly, $ l=l_{1}\oplus l_{-1} $, $ l_{1}\subset V_{1} $, $ l_{-1}\subset V_{-1} $. Let $ l_{\lambda}' $ be any
complementary subspace to $ l_{\lambda} $ in $ V_{\lambda} $, $ \lambda\in\left\{1,-1\right\} $. Since the action of $ i $ in
$ V/l $ is isomorphic to the action of $ i $ in $ l_{1}'\oplus l'_{-1} $, we see that $ l'_{-1}=0 $,
thus $ l\supset V_{-1} $. Since $ l $ is Lagrangian, $ l=V_{-1} $, which finishes the proof. \end{proof}

Applying this lemma to $ V={\mathcal T}_{f}M $, $ l={\mathcal T}_{f}L $ (here $ f\in M $ is such that $ i f=f $, $ L $
is the leaf of $ {\mathcal F} $ through $ f $) finishes the proof of Lemma~\ref{lm9.90}. \end{proof}

This finishes the proof of Proposition~\ref{prop9.52}. \end{proof}

\section{Real and complex bihamiltonian structures }\label{h92}\myLabel{h92}\relax 

Consider a real-analytic bihamiltonian structure $ \left\{,\right\}_{1,2} $ on $ M $. One
can construct a {\em small complex-analytic\/} neighborhood $ M_{{\mathbb C}}\supset M $ of the manifold
$ M $ such that the brackets $ \left\{,\right\}_{1,2} $ can be analytically continued on $ M_{{\mathbb C}} $. If
the bihamiltonian structure on $ M $ is micro-Kronecker, one may assume the
same about $ M_{{\mathbb C}} $ (possibly, after decreasing $ M_{{\mathbb C}} $).

\begin{proposition} \label{prop95.10}\myLabel{prop95.10}\relax  Let $ m\in M\subset M_{{\mathbb C}} $, suppose that $ U\subset M_{{\mathbb C}} $, $ m\in U $, allows an
anti-involution $ \iota\colon U \to U $ with $ m=\iota m $ and defect 0 near $ m $. Then there is a
neighborhood $ U_{1} $ of $ m $ in $ M $ which allows an anti-involution $ \iota'\colon U' \to U' $
with $ m=\iota'm $ and defect 0 near $ m $. \end{proposition}

\begin{proof} Recall that

\begin{definition} An {\em antiholomorphic mapping\/} $ C\colon N \to N' $ of complex analytic
manifolds $ N $ and $ N' $ is such a mapping of sets that $ \overline{C^{*}\varphi} $ is a holomorphic
function on $ N $ if $ \varphi $ is a holomorphic function on $ N' $. \end{definition}

Clearly, if $ C $ is antiholomorphic, and $ Z\subset N' $ is a complex-analytic
subvariety, then $ C^{-1}Z $ is a complex analytic subvariety. Indeed, if $ \varphi=0 $ is
an equation of $ Z $, then $ \overline{C^{*}\varphi}=0 $ is an equation of $ C^{-1}Z $. Similarly, one can
transfer a Kronecker structure and a bihamiltonian structure via an
antiholomorphic bijection. If $ N_{{\mathbb C}} $ is a complexification of a real-analytic
manifold $ N $, then $ N_{{\mathbb C}} $ is equipped with an antiholomorphic involution $ C $ such
that $ \operatorname{Fix}\left(C\right)=N $. If $ N $ has a Kronecker or bihamiltonian structure, so has
$ N_{{\mathbb C}} $, and the structures on $ N_{{\mathbb C}} $ are invariant w.r.t.~$ C $.

Consider the antiholomorphic involution $ C $ for $ M_{{\mathbb C}} $, $ \operatorname{Fix}\left(C\right)=M $. Consider
a cross-section $ \sigma $ of the projection $ \pi_{{\mathbb C}} $ which passes through $ m $. Then $ C\sigma $ is
also a cross-section which passes through $ m $. Since fibers of $ {\mathcal F} $ have an
locally-affine structure, the section $ \frac{\sigma+C\sigma}{2} $ is well defined on an
appropriate neighborhood of $ \pi m\in{\mathcal B} $. Moreover, this section is a
cross-section. Indeed, this immediately follows

\begin{lemma} \label{lm9.87}\myLabel{lm9.87}\relax  Consider a Lagrangian foliation $ {\mathcal F} $ on a symplectic manifold $ M $,
and submanifolds $ S_{1} $, $ S_{2} $, $ S_{3} $ of dimension $ \frac{\dim  M}{2} $ which are transversal
to $ {\mathcal F} $. Suppose that for any leaf $ L $ of $ {\mathcal F} $ the intersection $ S_{i}\cap L $, $ i=1,2,3 $,
consists of one point $ p_{i} $, and $ p_{2} $ is the midpoint of the segment $ p_{1}p_{3} $ in
the locally affine structure on $ L $. If $ L_{1} $ and $ L_{3} $ are Lagrangian, so is $ L_{2} $.
\end{lemma}

\begin{proof} Indeed, we may assume that $ M={\mathcal T}^{*}N $, leaves of $ {\mathcal F} $ are fibers of $ \pi:
M \to N $, $ S_{1} $ is $ 0 $-section of $ {\mathcal T}^{*}N $, $ S_{2,3} $ are graphs of sections $ \varepsilon_{2,3} $ of
$ {\mathcal T}^{*}N=\Omega^{1}N $, and $ \varepsilon_{3}=2\varepsilon_{2} $.

\begin{lemma} \label{lm9.88}\myLabel{lm9.88}\relax  Consider a section $ \varepsilon $ of $ {\mathcal T}^{*}B $. The graph of $ \varepsilon $ (which is a
submanifold of $ {\mathcal T}^{*}B $) is Lagrangian iff $ d\varepsilon=0\in\Omega^{2}B $. \end{lemma}

Application of this obvious lemma finishes the proof of Lemma
~\ref{lm9.87}. \end{proof}

Now the section $ \frac{\sigma+C\sigma}{2} $ of $ \pi $ is $ C $-invariant, thus it is a
complexification of the section $ \sigma_{{\mathbb R}} $ for $ M $. Since the complexification of
$ \sigma_{{\mathbb R}} $ is a cross-section, so is $ \sigma_{{\mathbb R}} $. This finishes the proof of the
proposition. \end{proof}

\section{Endomorphisms of bihamiltonian structures and the principal theorem }\label{h10}\myLabel{h10}\relax 

The following statement is widely known:

\begin{proposition} \label{prop10.30}\myLabel{prop10.30}\relax  Suppose that $ \left(\lambda_{1},\lambda_{2}\right) $ and $ \left(\lambda'_{1},\lambda'_{2}\right) $ are two
non-proportional vectors in $ {\mathbb K}^{2} $. Consider a bihamiltonian structure $ \left\{,\right\}_{1,2} $
on a manifold $ M $. Suppose that a function $ F $ on $ M $ is a Casimir function for
the bracket $ \lambda_{1}\left\{,\right\}_{1}+\lambda_{2}\left\{,\right\}_{2} $. Let $ \chi\in\operatorname{Vect} M $ be the Hamiltonian vector field\footnote{I.e., $ \chi={\text H}\left(dF\right) $; here $ {\text H} $ is the Hamiltonian mapping $ {\mathcal T}^{*}\left(M\right) \to {\mathcal T}\left(M\right) $.}
of $ F $ w.r.t.~the bracket $ \lambda'_{1}\left\{,\right\}_{1}+\lambda'_{2}\left\{,\right\}_{2} $. Consider the flow $ \alpha_{t} $ of $ \chi $ in
time $ t\in{\mathbb K} $, as a mapping $ \alpha_{t}\colon U_{1} \to U_{2} $; here $ U_{1,2} $ are open subsets of $ M $.
Consider $ U_{1,2} $ with restriction bihamiltonian structures. Then $ \alpha_{t} $ is an
isomorphism of bihamiltonian structures. \end{proposition}

\begin{proof} Indeed, the Hamiltonian vector field of any function w.r.t.~a
Poisson bracket $ \left\{,\right\} $ preserves $ \left\{,\right\} $. Thus the bracket $ \lambda'_{1}\left\{,\right\}_{1}+\lambda'_{2}\left\{,\right\}_{2} $ is
preserved by $ \chi $, thus by $ \alpha_{t} $. On the other hand, due to the condition on $ F $,
$ \chi $ is proportional to the Hamiltonian flow of $ F $ w.r.t.~any bracket
$ \lambda''_{1}\left\{,\right\}_{1}+\lambda''_{2}\left\{,\right\}_{2} $ (as far as $ \lambda''_{1}:\lambda''_{2}\not=\lambda_{1}:\lambda_{2} $), thus $ \chi $ preserves the
Poisson structure $ \lambda''_{1}\left\{,\right\}_{1}+\lambda''_{2}\left\{,\right\}_{2} $ as well. By linearity, it preserves
$ \lambda_{1}\left\{,\right\}_{1}+\lambda_{2}\left\{,\right\}_{2} $ too. \end{proof}

\begin{definition} For $ m\in M $ call $ v\in{\mathcal T}_{m}M $ a {\em biflow\/} vector if on a small open subset
$ U\subset M $, $ U\supset m $, the vector $ v $ can be represented as a value of $ \chi $ at $ m $; here $ \chi $ is
a vector field from Proposition~\ref{prop10.30}. \end{definition}

\begin{theorem} \label{th10.50}\myLabel{th10.50}\relax  Consider a micro-Kronecker bihamiltonian structure on a
manifold $ M $. Let $ m\in M $, let $ L $ be the leaf of action foliation on $ M $ which
passes through $ m $. If $ m'\in L $, then there are neighborhoods $ U $, $ U' $ of $ m $ and $ m' $
and a diffeomorphism $ \alpha\colon U \to U' $ which
\begin{enumerate}
\item
sends the restriction bihamiltonian structure on $ U $ to the
restriction bihamiltonian structure on $ U' $;
\item
for $ n\in U_{1} $ the point $ \alpha\left(n\right) $ is on the same leaf of action foliation as
$ n $.
\end{enumerate}
\end{theorem}

\begin{proof} Due to Proposition~\ref{prop10.30} it is enough to show that the
span of biflow vectors in $ {\mathcal T}_{m}M $ coincides with $ {\mathcal T}_{m}L $. Decrease $ M $ so that the
action foliation becomes a fibration $ \pi\colon M \to {\mathcal B} $ with a base $ {\mathcal B} $. Let $ b=\pi\left(m\right) $.

\begin{lemma} Consider the pencil $ \widetilde{{\mathcal P}}_{1,2}\colon {\mathcal T}_{b}^{*}{\mathcal B} \to {\mathcal T}_{m}L $ from the proof of
Proposition~\ref{prop8.40}. A vector $ v\in{\mathcal T}_{m}M $ is a biflow vector iff $ v $ can be
written as $ \left(\lambda_{1}'\widetilde{{\mathcal P}}_{1}+\lambda_{2}'\widetilde{{\mathcal P}}_{2}\right)\alpha $ with $ \left(\lambda_{1}\widetilde{{\mathcal P}}_{1}+\lambda_{2}\widetilde{{\mathcal P}}_{2}\right)\alpha=0 $, $ \lambda_{1},\lambda_{2},\lambda_{1}',\lambda_{2}'\in{\mathbb K} $,
$ \left(\lambda_{1},\lambda_{2}\right)\not=\left(0,0\right) $, $ \alpha\in{\mathcal T}_{b}^{*}{\mathcal B} $. \end{lemma}

\begin{proof} If $ f $ is a function on $ {\mathcal B} $, and $ \chi $ is the Hamiltonian flow of $ f\circ\pi $
w.r.t.~the Poisson bracket $ \lambda'_{1}\left\{,\right\}_{1}+\lambda'_{2}\left\{,\right\}_{2} $, then the value of $ \chi $ at $ m\in M $
coincides with $ \left(\lambda_{1}'\widetilde{{\mathcal P}}_{1}+\lambda_{2}'\widetilde{{\mathcal P}}_{2}\right)\left(df|_{b}\right) $. Since locally any Casimir function
for $ \lambda_{1}\left\{,\right\}_{1}+\lambda_{2}\left\{,\right\}_{2} $ can be written as $ F=f\circ\pi $ for an appropriate $ f $ such that
$ \left(\lambda_{1}\widetilde{{\mathcal P}}_{1}+\lambda_{2}\widetilde{{\mathcal P}}_{2}\right)\left(df|_{b}\right)=0 $, this implies the ``only if'' part of the lemma.

On the other hand, if $ \left(\lambda_{1}\widetilde{{\mathcal P}}_{1}+\lambda_{2}\widetilde{{\mathcal P}}_{2}\right)\alpha=0 $, then $ \alpha $ is normal to the leaf
of the integrating foliation $ {\mathcal F}_{\lambda_{1}:\lambda_{2}} $ on $ {\mathcal B} $. Decreasing $ {\mathcal B} $, one can find a
function $ f $ on $ {\mathcal B} $ such that $ f $ is constant on leaves of $ {\mathcal F}_{\lambda_{1}:\lambda_{2}} $, and $ df|_{b}=\alpha $.
This implies the ``if'' part of the lemma. \end{proof}

\begin{lemma} Consider a Kronecker relation in a vector space $ V $, and the
associated pencil $ {\mathcal P}_{1,2}\colon V \to V' $. Let $ \widetilde{V}' $ be the span of vectors $ v'\in V' $ of
the form $ \left(\lambda_{1}'{\mathcal P}_{1}+\lambda_{2}'{\mathcal P}_{2}\right)v $ with $ v\in V $ and $ \left(\lambda_{1}{\mathcal P}_{1}+\lambda_{2}{\mathcal P}_{2}\right)v=0 $, $ \lambda_{1},\lambda_{2},\lambda_{1}',\lambda_{2}'\in{\mathbb K} $,
$ \left(\lambda_{1},\lambda_{2}\right)\not=\left(0,0\right) $. Then $ \widetilde{V}'=V' $. \end{lemma}

\begin{proof} We may assume $ \lambda'_{1}=1 $, $ \lambda'_{2}=0 $. Since $ {\mathcal P}_{1}V=V' $, it is enough to show
that vectors $ v\in V $ such that $ \left(\lambda_{1}{\mathcal P}_{1}+\lambda_{2}{\mathcal P}_{2}\right)v=0 $ for an appropriate $ \lambda_{1},\lambda_{2}\in{\mathbb K} $ span
$ V $. But this is a corollary of Lemma~\ref{lm35.30}. \end{proof}

\begin{lemma} In conditions of Proposition~\ref{prop10.30} assume that the
bihamiltonian structure on $ M $ is micro-Kronecker, and $ {\mathcal F} $ is the action
foliation. Then for any $ m\in U_{1} $ points $ m $ and $ \alpha_{t}m $ are on the same leaf of $ {\mathcal F} $.
\end{lemma}

\begin{proof} This follows from the fact that biflow vectors are tangent to
leaves of action foliation. \end{proof}

This finishes the proof of Theorem~\ref{th10.50}. \end{proof}

\begin{theorem} \label{th10.60}\myLabel{th10.60}\relax  Consider a micro-Kronecker bihamiltonian structure on a
manifold $ M $ and an anti-involution $ \alpha $ of $ M $. Let $ Z $ be the submanifold formed
by the fixed points of $ \alpha $ of defect 0. Let $ U $ be the union of leaves of the
action foliation of $ M $ which intersect $ Z $. Then
\begin{enumerate}
\item
The subset $ U\subset M $ is open;
\item
If $ m\in U $, then Conjecture~\ref{con8.50} holds for $ m $ and $ M $.
\end{enumerate}
\end{theorem}

\begin{proof} This is an immediate corollary of Theorems~\ref{th9.50} and
~\ref{th10.50}. \end{proof}

Proposition~\ref{prop9.20} and Theorem~\ref{th10.60} show that in fact
Conjectures~\ref{con8.50} and~\ref{con9.45} are equivalent to each other. However,
the r\^oles of these conjecture are very dissimilar. As we will show in
Section~\ref{h11}, for some particular bihamiltonian structures of
mathematical physics Conjecture~\ref{con9.45} is easy to verify by an explicit
construction (see Theorem~\ref{th11.50}), thus for these structures Conjecture
~\ref{con8.50} {\em follows\/} from this construction.

On the other hand, as it was shown in \cite{GelZakhWeb}, in the case an
arbitrary bihamiltonian structure of rank 1 it is possible to prove
Conjecture~\ref{con8.50} using some ``hard'' cohomological statements. One can
expect that a similar approach may succeed in the case of higher rank as
well. But currently it is not clear how one could prove Conjecture~\ref{con9.45}
without a reference to Conjecture~\ref{con8.50} (or the calculation of
some cohomology which will immediately prove Conjecture~\ref{con8.50}).

This suggests that in some particular cases it is easier to directly
deduce the statement of Conjecture~\ref{con9.45}, but in general case
Conjecture~\ref{con8.50} should be easier to tackle.

\section{Method of argument translation }\label{h11}\myLabel{h11}\relax 

Consider a Lie algebra $ {\mathfrak g} $ and a $ 2 $-cocycle $ c_{2} $ of $ {\mathfrak g} $. As in Example
~\ref{ex002.45}, such a pair induces a bihamiltonian structure $ \left\{,\right\}_{1,2} $ on $ {\mathfrak g}^{*} $.
Moreover, if $ c_{2} $ is a coboundary of a $ 1 $-chain $ c_{1} $, one may consider $ c_{1} $ as
an element of $ {\mathfrak g}^{*} $. Obviously, the Poisson structure $ \left\{,\right\}_{1} $ is a Lie
derivative of the Poisson structure $ \left\{,\right\}_{2} $ in the direction of a parallel
translation of $ {\mathfrak g}^{*} $ in the direction of $ c_{1}\in{\mathfrak g}^{*} $, and $ \left\{,\right\}_{1} $ is
translation-invariant. Thus $ \left\{,\right\}_{2} +\lambda\left\{,\right\}_{1} $ is a parallel translation of $ \left\{,\right\}_{2} $
by $ \lambda c_{1} $. Due to this observation consideration of the pair $ \left\{\right\}_{1,2} $ when $ c_{2} $
is a coboundary is often called the {\em method of argument translation}.

Moreover, if $ {\mathfrak g} $ is semisimple, then $ c_{2} $ is automatically a coboundary.

\begin{definition} Given $ c_{1}\in{\mathfrak g}^{*} $, let $ \left\{,\right\}_{1,2} $ be the bihamiltonian structure on $ {\mathfrak g}^{*} $
constructed based on $ 2 $-coboundary $ dc_{1} $. Call this structure the {\em associated\/}
to $ c_{1}\in{\mathfrak g}^{*} $ structure. \end{definition}

\begin{definition} Given a Lie algebra $ {\mathfrak g} $ over $ {\mathbb K} $, let $ {\mathfrak g}_{{\mathbb C}}\buildrel{\text{def}}\over{=}{\mathfrak g} $ if $ {\mathbb K}={\mathbb C} $, and
$ {\mathfrak g}_{{\mathbb C}}\buildrel{\text{def}}\over{=}{\mathfrak g}\otimes{\mathbb C} $ if $ {\mathbb K}={\mathbb R} $. \end{definition}

\begin{lemma} Consider the Poisson structure $ \left\{,\right\}_{2} $ on $ {\mathfrak g}^{*} $ and $ \alpha\in{\mathfrak g}^{*} $. Let $ L $ be
the symplectic leaf of $ \left\{,\right\}_{2} $ through $ \alpha $. Then $ {\mathcal N}_{\alpha}^{*}L=\operatorname{Stab}_{\operatorname{ad}^{*}}\alpha $. \end{lemma}

The proof of this lemma is a direct calculation.

Recall description of the geometry of the set of regular elements in
a dual space to a Lie algebra.

\begin{definition} The {\em rank\/} $ \operatorname{rk}\left({\mathfrak g}\right) $ of a Lie algebra $ {\mathfrak g} $ is $ \min _{\alpha\in{\mathfrak g}_{{\mathbb C}}^{*}}\dim  \operatorname{Stab}_{\operatorname{Ad}^{*}}\alpha $.
An element $ \alpha\in{\mathfrak g}^{*} $ is {\em regular\/} if $ \dim  \operatorname{Stab}_{\operatorname{Ad}^{*}}\alpha=\operatorname{rk}\left({\mathfrak g}\right) $, and {\em irregular\/}
otherwise. \end{definition}

Obviously, regular elements form a (Zariski) open and dense subset
of $ {\mathfrak g}^{*} $.

\begin{definition} \label{def11.08}\myLabel{def11.08}\relax  Let $ \alpha\in{\mathfrak g}^{*} $, $ \beta\in{\mathfrak g}_{{\mathbb C}}^{*} $. Call $ \beta $ {\em compatible\/} with $ \alpha $ if $ \beta+\lambda\alpha $ is
regular for any $ \lambda\in{\mathbb C} $. Call a regular element $ \alpha\in{\mathfrak g}^{*} $ {\em strongly regular\/} if
there exists a compatible with $ \alpha $ element of $ {\mathfrak g}_{{\mathbb C}}^{*} $. \end{definition}

\begin{definition} Call $ {\mathfrak g} 2 $-{\em regular}, if the algebraic subvariety $ {\mathcal I}\subset{\mathfrak g}_{{\mathbb C}}^{*} $ of
irregular elements has codimension 2 or more. \end{definition}

\begin{proposition} \label{prop11.07}\myLabel{prop11.07}\relax  Suppose that there exists a strongly regular
element $ \alpha $ in $ {\mathfrak g}_{{\mathbb C}}^{*} $. Then $ {\mathfrak g} $ is $ 2 $-regular.

Suppose that $ {\mathfrak g} $ is $ 2 $-regular, and $ \alpha $ is a regular element of $ {\mathfrak g}^{*} $. Then
$ \alpha $ is strongly regular, and the set of elements of $ {\mathfrak g}^{*} $ which are compatible
with $ \alpha $ is non-empty and Zariski open. \end{proposition}

\begin{proof} Let $ {\mathcal I}\subset{\mathfrak g}_{{\mathbb C}}^{*} $ be the subvariety of irregular elements. Let $ \pi $ be the
projection of $ {\mathfrak g}_{{\mathbb C}}^{*} $ to $ {\mathfrak g}_{{\mathbb C}}^{*}/{\mathbb K}\alpha $. The existence of an $ \alpha $-compatible element is
equivalent to $ \pi{\mathcal I}\not=\pi{\mathfrak g}_{{\mathbb C}}^{*} $.

If $ \beta $ is regular, then any non-zero scalar multiple of $ \beta $ is also
regular. Thus one can consider a closed subvariety $ {\mathbb P}{\mathcal I} $ of irregular
elements in the projectivization $ {\mathbb P}{\mathfrak g}^{*} $ of $ {\mathfrak g}^{*} $. Given a strongly regular
element $ \alpha $ and a compatible element $ \beta $, one obtains a line $ l={\mathbb P}\left<\alpha,\beta\right> $ in
$ {\mathbb P}{\mathfrak g}^{*} $, and $ l\cap{\mathbb P}{\mathcal I}=\varnothing $. Clearly, any nearby line $ l' $ will also not intersect $ {\mathbb P}{\mathcal I} $.
Thus the set of elements $ \beta $ which are compatible with $ \alpha $ is Zariski open.
Thus the intersection of this set with $ {\mathfrak g}^{*}\subset{\mathfrak g}_{{\mathbb C}}^{*} $ is non-empty.

Since $ l\cap{\mathbb P}{\mathcal I}=\varnothing $, $ {\mathbb P}{\mathcal I} $ has codimension 2 or more, thus the same is true
for $ {\mathcal I} $. On the other hand, if $ {\mathcal I} $ has codimension 2 or more, then the
projection of $ {\mathcal I} $ to $ {\mathfrak g}^{*}/{\mathbb K}\alpha $ is not surjective; here $ \alpha\in{\mathfrak g}^{*} $ is arbitrary. This
implies that any regular element of $ {\mathfrak g}^{*} $ is strongly regular. \end{proof}

\begin{proposition} \label{prop11.10}\myLabel{prop11.10}\relax  Suppose that $ c_{1}\in{\mathfrak g}^{*} $ is strongly regular. Then there
is a dense open subset $ U\subset{\mathfrak g}^{*} $ such that the restriction on $ U $ of the pair
$ \left\{,\right\}_{1,2} $ associated to $ c_{1} $ is micro-Kronecker of rank $ \operatorname{rk}\left({\mathfrak g}\right) $. \end{proposition}

\begin{proof} By Proposition~\ref{prop11.07}, the set $ U $ of compatible with $ c_{1} $
elements of $ {\mathfrak g}_{{\mathbb C}}^{*} $ is Zariski open (thus dense). Show that $ U $ (or $ U\cap{\mathfrak g} $ in the
case $ {\mathbb K}={\mathbb R} $) satisfies the conditions of the proposition. One may assume
that $ {\mathbb K}={\mathbb C} $.

We need to show that for $ \beta\in U $ the symplectic leaf through $ \beta $ of
$ \lambda_{1}\left\{,\right\}_{1}+\lambda_{2}\left\{,\right\}_{2} $, $ \left(\lambda_{1},\lambda_{2}\right)\not=\left(0,0\right) $, has codimension $ \operatorname{rk}\left({\mathfrak g}\right) $. For $ \lambda_{2}=0 $ the normal
space to this leaf coincides with $ \operatorname{Stab}_{\operatorname{Ad}^{*}}c_{1} $, thus regularity of $ c_{1} $
implies the statement. Thus we may assume $ \lambda_{2}=1 $. The Poisson structure
$ \lambda\left\{,\right\}_{1}+\left\{,\right\}_{2} $ is a $ \lambda c_{1} $-translation of $ \left\{,\right\}_{2} $. If $ \beta\in U $, then $ \beta_{1}=\lambda c_{1}+\beta $ is
regular.

Thus it is enough consider $ \lambda=0 $. But the normal space to the leaf of
symplectic foliation through $ \beta_{1} $ is $ \operatorname{Stab}_{\operatorname{ad}^{*}}\beta_{1} $, which finishes the proof. \end{proof}

\begin{proposition} If $ {\mathfrak g} $ is reductive, then any regular element $ \alpha\in{\mathfrak g}^{*} $ is
strongly regular. \end{proposition}

\begin{proof} Indeed \cite{Ada69Lec}, irregular elements of a semisimple Lie
algebra form a Zariski closed subvariety of codimension 3. This implies
that the same statement for reductive algebras, thus the proposition. \end{proof}

The arguments above are not new, see \cite{Bol91Com,Pan98Sym}.

The last proposition cannot be inverted:

\begin{example} \label{ex11.22}\myLabel{ex11.22}\relax  (Proposed by V.~Serganova) For any vector space $ V $
there is a canonical symmetric pairing on $ V\oplus V^{*} $. For any Lie algebra $ {\mathfrak g} $
this pairing on $ {\mathfrak g}\oplus{\mathfrak g}^{*} $ is an invariant pairing on the Lie algebra
$ {\mathfrak G}={\mathfrak g}\ltimes \operatorname{ad}_{{\mathfrak g}}^{*} $; here $ \operatorname{ad}_{{\mathfrak g}}^{*} $ is the adjoint representation of $ {\mathfrak g} $ with
trivial structure of Lie algebra. This gives an isomorphism $ {\mathfrak G}\simeq{\mathfrak G}^{*} $ of
$ {\mathfrak G} $-modules.

This allows one to replace the study of dimensions of stabilizers of
elements of $ {\mathfrak G}^{*} $ by the study of dimensions of stabilizers of
elements of $ {\mathfrak G} $. If $ {\mathfrak g}={\mathfrak s}{\mathfrak l}_{2} $, it is easy to show that the set of
irregular elements coincides with the radical of $ {\mathfrak G} $, which has codimension
3. Later, in Example~\ref{ex95.75}, we will see that this statement on
codimension holds for other reductive algebras too. \end{example}

Now show how one can refine the description of bihamiltonian
structure on $ {\mathfrak g} $ by applying the general machinery of this paper.

\begin{definition} \label{def11.25}\myLabel{def11.25}\relax  A linear mapping $ \iota\colon {\mathfrak g} \to {\mathfrak g} $ is an
{\em antiinvolution\/} of $ {\mathfrak g} $ if $ \iota $ is an involution of a vector space, and
$ \left[\iota X,\iota Y\right]=-\iota\left[X,Y\right] $ for any $ X,Y\in{\mathfrak g} $. An anti-involution $ \iota $ is
{\em admissible}, if
\begin{enumerate}
\item
The irregular elements in the vector subspace $ \operatorname{Fix}\left(\iota^{*}\right)\subset{\mathfrak g}^{*} $ of fixed
points of $ \iota $ in $ {\mathfrak g}^{*} $ form a subvariety of $ \operatorname{Fix}\left(\iota^{*}\right) $ of codimension 2 or more.
\item
The subset $ {\mathcal U}\left(\iota\right)\subset\operatorname{Fix}\left(\iota^{*}\right) $ consisting of points $ \alpha\in\operatorname{Fix}\left(\iota^{*}\right) $ such that the
$ \operatorname{Ad}^{*} $-orbit of $ \alpha $ is transversal to $ \operatorname{Fix}\left(\iota^{*}\right) $ is not empty;
\end{enumerate}

Call $ \alpha\in{\mathfrak g}^{*} $ {\em admissible}, if $ \alpha $ is regular, and there is an admissible
anti-involution $ \iota $ such that $ \iota\alpha=\alpha $. \end{definition}

\begin{remark} \label{rem11.28}\myLabel{rem11.28}\relax  Obviously, admissible elements exists only in $ 2 $-regular
Lie algebras, and are strongly regular. Clearly, for an admissible
anti-involution $ \iota $ the set $ {\mathcal U}\left(\iota\right) $ is Zariski open in $ \operatorname{Fix}\left(\iota^{*}\right) $. Moreover, if $ \alpha $
is admissible, and $ \iota $ is the corresponding admissible anti-involution, then
elements $ \beta\in\operatorname{Fix}\left(\iota^{*}\right) $ which are compatible with $ \alpha $ form a non-empty Zariski
open subset of $ \operatorname{Fix}\left(\iota^{*}\right) $. \end{remark}

\begin{theorem} \label{th11.40}\myLabel{th11.40}\relax  Suppose that $ c_{1}\in{\mathfrak g}^{*} $ is admissible. Consider the
bihamiltonian structure $ \left\{,\right\}_{1,2} $ associated to $ c_{1} $. Then there is an open
subset $ M\subset{\mathfrak g}^{*} $ such that for any $ m\in M $ Conjecture~\ref{con8.50} holds. \end{theorem}

\begin{proof} Consider an admissible anti-involution $ \iota $ such that $ \iota^{*}c_{1}=c_{1} $. Let
$ U $ be the Zariski open subset of $ {\mathfrak g}_{{\mathbb C}}^{*} $ where $ \left\{,\right\}_{1,2} $ is micro-Kronecker.
Since $ \iota $ is an anti-involution of $ {\mathfrak g} $, $ \iota^{*} $ is an anti-involution of the
Poisson structure $ \left\{,\right\}_{2} $ on $ {\mathfrak g}^{*} $. Since $ \iota^{*} $ preserves $ c_{1} $, and $ \left\{,\right\}_{1} $ is the
derivative of $ \left\{,\right\}_{2} $ w.r.t.~translations in the direction of $ c_{1} $, $ \iota^{*} $ is an
anti-involution of $ \left\{,\right\}_{1} $ as well. Thus $ U $ is $ \iota^{*} $-invariant, and $ \iota^{*}|_{U} $ is an
anti-involution of a micro-Kronecker bihamiltonian structure.

Note that the same arguments as in the proof of Proposition
~\ref{prop11.07} show that $ U\cap\operatorname{Fix}\left(\iota^{*}\right)\not=\varnothing $. Let $ \widetilde{U}=U\cap{\mathcal U}\left(\iota^{*}\right) $. This is a non-empty
Zariski open subset of $ \operatorname{Fix}\left(\iota^{*}\right) $. Let $ \beta\in\widetilde{U} $. The principal step is to show
that the defect of $ \iota^{*} $ at $ \beta $ is 0.

Let $ U_{1} $ be a neighborhood of $ \beta $ in $ U $ such that the action foliation $ {\mathcal F} $
of the bihamiltonian structure becomes a fibration $ \pi\colon U_{1} \to {\mathcal B} $. It is
enough to show that for any function $ \varphi $ on $ {\mathcal B} $ the function $ \varphi\circ\pi $ on $ U_{1} $ is
preserved by $ \iota^{*} $. In turn, it is enough to do the same for a large enough
collection of functions $ \varphi $ on $ {\mathcal B} $. Take as such collection functions $ \varphi $ which
are constant on fibers of the integrating foliation $ {\mathcal F}_{\lambda} $ on $ {\mathcal B} $, for each one
of $ \lambda\in\Lambda $ (for a sufficiently large $ \Lambda\subset{\mathbb K}{\mathbb P}^{1} $). We may suppose $ \Lambda\subset{\mathbb K} $.

Obviously, $ \varphi $ is constant on fibers of $ {\mathcal F}_{\lambda} $ iff $ \varphi\circ\pi $
is constant on fibers of the symplectic foliation $ \widetilde{{\mathcal F}}_{\lambda} $ of the Poisson
structure $ \lambda\left\{,\right\}_{1}+\left\{,\right\}_{2} $. Thus it is enough to show that $ \iota^{*} $ preserves such
functions. Since $ \lambda\left\{,\right\}_{1}+\left\{,\right\}_{2} $ is the result of translation of $ \left\{,\right\}_{2} $ by $ \lambda c_{1} $,
it is enough to show this for $ \lambda=0 $ and $ \beta+\lambda c_{1} $, $ \lambda\in\Lambda $, taken instead of $ \beta $.

However, $ \beta'=\beta+\lambda c_{1} $ is in $ \widetilde{U} $ for $ \lambda $ in an open subset of $ {\mathbb K} $, thus we can
restrict our attention to a given $ \beta'\in\widetilde{U} $ and $ \lambda=0 $. Since $ \beta'\in{\mathcal U}\left(\iota^{*}\right) $, any
$ \operatorname{Ad}^{*} $-orbit which passes near $ \beta' $ intersects $ \operatorname{Fix}\left(\iota^{*}\right) $. Since $ \iota^{*} $ sends an
$ \operatorname{Ad}^{*} $-orbit to an $ \operatorname{Ad}^{*} $-orbit, this implies that $ \iota^{*} $ preserves any $ \operatorname{Ad}^{*} $-orbit
which passes near $ \beta' $. Decrease $ U_{1} $ so that $ \operatorname{Ad}^{*} $-orbit of any $ \gamma\in U_{1} $
intersects $ \operatorname{Fix}\left(\iota^{*}\right) $. However, symplectic leaves of $ \left\{,\right\}_{2} $ coincide with
$ \operatorname{Ad}^{*} $-orbits in $ {\mathfrak g}^{*} $. Thus on a neighborhood of $ \beta' $ any function which is
constant on symplectic leaves of $ \left\{,\right\}_{2} $ is $ \iota^{*} $-invariant.

This implies that the defect of $ \iota^{*} $ near $ \beta $ is indeed 0, and we are in
conditions of Theorem~\ref{th10.60}. This implies the theorem for $ M $ being the
union of fibers of the action foliation which intersect $ \widetilde{U} $. \end{proof}

\begin{theorem} \label{th11.50}\myLabel{th11.50}\relax  Suppose that $ {\mathfrak g} $ is semisimple, and $ c_{1}\in{\mathfrak g}^{*} $ is regular and
semisimple. Then there is a dense open subset $ U\subset{\mathfrak g}^{*} $ such that the
restriction on $ U $ of the pair $ \left\{,\right\}_{1,2} $ associated to $ c_{1} $ satisfies Conjecture
~\ref{con8.50}. \end{theorem}

\begin{proof} Identify $ {\mathfrak g}^{*} $ with $ {\mathfrak g} $ using the Killing form. Due to Proposition
~\ref{prop95.10}, it is enough to consider the case $ {\mathbb K}={\mathbb C} $.

\begin{definition} Consider a semisimple Lie algebra $ {\mathfrak g} $. Define the {\em Cartan
antiinvolution\/} $ \iota $ by its action on standard generators $ e_{i},f_{i},h_{i} $,
$ i=1,\dots ,r $:
\begin{equation}
\iota\left(e_{i}\right)=f_{i}\text{, }\iota\left(f_{i}\right)=e_{i}\text{, }\iota\left(h_{i}\right)=h_{i}.
\notag\end{equation}
\end{definition}

\begin{lemma} \label{lm11.60}\myLabel{lm11.60}\relax  The Cartan anti-involution of a semisimple Lie algebra is
admissible. \end{lemma}

\begin{proof} Identification of $ {\mathfrak g} $ with $ {\mathfrak g}^{*} $ allows one to consider $ \iota $ instead of
$ \iota^{*} $. Any anti-involution sends an $ \operatorname{Ad} $-orbit to an $ \operatorname{Ad} $-orbit. Since $ \operatorname{Fix}\left(\iota\right)\supset{\mathfrak h} $,
and an orbit of a regular element of $ {\mathfrak h} $ is transversal to $ {\mathfrak h} $, $ \iota $ satisfies
the second condition of Definition~\ref{def11.25}. Thus to prove the lemma it
is enough to show that irregular elements in $ \operatorname{Fix}\left(\iota\right) $ form a subvariety of
codimension 2 or more. (Note that this statement is {\em not\/} true if one
substitutes $ {\mathfrak h} $ instead of $ \operatorname{Fix}\left(\iota\right)! $)

Due to homogeneity of the set of regular elements, a translation to
algebraic geometry in a projective space show that it is enough to prove
this statement for an arbitrary vector subspace in $ \operatorname{Fix}\left(\iota\right) $ taken instead
of $ \operatorname{Fix}\left(\iota\right) $ (see Proposition~\ref{prop11.07}). Recall that \cite{Bour75Lie78}:

\begin{lemma} \label{lm11.82}\myLabel{lm11.82}\relax  There are numbers $ a_{i}\not=0 $, $ b_{i}\not=0 $, $ c_{i} $, $ i=1,\dots ,r $, such that
for the elements $ E=\sum a_{i}e_{i} $, $ F=\sum b_{i}f_{i} $, $ H=\sum c_{i}h_{i} $ the vector subspace
$ V=\left<E,F,H\right>\subset{\mathfrak g} $ is a Lie subalgebra isomorphic to $ {\mathfrak s}{\mathfrak l}_{2} $. The adjoint action
of $ V $ on $ {\mathfrak g} $ is a direct sum of $ r $ odd-dimensional irreducible
representations of $ {\mathfrak s}{\mathfrak l}_{2} $. \end{lemma}

Since the action of any non-zero element of $ {\mathfrak s}{\mathfrak l}_{2} $ in an
odd-dimensional irreducible representation has $ 1 $-dimensional null-space,
this shows that the stabilizer of any non-zero point of $ V $ has dimension
$ r $. But $ \operatorname{rk}{\mathfrak g}=r $, thus all the non-zero elements of $ V $ are regular. Moreover,
conjugating $ V $ with elements of the Lie group $ \exp \left({\mathfrak h}\right) $ of $ {\mathfrak h} $, one may assume
that $ a_{i}=b_{i} $, $ i=1,\dots ,r $. Thus the subspace $ V $ is $ \iota $-invariant, let $ V_{0} $ be the
$ 2 $-dimensional invariant subspace $ \left<H,E+F\right> $ of $ \iota|_{V} $.

Since $ V_{0}\subset\operatorname{Fix}\left(\iota\right) $ intersects irregular elements on $ \left\{0\right\} $, which is a
subvariety of codimension 2, irregular elements in $ \operatorname{Fix}\left(\iota\right) $ form a
submanifold of codimension 2 or more. This finishes the proof of Lemma
~\ref{lm11.60}. {}\end{proof}

Lemma~\ref{lm11.60} implies that any regular $ c_{1}\in{\mathfrak g}^{*} $ which is on an
$ \operatorname{Ad}^{*} $-orbit of $ \operatorname{Fix}\left(\iota^{*}\right) $ is admissible. But any regular semisimple element of
$ {\mathfrak g} $ is conjugate to an element of $ {\mathfrak h} $, it is on an $ \operatorname{Ad} $-orbit of an $ \iota $-invariant
element! Translating from $ {\mathfrak g} $ to $ {\mathfrak g}^{*} $, any regular semisimple $ c_{1}\in{\mathfrak g}^{*} $ is
admissible. Hence there is an open subset $ U\subset{\mathfrak g}^{*} $ on which Conjecture
~\ref{con8.50} holds.

Let us show that $ U $ can be taken Zariski open, thus dense. Let $ U_{0} $ be
the Zariski open subset where foliation $ {\mathcal F} $ makes sense. Recall that $ U $ is
the union of leaves of $ {\mathcal F} $ which intersect $ \operatorname{Fix}\left(\iota\right) $. Show that $ U $ contains a
Zariski open subset.

There is a Zariski open subset $ U_{1} $ of $ {\mathfrak g} $ such that $ U_{1} $ is $ \operatorname{Ad} $-invariant,
and $ \operatorname{Ad} $-invariant polynomials distinguish\footnote{In the same sense as in Remark~\ref{rem002.20}.} $ \operatorname{Ad} $-orbits in $ U_{1} $. On the other
hand, locally the fibers of the action foliation $ {\mathcal F} $ are intersections of a
finite number of $ \lambda c_{1} $-shifted $ \operatorname{Ad} $-orbits, $ \lambda\in\Lambda $, $ \operatorname{card}\left(\Lambda\right)<\infty $. Taking a large
enough finite collection $ p_{i} $, $ i\in I $, of invariant polynomials on $ {\mathfrak g} $, we see
that fibers of $ {\mathcal F} $ coincide with connected components of level sets of
$ p_{i,\lambda} $, $ i\in I $, $ \lambda\in\Lambda $; here $ p_{i,\lambda}\left(X\right)\buildrel{\text{def}}\over{=}p_{i}\left(X+\lambda c_{1}\right) $. Let $ U_{2}=U_{0}\cap\bigcap_{\lambda\in\Lambda}\left(U_{1}-\lambda c_{1}\right) $; here
$ U_{1}-\lambda c_{1} $ is the parallel translation of $ U_{1} $ by $ -\lambda c_{1} $. Obviously, $ U_{2} $ is
Zariski open and is a union of fibers of foliation $ {\mathcal F} $.

Let $ \Pi $ be the polynomial mapping of $ {\mathfrak g} $ to $ {\mathbb C}^{N} $, $ N=\operatorname{card}\left(I\right)\operatorname{card}\left(\Lambda\right) $, with
components $ p_{i,\lambda} $. Let $ Y=\overline{\operatorname{Im}\Pi} $, clearly $ \Pi U_{2} $ contains a Zariski open subset
$ Y_{0} $ of $ Y $. Decreasing $ Y_{0} $, one may assume that $ \Pi|_{\Pi^{-1}Y_{0}} $ is a submersion to
$ Y_{0} $. Since $ \Pi $ is constant on leaves of $ {\mathcal F} $, $ \Pi\operatorname{Fix}\left(\iota\right)=\Pi U $. Since $ U $ contains an
open subset of $ {\mathfrak g} $, $ \Pi\operatorname{Fix}\left(\iota\right) $ contains an open subset of $ Y $, thus a Zariski
open subset of $ Y $. Thus one can assume that $ Y_{0}\subset\Pi\operatorname{Fix}\left(\iota\right) $. Let $ U_{3}=\Pi^{-1}Y_{0} $,
$ Z=U_{3}\cap\operatorname{Fix}\left(\iota\right) $.

We obtain the following mappings: $ Z \buildrel{i}\over{\hookrightarrow} U_{3} \xrightarrow[]{\Pi} Y_{0} $; here $ U_{3} $ is
Zariski open in $ {\mathfrak g} $, $ Z $ is Zariski closed in $ U_{3} $, $ \Pi $ is (a restriction of) a
polynomial mapping which is a submersion onto $ Y_{0} $, and $ \Pi\circ i $ is surjective.
Let $ U_{4} $ be the union of connected components of fibers of $ \Pi $ which
intersect $ Z $. Since $ U_{4}\subset U $, it is enough to show that $ U_{4} $ coincides with $ U_{3} $
(which is a union of fibers of $ \pi $ which intersect $ Z $).

Since $ \Pi $ is a submersion, the number of connected components of
fibers of $ \Pi $ can only jump up during specialization. This implies that $ U_{4} $
is open in $ U_{3} $. Decreasing $ Y_{0} $, we may assume that the number of connected
components of the fiber of $ \Pi $ over $ z $ does not depend on $ z\in Y_{0} $. This implies
that $ U_{4} $ is closed in $ U_{3} $. Since $ U_{3} $ is connected (as a Zariski open subset
of a vector space), $ U_{3}=U_{4} $. This finishes the proof of Theorem~\ref{th11.50}. \end{proof}

\begin{amplification} In Theorem~\ref{th11.50} one can take $ {\mathfrak g} $ being
reductive, and drop the condition of semisimplicity on $ c_{1} $. \end{amplification}

\begin{proof} As the proof of Theorem~\ref{th11.50} shows, it is enough to show
that

\begin{proposition} \label{prop11.75}\myLabel{prop11.75}\relax  Consider a reductive Lie algebra $ {\mathfrak g} $ and its Cartan
anti-involution $ \iota $. Then any element $ X $ of $ {\mathfrak g} $ is on $ \operatorname{Ad} $-orbit of $ \iota $-invariant
element. \end{proposition}

\begin{proof}[Proof (V.~Serganova) ] First of all, the statement for reductive
algebras follows momentarily from the case of semisimple Lie algebras.
For semisimple elements the statement is obvious (since $ {\mathfrak h} $ is
$ \iota $-invariant). Consider the case when $ X\in{\mathfrak g} $ is nilpotent.

One may assume that $ X\not=0 $. If $ {\mathfrak g}={\mathfrak s}{\mathfrak l}_{2} $, then any nilpotent element is
$ \operatorname{SL}_{2} $-conjugate to $ \left(
\begin{matrix}
1 & i \\ i & -1
\end{matrix}
\right) $ which is symmetric, thus $ \iota $-invariant. Thus to
prove the proposition for the case of nilpotent $ X $ is enough to show

\begin{lemma} \label{lm11.78}\myLabel{lm11.78}\relax  For any nilpotent $ X\in{\mathfrak g} $ with a reductive Lie algebra $ {\mathfrak g} $
there is a subalgebra $ {\mathfrak g}_{0}\subset{\mathfrak g} $ such that $ {\mathfrak g}_{0}\simeq{\mathfrak s}{\mathfrak l}_{2} $, $ \iota{\mathfrak g}_{0}={\mathfrak g}_{0} $, $ \iota|_{{\mathfrak g}_{0}} $ is the Cartan
anti-involution of $ {\mathfrak s}{\mathfrak l}_{2} $, and $ X $ is conjugate to $ X'\in{\mathfrak g}_{0} $. Here $ \iota $ is the Cartan
anti-involution of $ {\mathfrak g} $. \end{lemma}

\begin{proof} It is enough to prove this for a semisimple $ {\mathfrak g} $. The
classification of nilpotent elements up to conjugation is well-known
(\cite{Dyn52Max}, or it can be deduced from \cite{Bour75Lie78}):

\begin{lemma} For any nilpotent element $ X\in{\mathfrak g} $, $ X\not=0 $, there is a reductive
subalgebra $ \widetilde{{\mathfrak g}} $ of $ {\mathfrak g} $ with $ \operatorname{rk}\left(\widetilde{{\mathfrak g}}\right)=\operatorname{rk}\left({\mathfrak g}\right) $ and a Cartan set of generators $ \widetilde{e}_{i} $, $ \widetilde{f}_{i} $,
$ i=1,\dots ,\widetilde{r} $, $ \widetilde{h}_{j} $, $ j=1,\dots ,\operatorname{rk}\left(\widetilde{{\mathfrak g}}\right) $, of $ \widetilde{{\mathfrak g}} $ such that $ X=\sum_{i=1}^{\widetilde{r}}\widetilde{e}_{i} $. Here $ \widetilde{r}=\operatorname{rk}\left(\widetilde{{\mathfrak g}}\right)-\dim 
Z\left(\widetilde{{\mathfrak g}}\right) $. \end{lemma}

The Cartan subalgebra of $ \widetilde{{\mathfrak g}} $ is a Cartan subalgebra of $ {\mathfrak g} $, thus after
conjugation one can assume that $ \widetilde{h}_{j} $ generate $ {\mathfrak h}\subset{\mathfrak g} $. Then $ \iota\widetilde{h}_{j}=\widetilde{h}_{j} $ for any
$ j=1,\dots ,\operatorname{rk}\left(\widetilde{{\mathfrak g}}\right) $. Let $ a_{ij} $ be coefficients in relations $ \left[\widetilde{h}_{j},\widetilde{e}_{i}\right]=a_{ij}\widetilde{e}_{i} $,
$ \left[\widetilde{h}_{j},\widetilde{f}_{i}\right]=-a_{ij}\widetilde{f}_{i} $. Given $ a_{ij} $ and $ \widetilde{h}_{i} $, these relations determine $ \widetilde{e}_{i} $, $ \widetilde{f}_{i} $
uniquely up to proportionality. This implies that $ \iota\widetilde{e}_{i} $ is proportional to
$ \widetilde{f}_{i} $. Thus $ \iota\widetilde{{\mathfrak g}}=\widetilde{{\mathfrak g}} $, moreover, after a rescaling $ \iota|_{\widetilde{{\mathfrak g}}} $ may be supposed to be the
Cartan anti-involution of $ \widetilde{{\mathfrak g}} $.

Substituting $ \widetilde{{\mathfrak g}} $ instead of $ {\mathfrak g} $, it follows that it is enough to prove
the statement of Lemma~\ref{lm11.78} for $ X=\sum_{i=1}^{\operatorname{rk}\left({\mathfrak g}\right)}e_{i} $. In this case Lemma
~\ref{lm11.82} implies that $ X $ is $ \exp \left({\mathfrak h}\right) $-conjugate to the element $ E $ of Lemma
~\ref{lm11.82}. Moreover, doing another $ \exp \left({\mathfrak h}\right) $-conjugation one can ensure that
the vector space $ \left<E,F,H\right> $ of Lemma~\ref{lm11.82} is $ \iota $-invariant. \end{proof}

Now the proposition is proven for semisimple and for nilpotent
elements $ X\in{\mathfrak g} $. For an arbitrary $ X\in{\mathfrak g} $, there is unique representation
$ X=X_{ss}+X_{\text{nil}} $ as a sum of commuting semisimple and nilpotent elements. Doing
conjugation, we may assume $ X_{ss}\in{\mathfrak h} $. Let $ {\mathfrak g}_{0}=\operatorname{Stab}_{\operatorname{ad}}X_{ss} $. Obviously, $ \iota{\mathfrak g}_{0}={\mathfrak g}_{0} $, $ {\mathfrak g}_{0} $
is a reductive Lie algebra, $ {\mathfrak g}_{0}\supset{\mathfrak h} $, and $ \iota|_{{\mathfrak g}_{0}} $ is the Cartan anti-involution
of $ {\mathfrak g}_{0} $. Since $ X_{\text{nil}}\in{\mathfrak g}_{0} $, we know that $ X_{\text{nil}} $ is $ G_{0} $-conjugate to an element of
$ \operatorname{Fix}\left(\iota|_{{\mathfrak g}_{0}}\right) $; here $ G_{0} $ is the $ \operatorname{Ad} $-group of $ {\mathfrak g}_{0} $. Since $ X_{ss} $ is $ G_{0} $-invariant and
$ \iota $-invariant, the conjugation by $ g\in G_{0} $ above sends $ X_{ss}+X_{\text{nil}} $ to
$ \operatorname{Fix}\left(\iota|_{{\mathfrak g}_{0}}\right)\subset\operatorname{Fix}\left(\iota\right) $. This finishes the proof of the proposition. \end{proof}

This finishes the proof of the amplification. \end{proof}

Finally, one can apply the accumulated information to prove a
generalization of one of conjectures of \cite{GelZakh99Web}:

\begin{theorem} \label{th11.70}\myLabel{th11.70}\relax  Suppose that $ {\mathfrak g} $ is reductive, and $ c_{1}\in{\mathfrak g}^{*} $ is regular. Then
there is a dense open subset $ U\subset{\mathfrak g}^{*} $ such that the restriction on $ U $ of the
pair $ \left\{,\right\}_{1,2} $ is flat\footnote{Flat bihamiltonian structures were introduced in Definition~\ref{def00.26}.}. \end{theorem}

\begin{proof} Due to Theorem~\ref{th11.50}, it is enough to show that on a dense
open subset of $ {\mathfrak g}^{*} $ the Kronecker web which is a (local) base of the action
foliation is {\em flat\/}\footnote{I.e., locally isomorphic to a translation-invariant structure.}. This was proven in \cite{Pan99Ver} (in less generality);
here we reproduce a more general form of these arguments:

\begin{lemma} Consider a Lie algebra $ {\mathfrak g} $, $ \operatorname{rk}{\mathfrak g}=r $, and $ \operatorname{Ad}^{*} $-invariant polynomials
$ p_{1},\dots ,p_{r} $ on $ {\mathfrak g}^{*} $. Let $ U\subset{\mathfrak g}^{*} $ consists of points $ \alpha\in{\mathfrak g}^{*} $ such that
$ dp_{1}|_{\alpha},\dots ,dp_{r}|_{\alpha} $ are linearly independent. Suppose that $ {\mathfrak g} $ is $ 2 $-regular,
$ U\not=\varnothing $, $ \dim {\mathfrak g} = 2\sum_{i=1}^{r}\deg  p_{i} + r $, and that $ c_{1}\in{\mathfrak g}^{*} $ is regular. Consider the
Kronecker structure on the local base $ {\mathcal B} $ of the action foliation $ {\mathcal F} $ of the
bihamiltonian structure on $ {\mathfrak g}^{*} $ associated to $ c_{1} $. This Kronecker structure
is flat on an open subset. \end{lemma}

\begin{proof} Argue as in the end of the proof of Theorem~\ref{th11.50}. On a
Zariski open subset $ U_{0} $ of $ {\mathfrak g}^{*} $ the polynomials $ p_{i} $, $ i=1,\dots ,r $, locally
distinguish $ \operatorname{Ad}^{*} $-orbits (thus symplectic leaves of $ \left\{,\right\}_{2} $) on $ {\mathfrak g}^{*}. $Thus on a
Zariski open subset $ U_{1} $ of $ {\mathfrak g}^{*} $ the polynomials $ p_{i,\lambda} $, $ i=1,\dots ,r $, $ \lambda\in{\mathbb K} $,
locally distinguish leaves of $ {\mathcal F} $ (here $ p_{i,\lambda}\left(\alpha\right)\buildrel{\text{def}}\over{=}p_{i}\left(\alpha+\lambda c_{1}\right) $). Associate to
a given $ \alpha\in{\mathfrak g}^{*} $ the coefficients $ a_{ij} $, $ i=1,\dots ,r $, $ j=0,\dots ,\deg  p_{i} $, of
polynomials $ p_{i}\left(\alpha+\lambda c_{1}\right) $ in $ \lambda $. This is a polynomial mapping $ a\colon {\mathfrak g}^{*} \to {\mathbb K}^{N} $,
$ N=\sum_{i=1}^{r}\deg  p_{i} + r $. We conclude that on $ U_{1} $ connected components of fibers
of $ a|_{U_{0}} $ coincide with leaves of $ {\mathcal F} $.

But the leaves of $ {\mathcal F} $ have codimension $ \frac{\dim {\mathfrak g}+r}{2} $, thus in the
conditions of the lemma the mapping $ a|_{U_{0}} $ is a submersion, and leaves of $ {\mathcal F} $
are connected components of fibers of $ a $. Thus on $ U_{1} $ the manifold $ a\left(U_{1}\right) $
may be considered as a base $ {\mathcal B} $ of the action foliation. Describe the
structure of Kronecker web on $ a\left(U_{1}\right) $.

Fix $ \lambda\in{\mathbb K} $. The symplectic leaves of $ \lambda\left\{,\right\}_{1}+\left\{,\right\}_{2} $ are $ -\lambda c_{1} $-translations
of $ \operatorname{Ad}^{*} $-orbits, thus on $ U_{1} $ they coincide with level sets of $ p_{i,\lambda} $. Thus the
projections of these leaves to $ {\mathcal B} $ may be described by equations
$ \sum_{j}a_{ij}\lambda^{j}=C_{i} $, $ i=1,\dots ,r $; here $ a_{ij} $ are coordinates on $ {\mathbb K}^{N} $. Thus fibers of
integrating foliations $ {\mathcal F}_{\lambda} $ on $ {\mathcal B} $ are parallel planes in $ {\mathbb K}^{N} $, hence the
Kronecker web structure on $ {\mathcal B} $ is translation-invariant. \end{proof}

To finish the proof of the theorem, it is enough to recall
\cite{Bour75Lie78} that reductive Lie algebras satisfy the lemma. \end{proof}

\begin{corollary} \label{cor11.80}\myLabel{cor11.80}\relax  In conditions of Theorem~\ref{th11.70} there is an open
dense subset $ U\subset{\mathfrak g}^{*} $ such that the restriction on $ U $ of the bihamiltonian
structure is locally isomorphic to a direct product of several copies of
the structure of Example~\ref{ex00.30}. \end{corollary}

\begin{proof} It is enough to show that a micro-Kronecker
translation-invariant bihamiltonian structure can be represented as a
product of structures of Example~\ref{ex00.30}. This a direct corollary of
Theorem~\ref{th4.35}. \end{proof}

Sum up conditions under which the statements about flatness can be
achieved (at least in a weak form).

\begin{definition} Consider an open subset $ U\subset V $ in a vector space $ V $, a foliation
$ {\mathcal F} $ of codimension $ r $ on $ U $, and an involution $ \psi $ of $ U $. Call $ {\mathcal F} $ {\em compatible\/} with
$ \psi $ if
\begin{enumerate}
\item
The leaves of $ {\mathcal F} $ are common level sets of polynomials $ p_{1},\dots ,p_{r} $ on $ U $;
\item
Differentials $ dp_{i} $ are linearly independent at any point of $ U $;
\item
$ 2\sum_{i=1}^{r}\deg  p_{i}+r=\dim  V $;
\item
The submanifold $ \operatorname{Fix}\left(\psi\right) $ of fixed points of $ \psi $ in $ U $ is nonempty and is
transversal to leaves of $ {\mathcal F} $.
\end{enumerate}
\end{definition}

The arguments we had so far lead to

\begin{corollary} \label{cor95.30}\myLabel{cor95.30}\relax  Consider an anti-involution $ \iota $ of a $ 2 $-regular Lie
algebra $ {\mathfrak g} $. Suppose that there is an open subset $ U\subset{\mathfrak g}^{*} $ such that the
$ \operatorname{Ad}^{*} $-orbits in $ U $ form a foliation which is compatible with the
involution $ \iota^{*} $ of $ {\mathfrak g}^{*} $. Then there is an open subset $ U_{1}\subset{\mathfrak g}^{*} $, and for
any $ c_{1}\in U_{1} $ there is an open subset $ U_{2}\left(c_{1}\right)\subset{\mathfrak g}^{*} $ on which the bihamiltonian
structure associated to $ c_{1} $ is flat. \end{corollary}

This reduces the question of flatness to a question on geometry of
$ \operatorname{Ad}^{*} $-action on $ {\mathfrak g}^{*} $. Let us sketch roughly how to construct new Lie algebras
which satisfy the conditions of the corollary.

First of all, if $ {\mathfrak g} $ satisfies these condition, then $ {\mathfrak g}^{n} $ satisfies
these conditions too. The subset $ U $ should be replaced by $ U^{n} $, and $ \iota $ by $ \iota^{n} $.
More generally, if $ Z $ is a set, or a topological space, or a manifold, or
a variety, or a scheme, then a following variation is possible. Denote by
$ {\mathfrak g}^{Z} $ the set (or topological space, etc) of mappings from $ Z $ to $ {\mathfrak g} $. It has a
natural Lie algebra structure, and it acts on $ \left({\mathfrak g}^{*}\right)^{Z} $. Moreover, $ \left({\mathfrak g}^{*}\right)^{Z} $ has
an involution $ \iota^{Z} $, and in many situations $ U^{Z} $ makes sense\footnote{If $ Z $ is an affine line, and $ U $ has many ``holes'' (complement to $ U $ is
large), then the set of algebraic mappings to $ U $ can consist of constant
mappings only. We want to avoid such a situation.} as well.

Suppose that $ U^{Z} $ makes sense, and $ G^{Z} $ makes sense too (here $ G $ is a Lie
group with the Lie algebra $ {\mathfrak g} $). Then the orbits of action of $ G^{Z} $ on $ \left({\mathfrak g}^{*}\right)^{Z} $
can be described in the same way as in the case of $ {\mathfrak g}^{n} $: given a linear
functional $ \varepsilon $ on functions on $ Z $ and an invariant polynomial $ p $ on $ U $, define
a $ G^{Z} $-invariant function $ p_{\varepsilon} $ on $ U^{Z} $ by $ p_{\varepsilon}\left(f\right)=\left< \varepsilon,p\circ f \right> $; here $ f\colon Z \to U $. If
$ Z $ is ``small''\footnote{It looks like it is enough to require that any \'etale covering of $ Z $ can
be refined to a usual covering.} intersections of $ G^{Z} $-orbits with $ U^{Z} $ coincide with common
level sets of functions $ p_{\varepsilon} $ on $ U^{Z} $. Thus all the needed qualitative
properties of orbits of $ G^{Z} $ on $ U^{Z} $ hold in this case if they hold for
the action of $ {\mathfrak g} $ on $ U $.

The property $ 2\sum_{i=1}^{r}\deg  p_{i}+r=\dim  V $ makes no sense if $ {\mathfrak g}^{Z} $ is
infinite-dimensional, but holds if $ {\mathfrak g}^{Z} $ is finite-dimensional. We conclude:

\begin{corollary} Suppose that a Lie algebra $ {\mathfrak g} $ with an anti-involution $ \iota $
satisfies conditions of Corollary~\ref{cor95.30}. If $ A $ is a
finite-dimensional commutative algebra, and $ Z=\operatorname{Spec} $ A, then the orbits of
action of $ G^{Z} $ on $ \left({\mathfrak g}^{*}\right)^{Z} $ are compatible with $ \iota^{Z} $ on an appropriate open
subset of $ G^{Z} $ (here $ G $ is the Lie group for $ {\mathfrak g} $). Moreover, $ {\mathfrak g}^{Z} $ is $ 2 $-regular.
\end{corollary}

\begin{proof} The only thing to prove is that $ {\mathfrak g}^{Z} $ is $ 2 $-regular. Constant
mappings give an inclusion $ {\mathfrak g}\hookrightarrow{\mathfrak g}^{Z} $, and regular elements of $ {\mathfrak g} $ go to regular
elements of $ {\mathfrak g}^{Z} $. This shows that $ {\mathbb P}{\mathfrak g} $ intersects irregular elements of $ {\mathbb P}{\mathfrak g}^{Z} $
over a subset of codimension 2 or more, thus $ {\mathfrak g}^{Z} $ is $ 2 $-regular. \end{proof}

Note that this statement does not allow us to show that $ {\mathfrak g}^{Z} $ satisfies
conditions of Corollary~\ref{cor95.30}, since $ \left({\mathfrak g}^{Z}\right)^{*} $ is not necessarily related
to $ \left({\mathfrak g}^{*}\right)^{Z} $. However, if $ A^{*} $ is isomorphic to $ A $ as an $ A $-module, then the
action of $ G^{Z} $ on $ \left({\mathfrak g}^{Z}\right)^{*} $ is isomorphic to the action of $ G^{Z} $ on $ \left({\mathfrak g}^{*}\right)^{Z} $.

If $ Z $ is infinitesimally small, then the sketches of arguments
outlined above can be made precise. This leads to the following:

\begin{theorem} If a Lie algebra $ {\mathfrak g} $ with an anti-involution $ \iota $ satisfies
conditions of Corollary~\ref{cor95.30}, and $ A $ is a finite-dimensional
commutative algebra such that $ A $-module $ A $ is self-dual, and $ Z=\operatorname{Spec} $ A,
then $ {\mathfrak g}^{Z} $ satisfies conditions of Corollary~\ref{cor95.30}. \end{theorem}

\begin{example} $ A_{a_{1}\dots a_{k}}={\mathbb K}\left[z_{1},\dots ,z_{k}\right]/z_{1}^{a_{1}+1}\dots z_{k}^{a_{k}+1} $ is self-dual. Let
$ Z_{a_{1}\dots a_{k}}=\operatorname{Spec} A_{a_{1}\dots a_{k}} $. Then $ Z_{a_{1}\dots a_{k}}=Z_{a_{1}}\times\dots \times Z_{a_{k}} $, thus
$ {\mathfrak g}^{Z_{a_{1}\dots a_{k}}}=\left(\left({\mathfrak g}^{Z_{a_{1}}}\right)^{\dots }\right)^{Z_{a_{k}}} $. In other words, the construction above with
$ Z=Z_{a_{1}\dots a_{k}} $ is equivalent to a repeated $ k $ times construction with $ Z=Z_{a_{l}} $,
$ l=1,\dots ,k $. \end{example}

\begin{example} For a different example of the self-dual case consider
\begin{equation}
B_{k}={\mathbb K}\left[z_{1},z_{2}\right]/\left(z_{1}z_{2},z_{1}^{k}-z_{2}^{k}\right).
\notag\end{equation}
There is a relation between $ B_{k} $ and $ A_{k}\times A_{k}={\mathbb K}\left[z_{1}\right]/z_{1}^{k+1}\times{\mathbb K}\left[z_{2}\right]/z_{2}^{k+1} $. Let
$ C_{k}=A_{k}\times A_{k}/\left(\left(z_{1}^{k},-z_{2}^{k}\right)\right) $. Then $ B_{k} $ can be included into $ C_{k} $ by $ z_{1} \mapsto \left(z_{1},1\right) $,
$ z_{2} \mapsto \left(1,z_{2}\right) $, $ 1 \to $ (1,1).

The scheme $ \operatorname{Spec} A_{k}\times A_{k} $ coincides with $ Z_{k}\coprod Z_{k} $. Thus $ Y_{k}\buildrel{\text{def}}\over{=}\operatorname{Spec} C_{k} $ is a
subscheme of $ Z_{k}\coprod Z_{k} $, and $ X_{k}\buildrel{\text{def}}\over{=}\operatorname{Spec} B_{k} $ is a quotient of $ Y_{k} $. The usual
``picture'' of $ Z_{k} $ is an ``interval of infinitesimal length'' $ k $. It is hard to
picture $ C_{k} $: it should be a subscheme of two copies of such an interval,
but the sense of the equation $ z_{1}^{k}-z_{2}^{k} $ of $ Y_{k} $ inside $ Z_{k}\coprod Z_{k} $ indicates
visualization via gluing, not via cutting-off unnecessary parts: indeed,
the function $ z^{k} $ vanishes on a subinterval $ Z_{k-1} $ of $ Z_{k} $, thus in some sense
$ z^{k} $ is non-zero only at the ``end'' of $ Z_{k} $. Thus in some sense the functions
on $ Y_{k} $ ``coincide'' with functions on $ Z_{k}\coprod Z_{k} $ which have the same values at
the ``ends'' of two copies of $ Z_{k} $. On the other hand, functions on $ X_{k} $ are
{\em exactly\/} the functions on $ Y_{k} $ values of which at two centers of two copies
of $ Z_{k} $ coincide.

Thus $ X_{k} $ may be visualized as two infinitesimal intervals $ Z_{k} $ glued
both at the centers and at the ``end'' (here gluing of the centers is done
in a precise algebro-geometric sense, gluing of ends---and ends
themselves---exist only as a figure of speech). This suggests an analogy
of $ X_{k} $ with a kind of ``infinitesimal loop''. As the results above show, $ {\mathfrak g}^{X_{k}} $
has ``nice'' properties if $ {\mathfrak g} $ does. It would be interesting to compare
properties of $ {\mathfrak g}^{X_{k}} $ with properties of the Lie algebra $ \widehat{{\mathfrak g}} $ of loops in $ {\mathfrak g} $. \end{example}

\begin{example} \label{ex95.75}\myLabel{ex95.75}\relax  Consider $ Z=Z_{1} $. Clearly, $ {\mathfrak g}^{Z_{1}}={\mathfrak g} \ltimes \operatorname{ad}_{{\mathfrak g}} $; here $ \operatorname{ad}_{{\mathfrak g}} $ is
the commutative Lie algebra which coincides with $ {\mathfrak g} $ as a vector space, and
with adjoint action of $ {\mathfrak g} $. If $ {\mathfrak g} $ is reductive, then $ {\mathfrak g}^{Z_{1}}\simeq{\mathfrak g} \ltimes \operatorname{ad}_{{\mathfrak g}}^{*} $ as
well. This shows that in the situation of Example~\ref{ex11.22} the
bihamiltonian structure is flat on an open subset if $ c_{1} $ belongs to an
open subset. \end{example}

Consider now a different point of view on some of the results of
this section.

\begin{remark} It is instructive to compare the statement of Proposition
~\ref{prop11.75} with the statement of Theorem~\ref{th4.35}. Let $ {\mathfrak g}={\mathfrak g}{\mathfrak l}\left(n,{\mathbb C}\right) $, $ \iota $ be the
Cartan anti-involution. Then $ -\iota $ is an involution. The fixed points of $ -\iota $
form $ {\mathfrak o}\left(n\right) $, or linear transformations preserving a non-degenerate
symmetric bilinear form $ \alpha\left(v,w\right) $ in $ {\mathbb C}^{n} $. The anti-involutions which differ
from $ \iota $ by a conjugation will lead to equivalent bilinear forms. Note that
any two non-degenerate symmetric bilinear form in $ {\mathbb C}^{n} $ are equivalent.

Fixed points of $ \iota $ are symmetric matrices. Given such a matrix $ X $,
consider the form $ \beta\left(v,w\right)=\alpha\left(Xv,w\right) $. It is a symmetric bilinear form, any
symmetric bilinear form can be written in this way for an appropriate $ X $.
A bilinear form $ \gamma $ in $ V $ can be considered as mappings $ \widetilde{\gamma}\colon V \to V^{*}\colon \left< \widetilde{\gamma}v,w
\right>=\gamma\left(v,w\right) $. Let $ \widetilde{\alpha} $, $ \widetilde{\beta} $ be mappings associated to $ \alpha $ and $ \beta $. Then $ X=\widetilde{\alpha} \widetilde{\beta}^{-1} $. Call
$ X $ the {\em associated operator\/} for bilinear forms $ \alpha $, $ \beta $.

Now the statement of Proposition~\ref{prop11.75} in the case $ {\mathfrak g}={\mathfrak g}{\mathfrak l}\left(n,{\mathbb C}\right) $
can be read as follows: any operator $ X $ in $ n $-dimensional complex vector
space $ V $ is\footnote{Contrast this with the real case and the case of positive $ \alpha $:
then any associated operator is diagonalizable!} an associated operator of two symmetric bilinear forms $ \alpha $, $ \beta $
(here $ \alpha $ is non-degenerate). On the other hand, if one does the same for
skew-symmetric bilinear forms, then by Theorem~\ref{th4.35} the Jordan blocks
$ X $ would come in pairs. \end{remark}

\section{Appendix on Kronecker decompositions }\label{h99}\myLabel{h99}\relax 

We know that any Kronecker relation $ W $ in a vector space $ V $ can be
decomposed into a direct sum of Kronecker blocks $ W_{i} $ in subspaces $ V_{i}\subset V $.
However, this decomposition is not unique. Here we sketch the
degree of arbitrariness of this decomposition.

Given such a decomposition $ V=\bigoplus V_{i} $, consider the following objects:

\begin{definition} The {\em isotypic component\/} $ {\mathcal I}_{k}\left(V\right) $ of type $ k $ of a decomposition
$ V=\bigoplus V_{i} $ is the sum of subspaces $ V_{i} $ of dimension $ k $. The {\em isotypic filtration\/}
$ F_{k} $ of a decomposition $ V=\bigoplus V_{i} $ is $ F_{k}V=\sum_{l\leq k}{\mathcal I}_{l}\left(V\right) $.

A vector subspace $ S\subset V $ is a $ k $-{\em isotypic block\/} if there is a
decomposition $ V=\bigoplus V_{i} $ such that $ S $ is the isotypic component of type $ k $. \end{definition}

\begin{theorem} \label{th4.20}\myLabel{th4.20}\relax  The isotypic filtration of a Kronecker relation in $ V $
does not depend on the choice of a decomposition $ V=\bigoplus V_{i} $ into Kronecker
blocks. \end{theorem}

\begin{proof} Start with

\begin{definition} Given a decomposition $ V=\bigoplus V_{i} $ of a relation $ W $ in $ V $ into
Kronecker blocks and $ \lambda\in{\mathbb K}{\mathbb P}^{1} $, let $ \lambda $-{\em filtration\/} be the filtration $ F_{k}\operatorname{Ker}_{\lambda}W $ of
$ \operatorname{Ker}_{\lambda}W $ by $ \operatorname{Ker}_{\lambda}W\cap F_{k}V $. \end{definition}

Due to Lemma~\ref{lm35.30}, if $ \lambda_{i} $, $ i\geq0 $, is a sequence of different
elements of $ {\mathbb K}{\mathbb P}^{1} $, then $ F_{k}V=\sum_{i=1}^{k}F_{k}\operatorname{Ker}_{\lambda_{i}}W $. Thus it is enough to show that
the $ \lambda $-filtration in $ \operatorname{Ker}_{\lambda}W $ does not depend on the choice of decomposition
into Kronecker blocks.

Now Lemmas~\ref{lm35.30} and~\ref{lm35.40} taken together imply that
$ F_{k}\operatorname{Ker}_{\lambda_{0}}W=\operatorname{Ker}_{\lambda_{0}}W\cap\sum_{i=1}^{k}\operatorname{Ker}_{\lambda_{i}}W $. Thus $ F_{k}\operatorname{Ker}_{\lambda_{0}}W $ does not depend on the choice
of decomposition. \end{proof}

\begin{remark} Note that this theorem implies the statement of Theorem~\ref{th2.50}
about uniqueness of the collection of dimension of Kronecker blocks in
the decomposition of a given Kronecker relation. \end{remark}

Obviously, a choice of a decomposition of a given isotypic component
into Kronecker blocks is extremely non-unique if the component has more
than one block.

\begin{proposition} \label{prop4.31}\myLabel{prop4.31}\relax  Consider a Kronecker relation $ W $ in $ V $, suppose that
all the Kronecker blocks of $ W $ have the same dimension. Fix $ \lambda_{0}\in{\mathbb K}{\mathbb P}^{1} $. Given
a decomposition $ V=\bigoplus V_{i} $ into Kronecker blocks $ V_{i} $ with relations $ W_{i} $, one
obtains a decomposition $ \operatorname{Ker}_{\lambda_{0}}W=\bigoplus\operatorname{Ker}_{\lambda_{0}}W_{i} $ into a direct sum sum of
$ 1 $-dimensional subspaces.

Given an arbitrary decomposition $ \operatorname{Ker}_{\lambda_{0}}W=\bigoplus Y_{i} $ into one-dimensional
subspaces, one can find a decomposition $ V=\bigoplus V_{i} $ into Kronecker blocks such
that $ \operatorname{Ker}_{\lambda_{0}}V_{i}=Y_{i} $. The subspace $ V_{i} $ is uniquely determined by the subspace
$ Y_{i} $. \end{proposition}

\begin{proof} Suppose that dimensions of Kronecker blocks of $ V $ are equal to
$ k $. Let $ \Lambda\subset{\mathbb K}{\mathbb P}^{1} $ has $ k+1 $ element, let $ K_{\lambda}=\operatorname{Ker}_{\lambda}W $, $ \lambda\in\Lambda $. By Amplification
~\ref{amp33.20} the collection $ \left\{K_{\lambda}\right\}_{\lambda\in\Lambda} $ uniquely determines $ W $. Express possible
decompositions of $ W $ into Kronecker blocks in terms of this collection.

Suppose that $ \Lambda=\Lambda_{0}\cup\left\{\widetilde{\lambda}\right\} $, $ \widetilde{\lambda}\notin\Lambda_{0} $. By Lemmas~\ref{lm35.30},~\ref{lm35.40},
$ V=\bigoplus_{\lambda\in\Lambda_{0}}K_{\lambda} $, denote by $ \pi_{\lambda} $, $ \lambda\in\Lambda_{0} $, the projection of $ V $ on $ K_{\lambda} $ according to
this decomposition. As one can easily check, $ \pi_{\lambda}K_{\widetilde{\lambda}}\not=0 $, thus $ \pi_{\lambda}K_{\widetilde{\lambda}}=K_{\lambda} $,
provided $ \lambda\in\Lambda_{0} $ and $ V $ is a Kronecker block with $ \dim  V=k $. Hence $ \pi_{\lambda}K_{\widetilde{\lambda}}=K_{\lambda} $ for
an arbitrary $ k $-isotypic $ V $. We see that projections $ \pi_{\lambda} $ identify all the $ K_{\lambda} $
with $ K_{\widetilde{\lambda}} $, thus one with another.

Assume $ \lambda_{0}\in\Lambda_{0} $. Due to the identifications above, a choice of a basis
$ v_{i}^{\left(\lambda_{0}\right)} $ in $ K_{\lambda_{0}} $ induces bases $ v_{i}^{\left(\lambda\right)} $ in each of the subspaces $ K_{\lambda} $, $ \lambda\in\Lambda_{0} $, thus
a basis in $ V $. Let $ V_{i} $ is spanned by $ v_{i}^{\left(\lambda\right)} $, $ \lambda\in\Lambda_{0} $.

\begin{lemma} Consider a vector space $ \widetilde{V} $, $ \dim  \widetilde{V}=k $, and one-dimensional
subspaces $ T_{\lambda}\subset\widetilde{V} $, $ \lambda\in\Lambda\subset{\mathbb K}{\mathbb P}^{1} $, $ \operatorname{card}\left(\Lambda\right)=k+1 $. Suppose that each collection of $ k $
subspaces out of $ \left\{T_{\lambda}\right\} $ spans the whole vector space. Then there is one and
only one Kronecker relation $ \widetilde{W} $ in $ \widetilde{V} $ such that $ \operatorname{Ker}_{\lambda}\widetilde{W}=T_{\lambda} $, $ \lambda\in\Lambda $, \end{lemma}

\begin{proof} The ``only one'' part follows from Amplification~\ref{amp33.20}. On
the other hand, if $ W_{0} $ is a Kronecker block in $ V_{0} $, $ \dim  V_{0}=k $, then there is
(exactly one up to proportionally) linear mapping $ f $ from $ V_{0} $ to $ \widetilde{V} $ such
that $ f\left(\operatorname{Ker}_{\lambda}W\right)=T_{\lambda} $, $ \lambda\in\Lambda $. Since $ f $ is invertible, putting $ \widetilde{W}=f_{!}W_{0} $ finishes the
proof. \end{proof}

Apply the lemma to $ \widetilde{V}=V_{i} $, $ T_{\lambda}=V_{i}\cap K_{\lambda} $, $ \lambda\in\Lambda $. By the construction of the
basis $ v_{i}^{\left(\lambda\right)} $, $ T_{\lambda} $ is one-dimensional, thus the conditions of the lemma
apply. This gives a Kronecker-block linear relation $ \widetilde{W}_{i} $ in $ V_{i} $. Then $ \widetilde{W}=\bigoplus\widetilde{W}_{i} $
is a Kronecker linear relation in $ V $ with all the Kronecker blocks having
dimension $ k $, and $ \operatorname{Ker}_{\lambda}\widetilde{W}=K_{\lambda} $ for $ \lambda\in\Lambda $. By Amplification~\ref{amp33.20}, $ W=\widetilde{W} $, thus
$ \bigoplus\widetilde{W}_{i} $ is the required decomposition of $ W $ into a direct sum of Kronecker
blocks. \end{proof}

Due to Theorem~\ref{th4.20}, the only arbitrariness in the choice of
$ k $-isotypic block is the choice of an appropriate complement of $ F_{k-1}V $ in
$ F_{k}V $. Study which complements may appear as isotypic blocks.

To simplify notations we may assume that $ V=F_{k}V $ (call such $ W $ a
{\em relation of type\/} $ \leq k $). This assumption holds until the end of this
section.

\begin{definition} Given a $ k $-isotypic block $ S\subset V $ of a Kronecker relation $ W $ of
type $ \leq k $ in $ V $, let $ \lambda $-{\em pivot space\/} of $ S $ be $ S\cap\operatorname{Ker}_{\lambda}W $; here $ \lambda\in{\mathbb K}{\mathbb P}^{1} $. \end{definition}

\begin{lemma} Suppose that $ W $ is a relation of type $ \leq k $ in $ V $. Then a $ \lambda $-pivot
space is a complement to $ F_{k-1}\operatorname{Ker}_{\lambda}W $ in $ \operatorname{Ker}_{\lambda}W $. \end{lemma}

\begin{proof} This follows directly from decomposability into Kronecker
blocks. \end{proof}

\begin{definition} Consider 3 subspaces $ V,V',V'' $ of a vector space $ W $. Say that
$ V'\equiv V'' \mod V $ if $ \dim  V' =\dim  V'' $ and images of $ V' $ and $ V'' $ in $ W/V $ coincide. \end{definition}

\begin{definition} Consider a finite subset $ \Lambda\subset{\mathbb K}{\mathbb P}^{1} $ and a Kronecker relation $ W $ in
$ V $ of type $ \leq k $. A collection of vector subspaces $ S_{\lambda}\subset\operatorname{Ker}_{\lambda}W $, $ \lambda\in\Lambda $, is
{\em admissible\/} if there is a $ k $-isotypic block $ S $ such that $ S_{\lambda}=S\cap\operatorname{Ker}_{\lambda}W $ for any $ \lambda\in\Lambda $.

Call $ \left\{S_{\lambda}\right\}_{\lambda\in\Lambda} l $-{\em admissible\/} if there is a $ k $-isotypic block $ S $
such that $ S_{\lambda}\equiv S\cap\operatorname{Ker}_{\lambda}W \mod F_{l}\operatorname{Ker}_{\lambda}W $ for any $ \lambda\in\Lambda $. \end{definition}

In particular, a collection $ \left\{S_{\lambda}\right\}_{\lambda\in\Lambda} $ is $ k-1 $-admissible iff $ S_{\lambda} $ is a
complement to $ F_{k-1}\operatorname{Ker}_{\lambda}W $ in $ \operatorname{Ker}_{\lambda}W $.

\begin{definition} A sequence $ v_{1},\dots ,v_{l} $ of elements of $ V $ forms a $ W $-{\em chain\/} if for
any two consecutive elements $ v $, $ \widetilde{v} $ of the sequence $ 0,v_{1},\dots ,v_{l},0 $ the pair
$ \left(v,\widetilde{v}\right)\in W $. \end{definition}

Consider the pencil $ {\mathcal P}_{1},{\mathcal P}_{2}\colon V \to V' $ which corresponds to $ W $ as in
Section~\ref{h1}. It is clear that $ v_{1},\dots ,v_{l} $ forms a $ W $-chain iff
$ v_{1}+\lambda v_{2}+\dots +\lambda^{l-1}v_{l} $ is in the kernel of $ \lambda{\mathcal P}_{1}-{\mathcal P}_{2} $ (here we consider $ \lambda $ as a new
variable, thus the relation holds over $ {\mathbb K}\left[\lambda\right] $). Each Kronecker block $ S_{i} $ of
dimension $ l $ in $ V $ with a basis $ {\mathbit f}_{1}^{\left(i\right)},\dots ,{\mathbit f}_{l}^{\left(i\right)} $ (as in Example~\ref{ex2.30})
gives a $ W $-chain $ {\mathbit f}_{1}^{\left(i\right)},\dots ,{\mathbit f}_{l}^{\left(i\right)} $.

\begin{definition} Given two $ W $-chains $ v_{1},\dots ,v_{l} $ and $ v'_{1},\dots ,v'_{m} $, $ m<l $, define
the $ n ${\em -th elementary operation\/} as a change of $ v_{i} $ to $ v_{i}+Cv'_{i-n} $. Here $ C\in{\mathbb K} $,
$ 0\leq n\leq l-m $, and we extend the sequence $ v'_{i} $ to $ i\leq0 $ and $ i>m $ by 0. The
{\em elementary operation of the first kind\/} is the $ 0 $-th elementary operation,
{\em elementary operation of the second kind\/} is the $ l-m $-th elementary
operation. \end{definition}

Note that an elementary operation transforms a $ W $-chain into a
$ W $-chain, and that the operations of the first kind do not change $ v_{l} $,
while operations of the second kind do not change $ v_{1} $.

\begin{remark} \label{rem4.90}\myLabel{rem4.90}\relax  Consider a $ k $-isotypic block $ S $ of $ V $. Taking a $ W $-chain $ v_{i} $
corresponding to a Kronecker block of $ S $, and a $ W $-chain $ v_{i}' $ corresponding
to a Kronecker block of $ F_{k-1}V $, one can perform elementary operations
using these chains. These operations will change the chain $ v_{i} $. The
following lemma implies that this change corresponds to a change of the
$ k $-isotypic block $ S $ into another $ k $-isotypic block: \end{remark}

\begin{lemma} Suppose that vectors $ v_{ij}\in L $, $ 1\leq i\leq k $, $ j\in J $, are linearly
independent, span a complement $ S $ to $ F_{k-1}V $ in $ V $, and for any $ j\in J $ the
sequence $ v_{ij} $, $ 1\leq i\leq k $, is a $ W $-chain in $ V $. Then $ S $ is a $ k $-isotypic block in
$ V $. \end{lemma}

\begin{proof} Consider again the pencil $ {\mathcal P}_{1},{\mathcal P}_{2}\colon V \to V' $. It is clear that
$ {\mathcal P}_{1}S={\mathcal P}_{2}S $ as subspaces in $ V' $. If we prove that $ {\mathcal P}_{1}\left(F_{k-1}V\right) $ (which coincides
with $ {\mathcal P}_{2}\left(F_{k-1}V\right) $) does not intersect $ {\mathcal P}_{1}S $, then one can split $ S $ and $ {\mathcal P}_{1}S $ into
direct summands (in $ V $ and $ V' $ correspondingly), which will prove the lemma.

Suppose that $ {\mathcal P}_{1}\left(F_{k-1}V\right) $ {\em does\/} intersect $ {\mathcal P}_{1}S $. Consider an arbitrary
$ k $-isotypic block $ S' $ in $ V $. Conditions of the lemma imply that $ S\equiv S' \mod
F_{k-1}V $. This implies that $ {\mathcal P}_{1}\left(F_{k-1}V\right) $ intersects $ {\mathcal P}_{1}S' $, which is a
contradiction. \end{proof}

\begin{definition} Given a $ k $-isotypic block $ S $ of $ V $, call the modifications of $ S $
resulting from elementary operations of Remark~\ref{rem4.90} the {\em elementary
operations\/} over $ k $-isotypic blocks. \end{definition}

Obviously:

\begin{lemma} \label{lm4.110}\myLabel{lm4.110}\relax  Using elementary operations of the first kind one can
change $ S\cap\operatorname{Ker}_{0:1}W $ to become an arbitrary complement to $ F_{k-1}\operatorname{Ker}_{0:1}W $ in
$ \operatorname{Ker}_{0:1}W $ without changing $ S\cap\operatorname{Ker}_{1:0}W $. Similarly, the operations of the
second kind will do the same with $ S\cap\operatorname{Ker}_{1:0}W $. \end{lemma}

This implies

\begin{lemma} \label{lm4.120}\myLabel{lm4.120}\relax  Fix $ \lambda',\lambda''\in{\mathbb K}{\mathbb P}^{1} $, $ \lambda'\not=\lambda'' $. Suppose that $ V=F_{k}V $. Then any
$ k-1 $-admissible pair of subspaces $ \left(S_{\lambda'},S_{\lambda''}\right) $ is admissible. \end{lemma}

Moreover, one can improve this statement by considering subsets
$ \Lambda\subset{\mathbb K}{\mathbb P}^{1} $ with more than two elements. Also, one can describe the degree of
arbitrariness in the choice of an isotypic block $ S $ with given
intersections $ S\cap\operatorname{Ker}_{\lambda}W $ for $ \lambda\in\Lambda $. Start with the following

\begin{proposition} Consider two\footnote{These blocks would correspond to different decomposition of $ V $ into
direct sums of Kronecker blocks.} $ k $-isotypic blocks $ S $ and $ S' $ in $ V $. There is a
sequence of elementary operations which transforms $ S $ into $ S' $. \end{proposition}

\begin{proof} By Lemma~\ref{lm4.110}, one may suppose that $ S\cap\operatorname{Ker}_{0:1}W=S'\cap\operatorname{Ker}_{0:1}W $.
Consider $ W $-chains which form bases in $ S $ and $ S' $. Denote these chains in $ S $
by $ v_{i,j} $, in $ S' $ by $ v_{i,j}' $ (here $ j $ enumerates chains, and $ i $ vectors inside a
chain). By Proposition~\ref{prop4.31} we may assume that $ v_{1,j}=v_{1,j}' $. Let
$ p_{j}=\sum_{i}v_{ij}\lambda^{i-1} $ be the polynomial in $ \lambda $ which corresponds to the $ W $-chain
$ v_{\bullet,j} $, similarly introduce $ p_{j}' $. Fix $ j $. Obviously, $ p_{j}'-p_{j} $ can be written as
$ \lambda q $, and the polynomial $ q $ corresponds to an appropriate $ W $-chain $ \widetilde{v}_{i} $.
Moreover, all the vectors $ \widetilde{v}_{i} $ are in $ F_{k-1}V $.

Decompose $ F_{k-1}V $ into direct sum of Kronecker components, consider
projections of $ \widetilde{v}_{i} $ to these components. Clearly, for any such projection
$ \pi_{\alpha} $ the vectors $ \pi_{\alpha}\widetilde{v}_{i} $ form a $ W $-chain. Now the proposition follows from the
following

\begin{lemma} Consider a Kronecker block $ \widetilde{W} $ in $ \widetilde{V} $ and a $ \widetilde{W} $-chain $ v_{i} $ which is a
basis of $ \widetilde{V} $. Let $ v_{i}' $ be an arbitrary $ \widetilde{W} $-chain in $ \widetilde{V} $. Let $ p=\sum_{i}v_{i}\lambda^{i-1} $,
$ p'=\sum_{i}v'_{i}\lambda^{i-1} $. Then there is a polynomial $ q\in{\mathbb K}\left[\lambda\right] $ such that $ p'=qp. $ \end{lemma}

\begin{proof} Write $ v_{j}' $ in the basis $ v_{i} $ and compare the coefficients using
the definition of a $ W $-chain. \end{proof}

This finishes proof of the proposition. \end{proof}

From this proposition one can immediately deduce

\begin{theorem} Let $ S $ be a $ k $-isotypic block in $ V $, and $ \Lambda\subset{\mathbb K}{\mathbb P}^{1} $.
\begin{enumerate}
\item
Let $ S' $ be another $ k $-isotypic block in $ V $. Suppose that $ l<k $,
\begin{equation}
S\cap\operatorname{Ker}_{\lambda}W\equiv S'\cap\operatorname{Ker}_{\lambda}W \mod F_{l}\operatorname{Ker}_{\lambda}W\text{ for }\lambda\in\Lambda.
\notag\end{equation}
If $ \operatorname{card}\left(\Lambda\right)=k-l $, then
\begin{equation}
S\cap\operatorname{Ker}_{\lambda}W\equiv S'\cap\operatorname{Ker}_{\lambda}W \mod F_{l}\operatorname{Ker}_{\lambda}W\text{ for any }\lambda.
\notag\end{equation}
\item
Suppose that for $ \lambda\in\Lambda $ a vector subspace $ S'_{\lambda}\subset\operatorname{Ker}_{\lambda}W $ is fixed, and
$ S'_{\lambda}\equiv S\cap\operatorname{Ker}_{\lambda}W \mod F_{l+1}\operatorname{Ker}_{\lambda}W $ for $ \lambda\in\Lambda $. Then if $ \operatorname{card}\left(\Lambda\right)=k-l $, then there exists
another $ k $-isotypic block $ S' $ in $ V $ such that
\begin{equation}
S'_{\lambda}\equiv S'\cap\operatorname{Ker}_{\lambda}W \mod F_{l}\operatorname{Ker}_{\lambda}W.
\notag\end{equation}

\end{enumerate}
\end{theorem}

This theorem gives a complete description of the arbitrariness in
the choice of the $ k $-isotypic block in $ V $. For example, consider a subset
$ \left\{\lambda_{i}\right\} $ of $ {\mathbb K}{\mathbb P}^{1} $. Subspaces $ S\cap\operatorname{Ker}_{\lambda_{1}}W $ and $ S\cap\operatorname{Ker}_{\lambda_{2}}W $ (which may be arbitrary
complements to $ F_{k-1}\operatorname{Ker}_{\lambda_{i}}W $ in $ \operatorname{Ker}_{\lambda_{i}}W $, $ i=1,2 $) completely determine $ S \mod
F_{k-2} $. In particular, they determine $ S\cap\operatorname{Ker}_{\lambda_{3}}W \mod F_{k-2}\operatorname{Ker}_{\lambda_{3}}W $. A choice of
an arbitrary subspace $ S_{\lambda_{3}} $ of $ \operatorname{Ker}_{\lambda_{3}}W $ with the same reduction $ \mod
F_{k-2}\operatorname{Ker}_{\lambda_{3}}W $ completely determines $ S \mod F_{k-3} $ by requiring $ S\cap\operatorname{Ker}_{\lambda_{2}}W=S_{\lambda_{3}} $,
etc., etc., etc. Together with Proposition~\ref{prop4.31} this describes all
possible decompositions of $ V $ into Kronecker blocks.

\bibliography{ref,outref,mathsci}

\providecommand{\bysame}{\leavevmode\hbox to3em{\hrulefill}\thinspace}
\begin{thebibliography}{10}

\bibitem{Ada69Lec}
J.~Frank Adams, \emph{Lectures on {L}ie groups}, W. A. Benjamin, Inc., New
  York-Amsterdam, 1969.

\bibitem{Arn89Math}
V.~I. Arnol{\cprime}d, \emph{Mathematical methods of classical mechanics},
  Springer-Verlag, New York, 199?, Translated from the 1974 Russian original by
  K. Vogtmann and A. Weinstein, Corrected reprint of the second (1989) edition.

\bibitem{ArnGiv85Sym}
V.~I. Arnol{\cprime}d and A.~B. Givental{\cprime}, \emph{Symplectic geometry},
  Current problems in mathematics. Fundamental directions, Vol.\ 4, Akad. Nauk
  SSSR Vsesoyuz. Inst. Nauchn. i Tekhn. Inform., Moscow, 1985, pp.~5--139, 291.

\bibitem{Bol91Com}
A.~V. Bolsinov, \emph{Compatible {P}oisson brackets on {L}ie algebras and the
  completeness of families of functions in involution}, Izv. Akad. Nauk SSSR
  Ser. Mat. \textbf{55} (1991), no.~1, 68--92.

\bibitem{Bour75Lie78}
N.~Bourbaki, \emph{Groupes et algebres de {L}ie}, Herman, Paris, 1975, Ch7-8.

\bibitem{Dyn52Max}
E.~B. Dynkin, \emph{Maximal subgroups of the classical groups}, Trudy Moskov.
  Mat. Ob\v s\v c. \textbf{1} (1952), 39--166.

\bibitem{FadTakh87Ham}
L.~D. Faddeev and L.~A. Takhtajan, \emph{Hamiltonian methods in the theory of
  solitons}, Springer-Verlag, Berlin, 1987, Translated from the Russian by A.
  G. Reyman [A. G. Re\u\i man].

\bibitem{FokFuch80Str}
A.~S. Fokas and B.~Fuchssteiner, \emph{On the structure of symplectic operators
  and hereditary symmetries}, Lett. Nuovo Cimento (2) \textbf{28} (1980),
  no.~8, 299--303.

\bibitem{Gan59The}
F.~R. Gantmacher, \emph{The theory of matrices. {V}ols. 1, 2}, Chelsea
  Publishing Co., New York, 1959, Translated by K. A. Hirsch.

\bibitem{GelDor79Ham}
I.~M. Gel{\cprime}fand and I.~Ja. Dorfman, \emph{Hamiltonian operators and
  algebraic structures associated with them}, Funktsional. Anal. i Prilozhen.
  \textbf{13} (1979), no.~4, 13--30, 96.

\bibitem{GelZakh99Web}
Israel~M. Gelfand and Ilya Zakharevich, \emph{Webs, {L}enard schemes, and the
  local geometry of bihamiltonian {T}oda and {L}ax structures}, Archived as
  math.DG/9903080. To appear in Selecta Math.

\bibitem{GelZakhFAN}
\bysame, \emph{Spectral theory for a pair of skew-symmetrical operators on {$
  S^1 $}}, Func. Anal. Appl. \textbf{23} (1989), no.~1, 85--93.

\bibitem{GelZakhWeb}
\bysame, \emph{Webs, {V}eronese curves, and bihamiltonian systems}, J. of Func.
  Anal. \textbf{99} (1991), 150--178.

\bibitem{GelZakh93}
\bysame, \emph{On the local geometry of bihamiltonian structures}, The Gelfand
  mathematical seminar, 1990--1992 (Boston), Birkh\"auser, 1993, pp.~51--112.

\bibitem{GelZakh94Spe}
\bysame, \emph{The spectral theory for a pencil of skewsymmetrical differential
  operators of the third order}, Comm. Pure Appl. Math. \textbf{47} (1994),
  no.~8, 1031--1041.

\bibitem{GuiSte77Geo}
Victor Guillemin and Shlomo Sternberg, \emph{Geometric asymptotics}, American
  Mathematical Society, Providence, R.I., 1977, Mathematical Surveys, No. 14.

\bibitem{Kir76Loc}
A.~A. Kirillov, \emph{Local {L}ie algebras}, Uspehi Mat. Nauk \textbf{31}
  (1976), no.~4(190), 57--76.

\bibitem{KosMag96Lax}
Y.~Kosmann-Schwarzbach and F.~Magri, \emph{Lax-{N}ijenhuis operators for
  integrable systems}, J. Math. Phys. \textbf{37} (1996), no.~12, 6173--6197.

\bibitem{Kos79Sol}
Bertram Kostant, \emph{The solution to a generalized {T}oda lattice and
  representation theory}, Adv. in Math. \textbf{34} (1979), no.~3, 195--338.

\bibitem{Lax76Alm}
Peter~D. Lax, \emph{Almost periodic solutions of the {K}d{V} equation}, SIAM
  Rev. \textbf{18} (1976), no.~3, 351--375.

\bibitem{Mag78Sim}
Franco Magri, \emph{A simple model of the integrable {H}amiltonian equation},
  Journal of Mathematical Physics \textbf{19} (1978), no.~5, 1156--1162.

\bibitem{Mag88Geo}
\bysame, \emph{On the geometry of soliton equations}, Preprint, 1988.

\bibitem{Mag95Geo}
\bysame, \emph{On the geometry of soliton equations}, Acta Appl. Math.
  \textbf{41} (1995), no.~1-3, 247--270, Geometric and algebraic structures in
  differential equations.

\bibitem{McKeanPC}
Henri McKean, private communication, 1990.

\bibitem{MishFom78Eul}
A.~S. Mi{\v{s}}{\v{c}}enko and A.~T. Fomenko, \emph{Euler equation on
  finite-dimensional {L}ie groups}, Izv. Akad. Nauk SSSR Ser. Mat. \textbf{42}
  (1978), no.~2, 396--415, 471.

\bibitem{MorPiz96Eul}
Carlo Morosi and Livio Pizzocchero, \emph{On the {E}uler equation:
  bi-{H}amiltonian structure and integrals in involution}, Lett. Math. Phys.
  \textbf{37} (1996), no.~2, 117--135.

\bibitem{Pan98Sym}
Andriy Panasyuk, \emph{Symplectic realizations of bihamiltonian structures},
  preprint, 1998.

\bibitem{Pan99Ver}
\bysame, \emph{Veronese webs for bihamiltonian structures of higher corank},
  {B}anach {C}enter {P}ublications, {P}oisson geometry, {P}roceedings of the
  conference dedicated to the memory of {S}tanislaw {Z}akrzewski, Warsaw 1998
  (Warszawa) (Pawel Urba\'{n}ski and Janusz Grabowski, eds.), {I}nstytut
  {M}atematyczny {PAN}, 1999, to appear.

\bibitem{Rig95Geo}
Marie-H{\'e}l{\`e}ne Rigal, \emph{G\'eom\'etrie globale des syst\`emes
  bihamiltoniens r\'eguliers de rang maximum en dimension 5}, C. R. Acad. Sci.
  Paris S\'er. I Math. \textbf{321} (1995), no.~11, 1479--1484.

\bibitem{Rig95Tis}
\bysame, \emph{Tissus de {V}\'eron\`ese en dimension 3}, S\'eminaire Gaston
  Darboux de G\'eom\'etrie et Topologie Diff\'erentielle, 1994--1995
  (Montpellier), Univ. Montpellier II, Montpellier, 1995, pp.~iv, 63--68.

\bibitem{Rig98Sys}
\bysame, \emph{Syst\`emes bihamiltoniens en dimension impaire}, Ann. Sci.
  \'Ecole Norm. Sup. (4) \textbf{31} (1998), no.~3, 345--359.

\bibitem{Thom91Pen}
Robert~C. Thompson, \emph{Pencils of complex and real symmetric and skew
  matrices}, Linear Algebra Appl. \textbf{147} (1991), 323--371.

\bibitem{Tur89Cla}
Francisco-Javier Turiel, \emph{Classification locale d'un couple de formes
  symplectiques {P}oisson-compatibles}, Comptes Rendus des Seances de
  l'Academie des Sciences. Serie I. Mathematique \textbf{308} (1989), no.~20,
  575--578.

\bibitem{Tur99Equi}
\bysame, \emph{${C}\sp \infty$-\'equivalence entre tissus de {V}eronese et
  structures bihamiltoniennes}, C. R. Acad. Sci. Paris S\'er. I Math.
  \textbf{328} (1999), no.~10, 891--894.

\bibitem{Tur99MemA}
\bysame, \emph{M\'emoire \^a l'appui du projet de {N}ote: ${C}\sp
  \infty$-\'equivalence entre tissus de {V}eronese et structures
  bihamiltoniennes}, Preprint, 1999.

\bibitem{TurAith61Int}
H.~W. Turnbull and A.~C. Aitken, \emph{An introduction to the theory of
  canonical matrices}, Dover Publications Inc., New York, 1961.

\bibitem{Wei83Loc}
Alan Weinstein, \emph{The local structure of {P}oisson manifolds}, J.
  Differential Geom. \textbf{18} (1983), no.~3, 523--557.

\end{thebibliography}
\end{document}